\documentclass[11pt]{amsart}
\usepackage{amsmath}
\usepackage{amsfonts,amssymb}
\usepackage{enumerate}
\usepackage{amsthm}
\theoremstyle{plain}
\newtheorem{theorem}{Theorem}[section]
\newtheorem{proposition}[theorem]{Proposition}

\newtheorem{lemma}[theorem]{Lemma}
\newtheorem{corollary}[theorem]{Corollary}
\newtheorem{definition}[theorem]{Definition}

\theoremstyle{definition}
\newtheorem{nota}[theorem]{Notation}
\newtheorem{remark}[theorem]{Remark}
\newtheorem{assu}[theorem]{Assumption}
\DeclareSymbolFont{pletters}{OT1}{cmr}{m}{sl}
\DeclareMathSymbol{s}{\mathalpha}{pletters}{`s}
\DeclareMathOperator{\diff}{\mathsf{d}\!}
\DeclareMathOperator{\cn}{\mathrm{div}}
\DeclareMathOperator{\curl}{\mathrm{curl}}

\DeclareMathOperator{\HF}{HF}
\DeclareMathOperator{\LF}{LF}

\newcommand{\Er}{\mathcal{E}}
\newcommand{\Fr}{\mathcal{F}}
\newcommand{\Gr}{\mathcal{G}}
\newcommand{\Hr}{\mathcal{H}}

\newcommand{\Kr}{\mathcal{K}}
\newcommand{\Or}{\mathcal{O}}
\newcommand{\Qr}{\mathcal{Q}}
\newcommand{\Tr}{\mathcal{T}}
\newcommand{\Ur}{\mathcal{U}}

\def\Bscal#1#2{\Bigl\langle{#1}\hspace{1pt},\hspace{1pt}{#2}\Bigr\rangle}
\def\cpe{\nu}
\def\defn{\mathrel{:=}}
\def\dom{\Omega_{t}}
\def\dscal#1#2{\langle\!\langle\,{#1}\hspace{1pt},\hspace{1pt}{#2}\,\rangle\!\rangle}
\def\edtt{(\varepsilon\partial_{t})}
\def\eps{\varepsilon}
\def\ffp#1{\partial_{t}{#1}+v\cdot\nabla {#1}}
\def\ffpp#1{\bigl(\partial_{t}{#1}+v\cdot\nabla {#1}\bigr)}
\def\fpll#1#2{(\partial{#1}/\partial #2)}
\def\Fi#1#2{\Lambda^{{#2}}_{{#1}}}
\def\Fl#1#2{{#1}\rightarrow{#2}}
\def\ge{\geqslant}
\def\ges{\gtrsim}
\def\id{I}
\def\ik{m}
\def\indexg{\sigma}
\def\index{\sigma_{0}}
\def\la{\left\lvert}
\def\lA{\left\lVert}
\def\le{\leqslant}
\def\les{\lesssim}
\def\L#1{\langle{#1}\rangle}

\def\Lip{W^{1,\infty}}
\def\loc{{\text loc}}
\def\ncp{\cpe}
\def\nh#1{{H^{{#1}}_{~}}}
\def\nhz{L^{2}}
\def\nhsc#1{{{H}^{{#1}}_{\cpe}}}
\def\norm#1{\left\lVert{\smash[t]{{#1}}}\right\rVert}

\def\ns#1{{H^{{#1}}_{~}}}
\def\ndscal#1#2{\bigl\langle\!\bigl\langle{#1}\hspace{1pt},\hspace{1pt}{#2}\bigr\rangle\!\bigr\rangle}
\def\nscal#1#2{\bigl\langle{#1}\hspace{1pt},\hspace{1pt}{#2}\bigr\rangle}
\def\pa{a}
\def\para{h}
\def\pe{\kappa}
\def\PA{A}
\def\r{\mu}
\def\ra{\right\rvert}
\def\rA{\right\rVert}
\def\scal#1#2{\langle \, {#1} \hspace{1pt},\hspace{1pt} {#2} \,\rangle}

\def\tvar{\widetilde{\var}}
\def\tvari{\widetilde{p}}
\def\tvard{\widetilde{v}}
\def\tvariii{\widetilde{\theta}}
\def\tVar{\widetilde{U}}
\def\ter{\psi}
\def\TER{\Psi}
\def\tvar{\widetilde{\var}}
\def\tVar{\widetilde{U}}
\def\var{U}
\def\vari{p}
\def\variii{\theta}
\def\vard{v}
\def\xC{\mathbb{C}}
\def\xD{\mathbb{D}}
\def\xN{\mathbb{N}}

\def\xR{\mathbb{R}}
\def\xS{\mathbb{S}}
\def\xT{\mathbb{T}}
\def\xZ{\mathbb{Z}}
\def\Ze#1{Z_{\eps,\cpe}^{{#1}}\,}
\begin{document}

\pagestyle{plain}

\title{Low Mach number Limit of the Full 
       Navier--Stokes Equations}
\author{Thomas Alazard}
\address{MAB Universit\'e de Bordeaux I\\ $33405$ Talence Cedex, France}
\email{thomas.alazard@math.u-bordeaux1.fr}

\begin{abstract}
The low Mach number limit for 
classical solutions to the full Navier Stokes equations is here studied. 
The combined effects of large temperature variations 
and thermal conduction are accounted. In particular 
we consider general initial data. 
The equations lead to a singular problem, depending on a 
small scaling parameter, whose linearized is not uniformly well-posed.
Yet, it is proved 
that the solutions exist and are uniformly bounded for 
a time interval which is independent of the Mach number~${\rm Ma}\in (0,1]$, 
the Reynolds number~${\rm Re}\in [1,+\infty]$ and the P\'eclet number 
${\rm Pe}\in [1,+\infty]$. Based on uniform estimates in Sobolev spaces, 
and using a Theorem of G.~M\'etivier and S.~Schochet \cite{MS1}, we next prove that the penalized terms 
converge strongly to zero. It allows us to 
rigorously justify, at least in the whole space case, the well-known
computations given in the introduction of the P.-L. Lions' book~\cite{Lions}.
\end{abstract}

\maketitle

\vspace{-0.5cm}
\tableofcontents

\section{Introduction}

\numberwithin{equation}{section}

There are five key physical assumptions that dictate the nature of the low Mach number limit: 
the equations may be isentropic or non-isentropic; 
the fluid may be viscous or inviscid; the fluid may be an efficient or
a poor thermal conductor; the domain may be bounded or unbounded; 
the temperature variations may be small or large. 
Yet, there are only two cases where the mathematical analysis of the low Mach number limit 
is well developed: first, 
in the isentropic regime~\cite{Dan1,Dan2,DG,DGLM,Hoff,LionsMas}; second,
for inviscid and non heat-conductive fluids~\cite{TA,MS1,MS3}. 

Our goal is to start a rigorous analysis 
of the general case where 
the combined effects of large temperature variations and thermal conduction are accounted. 
As first anticipated in~\cite{Maj}, 
it yields some new problems concerning the nonlinear coupling of the equations.

\subsection{Setting the problem up}
The full Navier--Stokes equations are:
\begin{equation}\label{system:NS}
\left\{
\begin{aligned} 
&\partial_{t}\rho+\cn(\rho\vard)=0,\\
&\partial_{t}(\rho \vard)+\cn(\rho\vard\otimes\vard)+\nabla P=\cn \sigma,\\
&\partial_{t}(\rho e)+\cn(\rho\vard e)+P\cn \vard= \cn(k\nabla\mathcal{T})+
\sigma\cdot D\vard,
\end{aligned}
\right.
\end{equation}
where $\rho$, $\vard=(\vard^{1},\ldots,\vard^{d})$, $P$, $e$ and $\mathcal{T}$
denote the fluid density, velocity, pressure, energy and temperature,
respectively. We consider Newtonian gases 
with Lam\'e viscosity coefficients~$\zeta$ and~$\eta$, so that
the viscous strain tensor~$\sigma$ is given by
$$
\sigma \defn 2\zeta D \vard + \eta \cn \vard \id_{d},
$$
where % $D\vard$ is defined by 
$2D\vard = \nabla\vard+(\nabla\vard)^{t}$ and $I_{d}$ is the $d\times d$ identity matrix.

\smallbreak
%The equations~\eqref{system:NS}, encompassing as they do such a rich variety of physical processes, 
%are far too complicated to solve. 
%Therefore, 
Considerable insight comes from 
being able to simplify the description of the 
governing equations~\eqref{system:NS} by introducing clever physical models and the use of judicious 
mathematical approximations. 
To reach this aim, a standard strategy is to 
introduce dimensionless numbers which determine the relative signif{i}cance of competing physical 
processes taking place in moving fluids. 
Not only this allows to derive simplified equations of motion, but also 
reveals the central feature of the phenomenon considered.

In this paper we distinguish three dimensionless parameters:
$$
\eps\in (0,1],\qquad \r\in [0,1]\quad\mbox{and}\quad \pe\in [0,1].
$$
The first parameter~$\eps$ is the Mach number, namely the ratio of a characteristic velocity in the 
flow with the sound speed in the fluid. The parameters~$\r$ and~$\pe$ are essentially the inverses of the Reynolds and P\'eclet numbers, they 
measure the importance 
of viscosity and heat-conduction. 

To rescale the equations, there are basically two 
approaches which are available. The first is to cast 
equations in dimensionless form by scaling 
every variable by its characteristic 
value~\cite{Maj,Munz}. The second is to consider 
one of the three changes of variables:
\begin{alignat}{11}
t&\rightarrow \eps^{2}t,&&\quad x&&\rightarrow \eps x,&&\quad
\vard&&\rightarrow \eps \vard, 
&&\quad \zeta&&\rightarrow\r\zeta,&&\quad\eta&&\rightarrow\r\eta,&&\quad k&&\rightarrow\pe k,\label{scaling}\\
t&\rightarrow \eps t, &&\quad x&&\rightarrow x,
&&\quad \vard&&\rightarrow \eps \vard,&&\quad\zeta&&\rightarrow\eps\r\zeta,
&&\quad\eta&&\rightarrow \eps\r\eta,&&\quad k&&\rightarrow \eps\pe
k,\notag \\
t&\rightarrow t, &&\quad x&&\rightarrow x/\eps,
&&\quad \vard&&\rightarrow \eps \vard,&&\quad\zeta&&\rightarrow\eps\r\zeta,
&&\quad\eta&&\rightarrow \eps\r\eta,&&\quad k&&\rightarrow \eps^{2}\pe
k. \notag
\end{alignat}
See~\cite{Lions,ZeytounianE} for comments on the first two changes of
variables. The third one is related to large-amplitude high-frequency
solutions (these
rapid variations are anomalous oscillations in the context of
nonlinear geometric optics~\cite{CGM}).

Anyway, these two approaches both yield the same result. 
The full Navier--Stokes equations, written in a non dimensional way, are:
\begin{equation}\label{system:ANS}
\left\{
\begin{aligned} 
&\partial_{t}\rho+\cn(\rho\vard)=0,\\
&\partial_{t}(\rho \vard)+\cn(\rho\vard\otimes\vard)+ \frac{\nabla P}{\eps^{2}}=\r\cn\sigma,\\
&\partial_{t}(\rho e)+\cn(\rho\vard e)+P\cn
\vard=\pe\cn(k\nabla\mathcal{T})+\eps \Qr,
\end{aligned}
\right.
\end{equation}
where $\Qr\defn \eps^{\alpha}\r\,\sigma\cdot D\vard$
with $\alpha \ge 0$ (namely, $\alpha=1$ for the first two changes of
variables and $\alpha=0$ for the last one). 
%From now, we deliberately omit the
%term~$\eps^{\alpha}\r\,\sigma\cdot D\vard$ to ease the readability.

Our study is concerned with the analysis of the low Mach number limit for 
classical solutions to the full Navier Stokes 
equations~\eqref{system:ANS} in the non-isentropic
general case and for general initial data. In particular, 
the combined effects of large temperature variations 
and thermal conduction are accounted. 
We are interested in the limit~$\eps\rightarrow 0$ and 
in proving results that are independent of $\r$
and~$\pe$. The analysis contains two parts. We first prove an existence 
and uniform boundedness result for a time interval 
independent of the parameters~$\eps$,~$\r$ and~$\pe$. 
%Since we consider large temperature variations, 
%the problem is linearly (uniformly) unstable. 
%Therefore, the previous result is {\em not} a 
%consequence of known results 
%for coupled hyperbolic/parabolic nonlinear systems. 
We next study the behavior of the solutions when $\eps$ tends to $0$.

\smallbreak
Many results have been obtained in the past two decades about the justification 
of the {\em incompressible limit\/}---which is a special case of the 
{\em low Mach number approximation}. 
Concerning the Euler equations ($\r=\pe=0$), 
the study began in earlies eighties with works of Klainerman and Majda~\cite{KM1,KM2}, 
Schochet~\cite{SchoE}, Isozaki~\cite{Iso,Iso2}, Ukai~\cite{Ukai}, and others. 
As regards the isentropic Navier-Stokes equations ($\r=1$, $\zeta$ and
$\eta$ constants, $\rho=\rho(P)$), 
the mathematical analysis of the low Mach number limit has come of age since the pioneering works. 
Recent progress are presented in Danchin~\cite{Dan1,Dan2}, 
Desjardins and Grenier~\cite{DG}, 
Desjardins, Grenier, Lions and Masmoudi~\cite{DGLM}, Hoff~\cite{Hoff} and 
Lions and Masmoudi~\cite{LionsMas}. 
They are also two very interesting 
earlier results concerning the group method: Grenier~\cite{Grenier}
and Schochet~\cite{SchoFast}. 

Concerning the non-isentropic Euler
equations having general initial data, 
M\'etivier and Schochet have recently proved a couple of theorems~\cite{MS1,MS2,MS3} 
that supersede a number of earlier results (a part of their study 
is extended in \cite{TA} to the boundary case). In particular they
have proved 
the existence of classical solutions on a time
interval independent of $\eps$. The aim of this paper is precisely to 
start a rigorous analysis of the corresponding problems 
for the full Navier--Stokes equations.

The study of the low Mach number limit is a vast subject of which we have barely scratched the 
surface here. To fill in this gap we recommend Desjardins and Lin~\cite{DL} 
and Gallagher~\cite{Gallagher} for well written survey papers.
For the reader who wishes to learn more about the physics,
P.-L. Lions~\cite{Lions}, Majda~\cite{Maj} 
and~Zeytounian~\cite{ZeytounianE,ZeytounianNS} 
are good places to start. 
Detailed historical accounts of the 
subject can be found in~\cite{Munz}, along with 
a broad number of references for further reading. In connection to 
the stability analysis performed below, let us point out that   the research of numerical 
algorithms valid for all flow speeds is a very active
field~\cite{HeVW,Klein,SBGK,VT}.% and the references therein. 

%Let us emphasize that the  main point of our analysis consists in
%proving estimates in Sobolev spaces that are independent of $\eps,\r$
%and $\pe$. 

%%%%%%%%%%%%%%%%%%%%%%%%%%%%%%%%%%%%%%%%%%%%%%%%%%%%%%%%%%%%%%%%
%%%%%%%%%%%%%%%%%%%%%%%%%%%%%%%%%%%%%%%%%%%%%%%%%%%%%%%%%%%%%%%%
%%%%%%%%%%%%%%%%%%%%%%%%%%%%%%%%%%%%%%%%%%%%%%%%%%%%%%%%%%%%%%%%
%%%%%%%%%%%%%%%%%%%%%%%%%%%%%%%%%%%%%%%%%%%%%%%%%%%%%%%%%%%%%%%%
%%%%%%%%%%%%%%%%%%%%%%%%%%%%%%%%%%%%%%%%%%%%%%%%%%%%%%%%%%%%%%%%
\subsection{Uniform stability}
We consider classical solutions, that is 
solutions valued in the Sobolev spaces $H^{s}(\xD)$ with 
$s$ large enough, the domain $\xD$ being either the whole space~$\xR^{d}$ or the torus~$\xT^{d}$. 
Our main result asserts that, for perfect gases, the classical
solutions exist and are uniformly bounded for a time interval independent of $\eps$, $\r$ and~$\pe$. 
We mention 
that the case of general gases involved additional difficulties (see
Remarks~\ref{rema:GSL1} and~\ref{rema:GSL2}) 
and will be addressed in a separate paper.

We choose to work with the
unknowns~$P$,~$\vard=(\vard^{1},\ldots,\vard^{d})$ and~$\mathcal{T}$. 
In order to be closed, the system must to be augmented with 
two equations of state, prescribing the density~$\rho$ and the energy~$e$
as given functions of~$P$ and ~$\mathcal{T}$. Here, we restrict ourselves to
perfect gases so that
$$
P=R\rho\mathcal{T}\quad\mbox{and}\quad e=C_{V}\mathcal{T},
$$	 
where $R$ and $C_{V}$ are given positive constants. We begin by rewriting 
equations~\eqref{system:ANS} in terms of
$(P,\vard,\mathcal{T})$. Set $\gamma=1+R/C_{V}$. Performing linear algebra, it
is found that
\begin{equation}\label{system:ANSF}
\left\{
\begin{aligned} 
&\partial_{t}P+\vard\cdot\nabla P +\gamma P\cn\vard
=(\gamma-1)\pe\cn(k\nabla\mathcal{T})+(\gamma-1)\eps\Qr,\\
&\rho(\partial_{t}\vard+\vard\cdot\nabla\vard)+\frac{\nabla P}{\eps^{2}}=\r\cn\sigma,\\
&\rho C_{V}(\partial_{t}\mathcal{T}+\vard\cdot\nabla\mathcal{T})+P\cn\vard=\pe\cn(k\nabla\mathcal{T})+\eps\Qr.
\end{aligned}
\right.
\end{equation}
where $\rho=P/(R\Tr)$. The equations~\eqref{system:ANSF} are supplemented with initial data:
\begin{equation}
P_{\arrowvert t=0} = P_{0}, \quad \vard_{\arrowvert t=0}=\vard_{0}  
\quad\mbox{and}\quad \mathcal{T}_{\arrowvert t=0} = \mathcal{T}_{0} .\label{system:CI}
\end{equation}

Finally, it is assumed that~$\zeta$,~$\eta$ and the coefficient of thermal
conductivity $k$ are $C^{\infty}$ functions of the temperature~$\mathcal{T}$,
satisfying
$$
k(\mathcal{T})>0,\quad \zeta(\mathcal{T})>0
\quad\mbox{and}\quad
\eta(\mathcal{T})+\zeta(\mathcal{T})>0.
$$

\begin{nota}\label{nota:A}
Hereafter,~$A$ denotes the set of adimensioned parameters:
$$A \defn \bigl\{ \,\pa=(\eps,\r,\pe)\,\arrowvert\, \eps \in (0,1],~\r\in [0,1],~\pe\in [0,1]\,\bigr\}.$$
\end{nota}
\begin{theorem}\label{mainresult:PG}
Let $d\ge 1$ and $\xD$ denotes 
either the whole space $\xR^{d}$ or the 
torus $\xT^{d}$. Consider an integer $s>1+d/2$. For all positive 
$\underline{P}$,~$\underline{\mathcal{T}}$ and~$M_{0}$, there is a positive time~$T$ 
such that for all~$\pa=(\eps,\r,\pe)\in\PA$ and all initial 
data~$(P_0^{\pa},\vard_{0}^{\pa},\mathcal{T}_{0}^{\pa})$ such that~$P_{0}^{\pa}$ 
and~$\mathcal{T}_{0}^{\pa}$ take positive values and such that
$$
\eps^{-1}\lA P_{0}^{\pa} -\underline{P}\rA_{\ns{s+1}(\xD)} + 
\lA \vard_{0}^{\pa}\rA_{\ns{s+1}(\xD)} + 
\lA\mathcal{T}_{0}^{\pa}-\underline{\mathcal{T}}\rA_{\ns{s+1}(\xD)} \le M_{0},
$$
the Cauchy problem for~\eqref{system:ANSF}--\eqref{system:CI} 
has a unique solution $(P^{\pa},\vard^{\pa},\mathcal{T}^{\pa})$ such that
$(P^{\pa}-\underline{P},\vard^{\pa},\mathcal{T}^{\pa}-\underline{\mathcal{T}})
\in C^{0}([0,T];H^{s+1}(\xD))$ and such that 
$P^{\pa}$ and~$\mathcal{T}^{\pa}$ take positive values. In addition there exists a positive 
$M$, depending only on~$M_{0}$, $\underline{P}$ and $\underline{\mathcal{T}}$, such that
$$
\sup_{\pa\in\PA} \sup_{t\in [0,T]} \Bigl\{ 
\eps^{-1}\lA P^{\pa}(t)-\underline{P}\rA_{\ns{s}(\xD)} + 
\lA\vard^{\pa}(t)\rA_{\ns{s}(\xD)} + 
\lA \mathcal{T}^{\pa}(t)-\underline{\mathcal{T}}\rA_{\ns{s}(\xD)}\Bigr\} \le M.
$$
\end{theorem}

The main obstacle is that the equations leads to a singular 
problem, depending on the small scaling parameter~$\eps$, 
whose linearized is not uniformly well-posed in Sobolev spaces. 
In other words, since we consider large temperature variations, 
the problem is linearly (uniformly) unstable (see~\cite{MS1} 
for comments on this instability). 
Therefore, we cannot obtain the nonlinear energy estimates by
differentiating the equations nor by localizing in the 
frequency space by means of Littlewood-Paley operators. 
In particular, the technical aspects are different from those present in the 
previous studies of the Cauchy problem for strong 
solutions of~\eqref{system:NS}.\footnote{The latter problem has been 
widely studied. Starting from~\cite{MN} and culminating 
in~\cite{Dan,DanE1} which investigate global strong solutions 
in spaces invariant by the scaling~\eqref{scaling} with $\r=\pe=1$,
following the approach initiated by Fujita and Kato.}

\begin{remark}
\textbf{i. }We will prove a more precise result which, in particular, 
exhibits some new smoothing effects for~$\cn\vard$ and $\nabla P$ (see Theorem~\ref{theo:uniform}).

\noindent \textbf{ii.} General initial data are here considered, and allow for large density and temperature variations. 
The hypothesis $P_{0}^{\pa}(x)-\underline{P}=\Or(\eps)$ is the natural scaling to
balance the acoustic components, see~\eqref{system:NSint}
and~\cite{Dan2,KM1,Maj,MS1}. 

\noindent \textbf{iii.} One technical reason why we are uniquely 
interested in the whole space $\xR^{d}$ 
or the Torus $\xT^{d}$ is that we will make use of the 
Fourier transform tools. A more serious obstacle is that, in the boundary case, there should 
be boundary layers to analyze~\cite{BDGL}. 
For the Euler equations (that is $\r=\pe=0$), however, 
Theorem~\ref{mainresult:PG} remains valid in the boundary case~\cite{TA}. 
\end{remark}

Before leaving this paragraph, let us say some words about the difficulties 
involved in the proof of Theorems~\ref{mainresult:PG}. 
After changes of variables (given in~\S\ref{subsection:COV}), 
we are led to study a mixed hyperbolic/parabolic system of nonlinear equations of the form:
\begin{equation}\label{system:NSint}
\left\{
\begin{aligned}
&g_{1}(\variii,\eps\vari)\bigl(\ffp{\vari}\bigr)
+\frac{1}{\eps}\cn\vard -\frac{\pe}{\eps}\chi_{1}(\eps\vari)\cn(\beta(\variii)\nabla\variii\bigr)
=\Upsilon_{1},\\
&g_{2}(\variii,\eps\vari)\bigl(\ffp{\vard}\bigr)  +\frac{1}{\eps}\nabla\vari - 
\r B_{2}(\variii,\eps\vari)\vard = 0,\\
&g_{3}(\variii,\eps\vari)\bigl(\ffp{\variii}\bigr) + \cn\vard - 
\pe\chi_{3}(\eps\vari)\cn(\beta(\variii)\nabla\variii) = \eps\Upsilon_{3},
\end{aligned}
\right.
\end{equation}
where $B_{2}$ is a second order elliptic differential operator, and 
$\Upsilon_{i}$ ($i=1,3$) are of no consequence.

One of the main differences between the Euler equations \cite{TA,MS1} 
and the full equations is the following. 
When~$\pe=0$, it is typically easy to obtain~$L^{2}$ estimates uniform in~$\eps$ 
by a simple integration 
by parts in which the large terms in~$1/\eps$ cancel out (see~\cite{KM1,MS1}). 
In sharp contrast (as observed in~\cite{Maj}), when~$\pe\neq 0$ and the 
initial temperature variations are large,  the problem is more involved 
because the penalization operator is no longer skew-symmetric. 
Several difficulties also specifically arise for proving estimates 
that are independent of~$\r$ and~$\pe$. 

Another main feature of the system~\eqref{system:NSint} is that $g_{1}$ and $g_{2}$ depend on $\variii$. 
As a consequence $\nabla g_{1}$ and $\nabla g_{2}$ are of typical size $\Or(1)$. 
The system~\eqref{system:NSint} does not enter into the classical 
framework of singular limit~\cite{Grenier,SchoFast} because of this strong coupling between the 
short time-scale and the fast time-scale. That is why we cannot derive estimates in 
Sobolev norm by standards methods using differentiation of the
equations. 
In this regard, an heuristic argument is worth
emphasis: 
For non-isentropic fluids, the pressure depends on the density as well as the temperature. 
Thus, as the Mach number goes to zero, the fact the variation in pressure converges to zero 
does not preclude large temperature fluctuations (provided that they are 
compensated by density fluctuations of opposite sign). Yet, even
small perturbations of the density might result in non-small velocity variations. 

We conclude this part with a remark concerning general gases.
\begin{remark}\label{rema:GSL1}
For general equations of state, we are led to study systems having the
form~\eqref{system:NSint} where the coefficients $\chi_{1}$ and
$\chi_{3}$ depend also on $\variii$. As a consequence, the singular
operator is {\em
  nonlinear}. Some computations indicate that {\em there are no\/} $L^{2}$
estimates uniform with respect to $\eps$. 
\end{remark}

\subsection{The low Mach number limit}
We now turn to consideration of the limit of solutions of~\eqref{system:ANS} 
in~$\xR^{d}$ as the Mach number $\eps$ goes to $0$. 
The purpose of the low Mach number approximation is to justify 
that the compression due to pressure variations can be neglected. 
This is a common assumption that is made when discussing the fluid dynamics of highly subsonic flows. 
In particular, provided the sound propagation is adiabatic, 
it is the same as saying that the flow is incompressible. 
On sharp contrast, this is no longer true if the combined 
effect of large temperature variations and heat conduction is
accounted. 
Indeed, going back
to~\eqref{system:ANSF} 
we compute that, formally, the limit system reads
\begin{equation}\label{system:limit}
\left\{
\begin{aligned}
&\gamma\underline{P}\cn \vard =(\gamma-1)\pe\cn(k\nabla\mathcal{T}),\\
&\rho \ffpp{\vard} + \nabla \pi =\r\cn\sigma,\\
&\rho C_{P} \ffpp{\Tr} =\pe\cn(k\nabla\mathcal{T}),
\end{aligned}
\right.
\end{equation}
where~$\rho=\underline{P}/(R\Tr)$,~$C_{P}=\gamma C_{V}$ and, in keeping with the notations of
Theorem~\ref{mainresult:PG}, $\underline{P}$ denotes the constant
value of the pressure at spatial infinity.

\begin{theorem}\label{resu:convergence}
Fix~$\r\in [0,1]$ and~$\pe \in [0,1]$. 
Assume that~$(P^\varepsilon, \vard^\varepsilon, \mathcal{T}^\varepsilon )$ 
satisfy~\eqref{system:ANSF} and 
$$
\sup_{\eps\in (0,1]} \sup_{t\in [0,T]} \norm{\eps^{-1}(P^\varepsilon(t)
-\underline{P})}_{\nh{s}}+\norm{\vard^\varepsilon(t)}_{\nh{s}}+\norm{\mathcal{T}^\varepsilon(t) -\underline{\mathcal{T}}}_{\nh{s}}<+\infty,
$$
for some fixed time~$T>0$, reference states
$\underline{P},\underline{\Tr}$ and index~$s$ large enough. Suppose in
addition that the initial data 
$\mathcal{T}^{\varepsilon}_{\arrowvert t=0}-\underline{\Tr}$ is compactly
supported. Then, for all $s'<s$, the pressure variations $\eps^{-1}(P^\varepsilon -\underline{P})$ converges 
strongly to~$0$ in~$L^{2}(0,T;H^{s'}_{\loc}(\xR^{d}))$. Moreover, for
all $s'<s$,  
$(\vard^{\eps},\Tr^{\eps})$ converges 
strongly in $L^{2}(0,T;H^{s'}_{\loc}(\xR^{d}))$ to a limit
$(\vard,\Tr)$ satisfying the limit system~\eqref{system:limit}.
\end{theorem}

Note that 
the convergence is not uniform in time for the oscillations on the acoustic time-scale prevent 
the convergence of the solutions on a small initial layer in time. 
The key to proving this convergence result is to prove the decay 
to zero of the local energy of the 
acoustic waves. To do so 
we will consider general systems which include~\eqref{system:NSint} as
a special case. In particular, the analysis of the general systems considered below
should apply for the study of the low Mach number combustion as
described in~\cite{Maj}. 
We mention that, in view of~\cite{TA}, it seems possible to consider the same problem 
for exterior domains (which is interesting 
for aeroacoustic~\cite{ZeytounianNS}). 
Yet, we will not address this question. The results proved 
in~\cite{BDGL,MS3} indicate that the periodic case involved important 
additional phenomena.

\begin{remark}\label{rema:GSL2}For perfect gases, the limit 
constraint is linear in the following sense: it is of the form
$\cn\vard_{e}=0$ with $\vard_{e}\defn \vard- C^{te}\pe k\nabla\Tr $. 
In sharp contrast, for general gases the constraint is
nonlinear. Indeed, it reads
~$\cn\vard=f(\underline{P},\mathcal{T})\pe\cn(k\nabla\mathcal{T})$. As a
consequence, it is not immediate that, in this case, 
the corresponding Cauchy problem
for~\eqref{system:limit} is well posed.
\end{remark}

\section{Main results}

%For the sake of notational clarity, from now we deliberately omit the term~$\eps\Qr$ which appears in the
%systems ~\eqref{system:ANS} and \eqref{system:ANSF}. 
We will see in \S\ref{subsection:COV} below that one can transform the equations~\eqref{system:ANSF} into a system of the
form%~\eqref{system:NS2}. 
%We consider systems having the form:
\begin{equation}\label{system:NS2}
\left\{
\begin{aligned}
&g_{1}(\variii,\eps\vari)\bigl(\ffp{\vari}\bigr)
+\frac{1}{\eps}\cn\vard
-\frac{\pe}{\eps}B_{1}(\variii,\eps\vari)\variii = \Upsilon_{1},\\
&g_{2}(\variii,\eps\vari)\bigl(\ffp{\vard}\bigr)
+\frac{1}{\eps}\nabla\vari - \r B_{2}(\variii,\eps\vari)\vard  = 0,\\
&g_{3}(\variii,\eps\vari)\bigl(\ffp{\variii}\bigr) + \cn\vard -  \pe
B_{3}(\variii,\eps\vari)\variii
= \eps \Upsilon_{3},
\end{aligned}
\right.
\end{equation}
where the 
unknown~$(\vari,\vard,\variii)$ is a function of the variables~
$(t,x)\in\xR\times\xD$ 
with values in~$\xR\times\xR^{d}\times\xR$. Recall that $\xD$
denotes either the whole space $\xR^{d}$ or the torus $\xT^{d}$. 
Moreover, the coefficients $g_{i}$, $i=1,2,3$, are real-valued and the $B_{i}$'s are second order differential operators given by:
\begin{align*}
B_{1}(\variii,\eps\vari)&\defn\chi_{1}(\eps\vari)\cn(\beta(\variii)\nabla\cdot),\\
B_{2}(\variii,\eps\vari)&\defn\chi_{2}(\eps\vari)\cn(2\zeta(\variii)D\cdot)+\chi_{2}(\eps\vari)\nabla(\eta(\variii)\cn\cdot),\quad
2D\defn\nabla+(\nabla\cdot)^{t},\\
B_{3}(\variii,\eps\vari)&\defn\chi_{3}(\eps\vari)\cn(\beta(\variii)\nabla\cdot).
\end{align*}
Finally, $\Upsilon_{i}\defn \chi_{i}(\variii,\eps\vari)F(\variii,\sqrt{\r}\nabla\vard)$ where
$F\in C^{\infty}$ is such that $F(0)=0$.%vanishing at the origin.
%We refer to \S\ref{subsection:COV} for transformation of 
%the equations~\eqref{system:ANSF} into a system of the
%form~\eqref{system:NS2}. 

\subsection{Structural assumptions}

\begin{assu}\label{assu:structural} 
To avoid confusion, denote by $(\vartheta,\wp)\in \xR^{2}$ the place holder of the unknown $(\variii,\eps\vari)$.
\begin{enumerate}[({A}1)]  
        \item The $g_{i}$'s, $i=1,2,3$, are $C^{\infty}$ 
positive functions of $(\vartheta,\wp)\in\xR^{2}$.  
        \item The coefficients $\beta$ and $\zeta$ are $C^{\infty}$
          positive functions of $\vartheta\in\xR$, and $\eta$ is a~$C^{\infty}$ function of~$\vartheta\in\xR$ satisfying 
          $\eta(\vartheta)+\zeta(\vartheta)>0$.
        \item The coef{f}icients $\chi_{i}$, $i=1,2,3$, are
          $C^{\infty}$ positive functions of $\wp\in\xR$. Moreover, for all $\wp\in\xR$, there holds
          $\chi_{1}(\wp)<\chi_{3}(\wp)$. 
\end{enumerate}
\end{assu}
%\begin{remark}The main assumption is the inequality
%  $\chi_{1}(\wp)<\chi_{3}(\wp)$. It plays a key role for the purpose of proving
%  $L^{2}$ estimates (see the next section and also the study of the
%  simplified system~\eqref{system:ex}). Moreover, given the
%  assumption $\beta(\vartheta)>0$, it implies that the operator $B_{1}(\variii,0)-B_{3}(\variii,0)$
%  [appearing in the last equation of the limit
%  system given below in~\eqref{system:limitNS2}] is positive. Which 
%  means nothing but the fact that the temperature evolves according to
%  the standard equation of heat diffusion!
%\end{remark}
\begin{assu} \label{assu:compatibility}
Use the notation $\diff f=\fpll{f}{\vartheta}\diff\vartheta+\fpll{f}{\wp}\diff\wp$.
There exists two $C^{\infty}$ diffeomorphisms $\xR^{2}\ni(\vartheta,\wp)\mapsto
          (S(\vartheta,\wp),\wp)\in\xR^{2}$ and $\xR^{2}\ni(\vartheta,\wp)\mapsto
          (\vartheta,\varrho(\vartheta,\wp))\in\xR^{2}$ such that
          $S(0,0)=\varrho(0,0)=0$ and
\begin{equation}\label{iden:comp}
\diff S = g_{3}\diff \vartheta - g_{1}\diff\wp \quad\mbox{and}\quad
\diff\varrho =
-\frac{\chi_{1}}{\chi_{3}}g_{3}\diff\vartheta+g_{1}\diff\wp.
\end{equation}
\end{assu}

\noindent {\em Brief discussion of the hypotheses.\/} 
\begin{enumerate}[(i)]
\item The main hypothesis in Assumption~\ref{assu:structural} is the inequality
  $\chi_{1}<\chi_{3}$. It plays a crucial role for the purpose of proving
  $L^{2}$ estimates (see the section~\ref{section:apriori} and
  especially the study of the
  simplified system~\eqref{system:ex}). Moreover, given the
  assumption $\beta(\vartheta)>0$, it ensures that the operator $B_{1}(\variii,0)-B_{3}(\variii,0)$
  [which appears in the last equation of the limit
  system given below in~\eqref{system:limitNS2}] is positive. Which 
  means nothing but the fact that the limit temperature evolves according to
  the standard equation of heat diffusion!
\item The identities given in~\eqref{iden:comp} are compatibility
  conditions between the penalization operator and the viscous
  perturbation. The reason for introducing $S$ [resp. $\varrho$] 
  will be clear in~\S\ref{subsection:incompressible} [resp. \S\ref{subsection:localexistence}].
\end{enumerate}
%From now we assume that
%Assumptions~\ref{assu:structural}--\ref{assu:compatibility} are 
%satisfied. 
%We refer to \S\ref{subsection:COV} for transformation of 
%the equations~\eqref{system:ANSF} into a system of the form~\eqref{system:NS2} 
%satisfying Assumptions~\ref{assu:structural}--\ref{assu:compatibility}.
%as well as for further comments on these assumptions.

\subsection{Uniform stability result}
%We begin by making some notational conventions. 
Given~$0\le t\le T$, a normed space~$X$
and a function~$\varphi$ defined on~$[0,T]\times \xD$ with values in~$X$, 
we denote by~$\varphi(t)$ the function~$\xD\ni x\mapsto
\varphi(t,x)\in X$. 
The usual Sobolev spaces are 
denoted~$H^{\indexg}(\xD)$. Recall that, when $\xD=\xR^{d}$, they are equipped with the norms
\begin{equation*}
\lA u \rA_\nh{\indexg}^{2} \defn (2\pi)^{-d} 
\int_{\xR^{d}} \L{\xi}^{2\indexg} |\widehat{u}(\xi)|^2\,d\xi,
\end{equation*}
where~$\widehat{u}$ is the Fourier transform of~$u$ and 
$\langle\xi\rangle \defn \bigl( 1+\la\xi\ra^2 \bigr)^{1/2}$. As
concerns the case $\xD=\xT^{d}$, replace the integrals in
$\xi\in\xR^{d}$ by sums on $k\in\xZ^{d}$.

Let us introduce a bit of notation we use continually in
the sequel.
\begin{nota}
Given $\sigma\in\xR$ and $\varrho\ge 0$, set 
$\norm{u}_{H^{\sigma}_{\varrho}}\defn\norm{u}_{\nh{\sigma-1}}+\varrho\norm{u}_{\nh{\sigma}}$.
\end{nota}

%Part of the analysis will be to circumvent 
%the fact that we have no uniform $L^{2}$ estimates in general. 
%What saves the day is that we can prove a partial decoupling of the system, 
%together with some new smoothing effects for the gradient of the pressure and the 
%divergence of the velocity field. More precisely, 

Recall that 
$\PA\defn\{\,(\eps,\r,\pe)\,\arrowvert\, 
\eps\in (0,1],\r\in [0,1],\pe\in [0,1]\,\}$. The norm that we will control is the following. 

\begin{definition}\label{defi:decoupling}
Let~$T >0$,~$\sigma\ge 0$,~$\pa=(\eps,\r,\pe)\in\PA$ and set $\cpe\defn\sqrt{\r+\pe}$. 
The space~$\Hr_{\pa}^{\sigma}(T)$ consists of these 
functions~$(\vari,\vard,\variii)$ defined on~$[0,T]\times\xD$ 
with values in~$\xR\times\xR^{d}\times\xR$ 
and such that $\norm{(\vari,\vard,\variii)}_{\Hr_{\pa}^{\sigma}(T)} <+\infty$, 
where
\begin{align*}
&\norm{(\vari,\vard,\variii)}_{\Hr_{\pa}^{\sigma}(T)}\defn
\sup_{t\in[0,T]} \,\bigl\{\, \norm{(\vari(t),\vard(t))}_{H^{\sigma+1}_{\eps\cpe}} 
+ \norm{\variii(t)}_{\nhsc{\sigma+1}} \,\bigr\} +\\
&+ \left(\int_{0}^{T}\r\norm{\nabla\vard}_{H^{\sigma+1}_{\eps\cpe}}^{2} 
+ \pe \norm{\nabla\variii}_{\nhsc{\sigma+1}}^{2}
+ \pe\norm{\cn\vard}_{\nh{\sigma}}^{2}
+(\r+\pe)\norm{\nabla\vari}_{\nh{\sigma}}^{2}\,dt\right)^{\frac{1}{2}}.%\\[1.5ex]
\end{align*}
%with
%$$
%\cpe\defn\sqrt{\r+\pe}.
%$$
\end{definition}

Similarly, the space~$\Hr_{\pa,0}^{\sigma}$ consists of these 
functions~$(\vari,\vard,\variii)$ defined on~$\xD$ 
with values in~$\xR\times\xR^{d}\times\xR$ 
and such that $\norm{(\vari,\vard,\variii)}_{\Hr_{\pa,0}^{\sigma}} <+\infty$, 
where
$$
\norm{(\vari,\vard,\variii)}_{\Hr_{\pa,0}^{\sigma}}\defn \norm{(\vari,\vard)}_{H^{\sigma+1}_{\eps\cpe}} 
+\norm{\variii}_{\nhsc{\sigma+1}} \quad\mbox{with}\quad \cpe\defn\sqrt{\r+\pe}.
$$
%with again $\cpe\defn\sqrt{\r+\pe}$.
%\end{definition}

\begin{remark}One may not expect at first to have to take these norms into account, 
but they come up on their own in~\S\ref{section:apriori}.
Moreover, not only the norms, 
but also the spaces depend on~$\pa$ (since we allow~$\r=0$ or~$\pe=0$). 
\end{remark}

\begin{nota} Given a normed space~$X$ and a non-negative~$M$, we denote by~$B(X;M)$ the 
ball of center~$0$ and radius~$M$ in~$X$.
\end{nota}

Here is our main result. In the context of
Assumptions~\ref{assu:structural}--\ref{assu:compatibility} we prove
that the solutions of~\eqref{system:NS2} exist and
are uniformly bounded for a time interval which is independent of $\pa\in\PA$.

\begin{theorem}\label{theo:uniform}
Suppose that system~\eqref{system:NS2} satisfies
Assumptions~$\ref{assu:structural}$--$\ref{assu:compatibility}$, and let~$d\ge 1$.
For all integer~$s>1+d/2$ and for all positive~$M_{0}$, there exists a positive~$T$ and a positive 
$M$ such that for all~$\pa\in\PA$ and all initial data in 
$B(\Hr_{\pa,0}^{s};M_{0})$, 
the Cauchy problem for~\eqref{system:NS2} 
has a unique classical solution in~$B(\Hr_{\pa}^{s}(T);M)$.
\end{theorem} 

The proof of this result will occupy us till
\S\ref{subsection:endproof1}. 
The crucial part is to obtain estimates in Sobolev norms that are
independent of~$\pa\in\PA$. 
Notable technical aspects include 
the proof of an energy estimate for linearized equations
[see~Theorem~\ref{theo:L2}] and the use of new tools to localize in
the frequency space [see~Propositions~\ref{prop:Product} and~\ref{prop:Friedrichs}].
With these results in hands, we begin in~$\S\ref{section:HFR}$ 
by analyzing the high frequency regime. The rest of the analysis is devoted to the proof of low frequencies estimates. 
We mention that we do not need specific estimates for medium frequencies.

\begin{remark}Up to numerous changes, a close inspection of the proof of
  Theorem~\ref{theo:uniform} indicates that, in fact: for all $M>M_{0}>0$, there
  exists $T>0$ such that for all~$\pa\in\PA$ and all initial data in 
$B(\Hr_{\pa,0}^{s};M_{0})$ 
the Cauchy problem for~\eqref{system:NS2} 
has a unique classical solution in~$B(\Hr_{\pa}^{s}(T);M)$. 
%Yet, we do not address this question.
\end{remark}

\subsection{Convergence toward the solution of the limit system}
We now turn to considering the behavior of the solutions of~\eqref{system:NS2} 
in~$\xR^{d}$ as the Mach number $\eps$ tends to zero. 
Fix $\r$ and $\pe$ and consider a family
of solutions of system~\eqref{system:NS2},
$(\vari^{\eps},\vard^{\eps},\variii^{\eps})$. It is assumed to be
bound in $C([0,T];H^{\sigma}(\xR^{d}))$ with $\sigma$ large enough and
$T>0$. Strong compactness of $\variii^{\eps}$ is clear from uniform 
bounds for $\partial_{t}\variii^{\eps}$. For the sequence
$(\vari^{\eps},\vard^{\eps})$, however, the uniform bounds 
imply only weak compactness, insufficient to insure that 
the limits satisfies the limit equations. We remedy this by 
proving that the penalized terms converge strongly to zero.

\begin{theorem}\label{theo:decay}
Suppose that the system~\eqref{system:NS2} satisfies
Assumption~$\ref{assu:structural}$. 
Fix $\r\in [0,1]$ and $\pe \in [0,1]$, and let $d\ge 1$. 
Assume that $(\vari^\varepsilon, \vard^\varepsilon, \variii^\varepsilon )$ 
satisfy~\eqref{system:NS2} and are uniformly bounded in~$\Hr_{(\eps,\r,\pe)}^{s}(T)$ 
for some fixed $T>0$ and $s>4+d/2$. 
Suppose that the initial data 
$\variii^{\varepsilon}(0)$ 
converge in $H^{s} (\xR^{d})$ to a function $\variii_{0}$ 
decaying sufficiently rapidly at infinity in the sense that 
$\L{x}^{\delta}\variii_{0}\in H^{s}(\xR^{d})$ for some given
$\delta>2$, where $\L{x}\defn (1+x^{2})^{1/2}$. 

Then, for all index $s'<s$,
$p^{\eps}\rightarrow 0$ strongly in $L^{2}(0,T;H^{s'}_{loc}(\xR^{d}))$
and $\cn\vard^{\eps}-\chi_{1}(\eps
p^{\eps})\cn(\beta(\variii^{\eps})\nabla\variii^{\eps})\rightarrow
0$ strongly in $ L^{2}(0,T;H^{s'-1}_{loc}(\xR^{d}))$.
%\begin{equation}\label{lastresult}
%\begin{split}
%&p^{\eps}\rightarrow 0 \text{  strongly in  } L^{2}(0,T;H^{s'}_{loc}(\xR^{d})),\\
%&\cn\vard^{\eps}-\chi_{1}(\eps p^{\eps})\cn(\beta(\variii^{\eps})\nabla\variii^{\eps})\bigr)\rightarrow 0
%\text{  strongly in  } L^{2}(0,T;H^{s'-1}_{loc}(\xR^{d})). 
%\end{split}
%\end{equation}
\end{theorem}

The proof is given in \S\ref{section:decay}. It is based on a Theorem of M\'etivier and Schochet \cite{MS1}, 
which we describe in Theorem~\ref{theo:MS}, 
about the decay to zero of the local energy for a class 
of wave operators with time dependent coefficients.

We have just seen that $\vari^{\eps}$ converges to $0$. The following
result states that $(\vard^{\eps},\variii^{\eps})$ converges toward
the solution of the limit system.

\begin{theorem}\label{theo:limit}
Using the same assumptions and notations 
as in Theorem~$\ref{theo:decay}$, $(\vard^\varepsilon,\variii^\varepsilon )$
converges weakly in $L^{\infty}(0,T;H^{s}(\xR^{d}))$ and strongly in 
$L^{2}(0,T;H^{s'}_{loc}(\xR^{d}))$ for all $s'<s$ to a limit
$(\vard,\variii)$ satisfying
\begin{equation}\label{system:limitNS2}
\left\{
\begin{aligned}
&\cn\vard =\pe \chi_{1}(0)\cn(\beta(\variii)\nabla\variii),\\
&g_{2}(\variii,0)\bigl(\ffp{\vard}\bigr)
+\nabla\pi - \r B_{2}(\variii,0)\vard  = 0,\\
&g_{3}(\variii,0)\bigl(\ffp{\variii}\bigr) -\pe(\chi_{3}(0)-\chi_{1}(0))\cn(\beta(\variii)\nabla\variii)=0,
\end{aligned}
\right.
\end{equation}
for some $\pi$ which can be chosen such that $\nabla\pi\in C^{0}([0,T];H^{s-1}(\xR^{d}))$.
\end{theorem}

Given Theorem~\ref{theo:decay}, the proof of
Theorem~\ref{theo:limit} follows from a close inspection of the proof
of Theorem~$1.5$ in~\cite{MS1}, and so will be omitted.

%\subsection{Lack of uniform estimates for general gases}

\subsection{Changes of variables}\label{subsection:COV}
To understand the role of the thermodynamics 
we rewrite equations~\eqref{system:ANS}
in terms of the pressure fluctuations~$\vari$, velocity~$\vard$ 
and temperature fluctuations~$\variii$; where~$\vari$ and~$\variii$ are defined by
\begin{equation}\label{defi:changeofvariables}
P=\underline{P}e^{\eps\vari}, \quad \Tr=\underline{\Tr}e^{\variii}
\quad \mbox{or}\quad \vari \defn \frac{1}{\eps}\log
\Bigl(\frac{P}{\underline{P}}\Bigr),\quad
\variii\defn \log \Bigl(\frac{\Tr}{\underline{\Tr}}\Bigr),
\end{equation}
where $\underline{P}$ and $\underline{\Tr}$ are given by the statement
of Theorem~\ref{mainresult:PG}. 

We can convert the pressure and temperature evolution equations into evolution equations for 
the fluctuations $\vari$ and $\variii$. Starting
from~\eqref{system:ANSF}, it is accomplished most readily by
logarithmic differentiation ($\partial_{t,x} P = \eps P \partial_{t,x}\vari$ and
$\partial_{t,x} \Tr = \Tr \partial_{t,x}\variii$). By so doing it is
found that $(\vari,\vard,\variii)$ satisfies
\begin{equation}\label{system:TrNS2}
\left\{
\begin{aligned}
&P\bigl(\ffp{\vari}\bigr)
+\frac{\gamma P}{\eps}\cn\vard -\frac{(\gamma-1)\pe}{\eps}\cn(k\Tr\nabla\variii) =(\gamma-1)\Qr,\\
&\frac{P}{R\Tr}\bigl(\ffp{\vard}\bigr)
+\frac{P}{\eps}\nabla\vari - \r \cn(2\zeta D\vard)+\nabla(\eta\cn\vard)  = 0,\\
&\frac{C_{V}P}{R}\bigl(\ffp{\variii}\bigr) + P\cn\vard -  \pe \cn(k\Tr\nabla\variii)=\eps\Qr.
\end{aligned}
\right.
\end{equation}

Therefore, the system~\eqref{system:NS2}
includes~\eqref{system:TrNS2} as a special case where 
%Indeed, dividing the first equation in ~\eqref{system:TrNS2} by
%$\gamma P$, and dividing the
%second and third ones by $P$, transforms~\eqref{system:TrNS2} into
%a system of the form~\eqref{system:NS2} with
\begin{equation*}
g_{1}^{*}\defn \frac{1}{\gamma}, ~~ g_{2}^{*}\defn\frac{1}{R \Tr}, ~~ g_{3}^{*}\defn
\frac{C_{V}}{R},~~ \chi_{1}^{*}\defn \frac{\gamma-1}{\gamma
  P},~~\chi_{2}^{*}=\chi_{3}^{*}\defn \frac{1}{P},
~~\beta^{*}=k(\Tr)\Tr,
\end{equation*}
where the $*$ indicates that the functions are evaluated at~$(\variii,\eps\vari)$. 
Moreover, for $i=1,3$, $\Upsilon_{i}\defn\chi_{i}^{*}F(\variii,\sqrt{\r}\nabla\vard)$ where
$F(\variii,\sqrt{\r}\nabla\vard)\defn\Qr$ is as in~\eqref{system:ANS}.

We easily verify that the Assumptions~\ref{assu:structural}--\ref{assu:compatibility} are 
satisfied in this case. Hence, Theorem~\ref{mainresult:PG} as stated in the
introduction is now a consequence of Theorem~\ref{theo:uniform} since
$P$ and $\Tr$ (given by~\eqref{defi:changeofvariables}) are obviously
positive functions and since 
$\|\cdot\|_{C^{0}([0,T];H^{s}(\xD))}\le K \|\cdot\|_{\Hr^{s}_{\pa}(T)}$ 
for some constant $K$ independent of~$\pa$. 
Similarly, Theorem~\ref{resu:convergence}
follows from Theorem~\ref{theo:decay} and Theorem~\ref{theo:limit}.

\section{Localization in the frequency space}\label{Nonlinear}
We now develop the analysis 
needed to localize in the frequency space. The first paragraph is a review 
consisting of various notations 
which serve as the requested background of what follows. 
The core of this section is $\S\S$\ref{par:weightprod}--\ref{par:Friedrichs}, 
in which we prove two technical ingredients needed to localize in the low frequency region. 
This will not be used before Section~\ref{section:BF} 
and could be omitted before the reader gets there. 

\subsection{Preliminaries} 
To fix matters, in this section we work on the whole space~$\xR^d$, 
yet all the results are valid {\em mutatis mutandis\/} in the Torus~$\xT^{d}$. 
All functions are assumed to be complex valued unless otherwise specif{i}ed. 
The notation~~$d$ always refers to the dimension, 
$\index$ always refers to a real number strictly greater than $d/2$, 
and~$\para$ stands for a small parameter.

\smallbreak
From now we use~$K$ to denote a generic constant 
whose value may change from line to line in the text. 
It always stand for a constant independent of~$\para$ 
(whenever~$\para$ can play a role). 
Within the proofs, we use the notation~$A\les B$ to say that~
$A\le K B$ for such a constant~$K$.

\smallbreak
We will often write the Sobolev spaces~$H^{\indexg}(\xR^d)$ as~$H^{\indexg}$. 
Given two normed vector spaces~$X_{1}$ and~$X_{2}$,~$\mathcal{L}(X_{1},X_{2})$ denotes the space of 
bounded operators from~$X_{1}$ to~$X_{2}$. We denote by~$\lA \cdot\rA_{\Fl{X_{1}}{X_{2}}}$ its norm and 
$\mathcal{L}(X)$ is a shorthand notation for~$\mathcal{L}(X,X)$.

\smallbreak
Recall the notation
$$
\langle\xi\rangle \defn \bigl( 1+\la\xi\ra^2 \bigr)^{1/2}.
$$
A function~$q$ belongs to the class~$S^m$ ($m$ is a given real number) 
if~$q(\xi)$ is a~$C^{\infty}$ function of~$\xi\in\xR^d$ and satisf\/ies the differential 
inequalities 
$|\partial_\xi^\alpha q (\xi) | \le Q_{\alpha} \L{\xi}^{m-|\alpha |}$, 
for all~$\xi\in\xR^{d}$ and for all multi-indices~$\alpha\in\xN^{d}$. 
Such a function is called a symbol. With the best constant~$Q_\alpha$ as semi-norms, 
$S^m$ is a Fr\'echet space.
Given a symbol~$q\in S^m$, 
the Fourier multiplier associated to~$q$ is given by the operator~$\Qr$ acting 
on tempered distribution~$u$ by 
$\widehat{\Qr u} \defn q \widehat{u}$. 
For~$m,k\in\xR$, let~$\Fi{k}{m}$ denotes 
the Fourier multiplier with symbol~$\L{k\xi}^{m}$. Put another way: 
\begin{equation}\label{defi:localizor}
\Fi{k}{m} = (\id-k^{2}\Delta)^{m/2},
\end{equation}
where~$\id$ denotes the identity operator.

\smallbreak
We now introduce the first of two families of operators which are used in the sequel to localize in the 
frequency space. Let~$\jmath$ be a~$C^{\infty}$ function of~$\xi\in\xR^d$, satisfying
\begin{equation*}
0\le \jmath \le 1, \quad \jmath(\xi) = 1 \text{ for } |\xi| \le 1, 
\quad \jmath(\xi)=0 \text{ for } |\xi| \ge 2,\quad \jmath(\xi)=\jmath\left(-\xi\right).
\end{equation*}
Set~$\jmath_{\para} (\xi) = \jmath(\para \xi)$, for~$0\le\para \le 1$ and~$\xi\in\xR^d$; so that 
$\jmath_{\para}$ is supported in the ball of radius~$2/\para$ about the origin. 
Then we def\/ine~$J_{\para}~$ as the Fourier multiplier with symbol 
$\jmath_{\para}$:
$$
J_{\para} = \jmath(\para D_x).
$$ 
Let us make a series of remarks on~$J_{\para}$. 
The operator~$J_{\para}$ is self-adjoint since~$\jmath_{\para}$ is a real-valued function. 
Using that~$\jmath_{\para}$ is even, we deduce that 
$J_{\para} u$ is real-valued for any real-valued tempered distribution~$u$. 
They are smoothing operators, 
and the family~$\{\, J_{\para} \mid 0<\para\le 1 \,\}$ is an approximate identity. 
In particular, we will often use the simple observation that 
for all~$r\ge 0$, there exists a non-negative constant~$K$ so that 
for all~$\para\in(0,1]$ and~$\indexg\in\xR$, one has
\begin{equation}\label{c2_gain}
\lA J_{\para}\rA_{\Fl{H^{\indexg}}{H^{\indexg+r}}}\le \frac{K}{\para ^r}\quad\mbox{and}\quad 
\lA I-J_{\para}\rA_{\Fl{H^{\indexg}}{H^{\indexg-r}}}\le \para ^r .
\end{equation}
One reason it is interesting to assume that~$\jmath$ has compact support is the following:
\begin{equation}\label{c2_presque_projection}
J_{\para} = J_{\para} J_{c\para}, \qquad \text{for all} \,\,\,\, 0\le c\le 2^{-1}. 
\end{equation}
%This is precisely one reason it is interesting to assume that~$\jmath$ has compact support. 

%%The first result will make the use of product estimate more convenient in subsequent arguments. 

%%Then,~\eqref{Com_Paradif_2} 
%%(resp.\~\eqref{Com_Paradif_3}) is still true, as is easily verif{i}ed 
%%by arguing separately on each entry in the matrix~$f$ (resp.\ $\Qr$). 

\subsection{A product estimate}\label{par:weightprod}
%The frequency space will be divided into two parts: 
%the first where~$|\xi|\les 1/\para$, the second where~$|\xi|\ges 1/\para$. 
As alluded to above, the Friedrichs mollifiers $J_{\para}$ are 
interesting because they are essentially projection 
operators (see~\eqref{c2_presque_projection}). 
On the opposite, it is also interesting to use a family of 
invertible smoothing operators. A good candidate is 
the family 
$$
\bigl\{\, \Fi{\para}{m}\defn (\id-\para^{2}\Delta)^{m/2}\,\arrowvert\, m\le 0,~\para\in (0,1]\,\bigr\}.
$$
To begin with, we discuss the properties of these operators 
which are related to the splitting of the frequency space into two parts: 
the first where~$|\xi|\les 1/\para$, the second where~$|\xi|\ges 1/\para$. 

\smallbreak
The operators~$\Fi{\para}{m}$ depend on the scale parameter~$\para$. 
Heuristically, when~$m\ge 0$, we expect that this~$\para$-dependence 
reduces to the following alternative
$$
J_{\para} \Fi{\para}{m} \les I \quad\mbox{and}\quad (I-J_{\para}) \Fi{\para}{m} \les \para^m |D_x|^m.
$$
As regards the case of negative powers of~$\Fi{\para}{}$, we will use that 
$$
\para^{m}\Fi{\para}{-m}\les \Fi{}{-m} \quad\mbox{and}\quad \Fi{\para}{-m}\les I.
$$
These statements are made precise by the following lemma.
\begin{lemma}
\label{c2_basses_hautes}
Let~$(m_{1},m_{2})\in\xR^{2}$ be such that~$0\le m_{1}\le m_{2}$. Then, for all 
$\para\in [0,1]$ and all~$\indexg\in\xR$, we have
\begin{equation}\label{esti:handled}
\lA \para^{m_{1}} \Fi{\para}{-m_{2}} \rA_{\Fl{H^{\indexg}}{H^{\indexg+{m_{1}}}}} \le 1.
\end{equation}

For all~$m\ge 0$ and all~$c\in (0,1]$, 
there exists a positive constant~$K$ such that 
for all~$\para\in [0,1]$, and all~$\indexg\in\xR$,
\begin{gather}
\bigl\lVert J_{c \para}  \Fi{\para}{m} \bigr\rVert_{\Fl{H^{\indexg}}{H^{\indexg}}} 
\le K,\label{c2_Lambda_at_low_frequencies}\\
\bigl\lVert (I-J_{c \para}) \Fi{\para}{m} \bigr\rVert_{\Fl{H^{\indexg}}{H^{\indexg-m}}} 
\le K \para ^m .\label{c2_Lambda_at_high_frequencies}
\end{gather}
\end{lemma}
\begin{proof}
To prove these results, we check that, on the symbol level
\begin{gather*}
0 \le \para^{m_{1}}\L{\para\xi}^{-m_{2}}\L{\xi}^{m_{1}}\le 1,\\
0\le \jmath(c\para\xi) \L{\para\xi}^{-m} \le \L{2 c^{-1}}^{m},\\
0 \le \bigl( 1-\jmath(c\para\xi) \bigr)\L{\para\xi}^{m}\L{\xi}^{-m}  \le 
\para^{m}\L{c}^2.
\end{gather*}
\renewcommand{\qedsymbol}{}

\end{proof}

We are now prepared to establish the following product estimates.
\begin{proposition}\label{prop:Product}
Let~$\index>d/2$,~$(\indexg_{1},\indexg_{2})\in\xR_+^{2}$ and~$(m_{1},m_{2})\in\xR_+^{2}$ 
be such that
\begin{equation}\label{c2_cond_coef_prod_normes_eps}
\indexg_{1} + \indexg_{2} + m_{1}+m_{2} \le 2 \index.
\end{equation}
Then, there exists K depending only on~$d,\index,\indexg_{i},m_{i}$ such that for all 
$\para\in [0,1]$ and~$u_i\in H^{\index-\indexg_{i}-m_{i}}$, 
\begin{equation*}
\lA\Fi{\para}{-m_{1}-m_{2}}(u_{1}u_{2})\rA_\nh{\index-\indexg_{1}-\indexg_{2}} 
\le K\lA\Fi{\para}{-m_{1}}u_{1}\rA_\nh{\index-\indexg_{1}}
\lA\Fi{\para}{-m_{2}}u_{2}\rA_\nh{\index-\indexg_{2}}.
\end{equation*}
This result extends to vector valued functions.
\end{proposition}
\begin{remark}
In words, this proposition just says that the smoothing effect of the operators $\Fi{\para}{-m}$ is
distributive.
\end{remark}

\begin{proof}The key point is that the operators~$\Fi{\para}{-m}$ are
  invertible. It allows us to derive the desired product estimates 
from the corresponding results in the usual setting $m_{1}=m_{2}=0$. To do so Proposition~\ref{prop:Product} 
is better formulated as follows: there exists~$K$ 
such that for all~$\para\in [0,1]$ and~$f_i\in H^{\index-\indexg_{i}}(\xR^d)$,
\begin{equation}\label{c2_produit}
\bigl\lVert \Fi{\para}{-m_{1}-m_{2}}
\left\{\bigl(\Fi{\para}{m_{1}} f_{1}\bigr)\bigl( \Fi{\para}{m_{2}} f_{2}\bigr) 
\right\}\bigr\rVert_{\nh{\index-\indexg_{1}-\indexg_{2}}} 
\le  K \lA f_{1} \rA_\nh{\index-\indexg_{1}}\lA f_{2}\rA_\nh{\index-\indexg_{2}}.
\end{equation}

The proof of this claim is based on the decomposition of each 
function~$f_i$ into two pieces: its low wave number part~$J_{\para}f_{i}$, and 
its high wave number part~$(I-J_{\para})f_{i}$. 
This leads to four products that are handled the same way introducing, 
for~$c \in\{0,1\}$ and~$m\in\xR$, the Fourier multipliers
$$
\Theta_{c,m }^{\para} \defn  \para^{- c m}\Fi{\para}{m} \left( J_{\para} - c I \right).
$$
That is, as is easily verified from the triangle inequality, 
the left-hand side of~\eqref{c2_produit} is less than
$$
\sum_{ 0\le c_{1},c_{2}\le 1} 
\Bigl\lVert \para^{ c _{1} m_{1}+ c _{2} m_{2}}\Fi{\para}{-m_{1}-m_{2}}
\bigl\{ \bigl( \Theta_{c _{1},m_{1}}^{\para} f_{1} \bigr)  
\bigl(\Theta_{c _{2},m_{2}}^{\para} f_{2} \bigr) \bigr\}\Bigr\rVert_{\nh{\index-\indexg_{1}-\indexg_{2}}}.
$$
Hence, to prove the claim~\eqref{c2_produit} it suffices now to
combine three ingredients:
\begin{align*}
&\lA \para^{ c _{1} m_{1}+ c _{2} m_{2}}\Fi{\para}{-m_{1}-m_{2}} v  \rA_\nh{\index-\indexg_{1}-\indexg_{2}}
\le\lA v\rA_\nh{\index-\indexg_{1}-\indexg_{2}- c_{1} m_{1} - c _{2} m_{2}} \,,\\
&\lA v_{1} v_{2} \rA_\nh{\index-\indexg_{1}-\indexg_{2}- c _{1} m_{1}- c _{2} m_{2}} \les
\lA v_{1}\rA_\nh{\index-\indexg_{1}- c _{1} m_{1}} \lA v_{2}\rA_\nh{\index-\indexg_{2}- c _{2} m_{2}}\, ,\\
&\lA \Theta_{c_{i},m_{i} }^{\para} v \rA_\nh{\index-\indexg_{i}-c_{i} m_{i}}\les
\lA v \rA_\nh{\index-\indexg_{i}} \quad (i\in\{1,2\}) \,.
\end{align*}
The first and last inequalities follow from 
Lemma~\ref{c2_basses_hautes}. In order to prove the second one, we first recall 
the classical rule of product in Sobolev spaces (see Theorem~$8.3.1.$ in~\cite{HorL}). 
For~$(r_{1} , r_{2}) \in\xR^2$, the product maps continuously 
$H^{r_{1}}(\xR^d)\times H^{r_{2}}(\xR^d)$ to~$H^{r}(\xR^d)$ whenever
\begin{equation}\label{product:HorL} 
     r_{1}+r_{2}\ge0, \quad r\le \min \{r_{1},r_{2}\} \quad\mbox{and}\quad r\le r_{1}+r_{2}-d/2,
\end{equation}
with the third inequality strict if~$r_{1}$ or~$r_{2}$ or~$-r$ is equal to~$d/2$. 
We next verify that~\eqref{product:HorL} applies with 
$$
r\defn \index-\indexg_{1}-\indexg_{2}- c _{1} m_{1}- c _{2} m_{2},\qquad r_{i}\defn \index-\indexg_{i}-c_{i}m_{i} \quad\mbox{for}~i=1,2. 
$$
\end{proof}

\subsection{A Friedrichs' Lemma}\label{par:Friedrichs}
In this paragraph we present a result which complements the standard Friedrichs' Lemma. 
To do that we first need a commutator estimate, 
stating that the commutator of a Fourier multiplier of order~$m$ 
and the multiplication by a function is an operator of order~$m-1$. 

\begin{lemma}\label{commutateur_S^1_f}
Let~$\index>d/2+1$ and~$m\in [0,+\infty)$. For any bounded subset~$\mathcal{B}$ of~$S^m$ and all 
$\indexg\in(-\index+m,\index-1]$, 
there exists a constant~$K$ 
such that for all symbol~$q\in\mathcal{B}$, all $f\in H^{\index}(\xR^d)$, and all $v\in H^{\indexg}(\xR^{d})$,
\begin{equation}\label{commutator:estimateusual}
\lA\Qr(fu)-f\Qr u\rA_{H^{\indexg-m+1}}\le K\lA f\rA_{\nh{\index}}\lA u\rA_{H^{\indexg}},
\end{equation}
where $\Qr$ is the Fourier multiplier with symbol~$q$. 
\end{lemma}
Lemma~\ref{commutateur_S^1_f} is classical. Yet, for the convenience
of the reader, we include a proof at the end of this paragraph. 

A word of caution: the estimate~\eqref{commutator:estimateusual}
carries over to matrix valued functions and symbols except for one key point. 
Suppose that~$f$ and~$q$ are matrix valued. 
In order that~\eqref{commutator:estimateusual} be true the following condition must be fulfilled: $q(\xi) f(x) = f(x) q(\xi)$ for all 
$(x,\xi)\in\xR^{2d}$.

\begin{proposition}\label{prop:Friedrichs}
Let~$\index>d/2+1$ and~$m\in [0,1]$. For all~$\indexg$ in the interval $(-\index+m,\index-1 ]$, 
there exists a constant~$K$, 
such that for all 
$\para\in (0,1]$, all~$f\in H^{\index}(\xR^d)$ and all~$u\in H^{-\index}(\xR^d)$,
\begin{equation}\label{Friedrichs}
\bigl\lVert J_{\para}(fu)-fJ_{\para}u\bigr\rVert _{\nh{\indexg-m+1}}\le\para^{m}K
\bigl\lVert f\bigr\rVert_{\nh{\index}}
\bigl\lVert\Fi{\para}{-(\index+\indexg)}u\bigr\rVert_\nh{\indexg}.
\end{equation}
\end{proposition}
\begin{remark} The thing of interest here is that the precise rate of
  convergence does not require much on the high wave number part of $u$.
\end{remark}
\begin{proof}
Given $\mu\in\xR$ and $g\in H^{\mu}(\xR^{d})$, we denote by~$g^{\flat}$ the multiplication operator 
$H^{-\mu}(\xR^{d})\ni u \mapsto g u\in \mathcal{S}'(\xR^{d})$. Then, 
Proposition~\ref{prop:Friedrichs} can be formulated concisely in the
following way. There exists a constant~$K$, depending only on~$\index$,~$m$ and~$\indexg$, 
such that for all~$\para\in (0,1]$, 
\begin{equation*}
\bigl\lVert \Fi{}{1-m}\bigl[ J_{\para} ,f^{\flat} \bigr] \Fi{\para}{\index+\indexg} 
\bigr\rVert_{\Fl{H^{\indexg}}{H^{\indexg}}} 
\le \para^{m} K \lA f\rA_\nh{\index}.
\end{equation*}
The proof of this claim makes use of the division of the frequency space into two pieces. 
The low frequencies region~$\la\xi\ra \les 1/\para$ and 
the high frequencies region~$\la\xi\ra\ges 1/\para$. We write
\begin{equation*}
\Fi{}{1-m} \bigl[ J_{\para} ,f^{\flat} \bigr]\Fi{\para}{\index+\indexg}  = 
C_{0}^{\para} + C_{\infty} ^{\para},\\ 
\end{equation*}
where
\begin{align*}
C_{0}^{\para}&\defn\Fi{}{1-m}\bigl[J_{\para},f^{\flat}\bigr]J_{\para}'\Fi{\para}{\index+\indexg},\\
C_{\infty}^{\para}&
\defn\Fi{}{1-m}\bigl[J_{\para},f^{\flat} \bigr](I-J_{\para}')\Fi{\para}{\index+\indexg},
\end{align*}
and~$J_{\para}'=J_{5^{-1}\para}$. 
The reason for the particular choice of the constant~$5^{-1}$ will be apparent in a moment.

F\/irstly, we estimate~$C_{0}^{\para}$. We rewrite~$C_{0}^{\para}$ as 
\begin{equation}\label{c2_C_1}
C_{0}^{\para}=\para^{m}\Fi{}{1-m}
\bigl[\para^{-m}(J_{\para}-I),f^{\flat}\bigr]J_{\para}'\Fi{\para}{\index+\indexg}.
\end{equation}
The key point to estimate~$C_{0}^{\para}$ is to notice 
that~$\left\{\, \para^{-m} (\jmath_{\para}-1) \mid 0<\para\le 1 \,\right\}$ 
is a bounded family in~$S^{m}$. Indeed, if~$q\in S^{m}$ and~$q(\xi)$ vanishes for small~$\xi$, then 
the symbols~$q_{\lambda}(\cdot)=\lambda^{-m}q(\lambda\cdot)$ belong uniformly to~$S^{m}$, 
for~$0<\lambda<\infty$. 
Once this is granted, Lemma~\ref{commutateur_S^1_f} implies that
the family 
$\left\{\, \Fi{}{1-m}\bigl[ \para^{-m} (J_{\para}-I) ,f^{\flat} \bigr] \mid 0<\para\le1 \,\right\}$ 
is bounded in~$\mathcal{L} (H^{\indexg})$. 
On the other hand~$\{\,J_{\para}'\Fi{\para}{\index+\indexg}\mid 0<\para\le 1\,\}$ 
is a bounded family in~$\mathcal{L}(H^{\indexg})$, see~\eqref{c2_Lambda_at_low_frequencies}. 
Consequently, in light of~\eqref{c2_C_1}, we end up with 
\begin{align*}
\bigl\lVert C_{0}^{\para}\bigr\rVert_{\Fl{H^{\indexg}}{H^{\indexg}}} 
&\le \para^{m} 
\bigl\lVert\Fi{}{1-m}\bigl[\para^{-m}(J_{\para}-I),f^{\flat}\bigr]
\bigr\rVert_{\Fl{H^{\indexg}}{H^{\indexg}}}
\lA J_{\para}'\Fi{\para}{\index+\indexg}\rA_{\Fl{H^{\indexg}}{H^{\indexg}}}\\
&\les \para^{m}\lA f\rA_\nh{\index}.
\end{align*}

Our next task is to show similar estimates for~$C_{\infty}^{\para}$. 
As regards~$C_{\infty}^{\para}$, 
the fact that the operators~$J_{\para}$ are essentially projection operators is the key to the proof. 
More precisely, we use the identity~\eqref{c2_presque_projection} written in the form 
$J_{\para}(1-J_{\para}')=0$ [recall that~$J_{\para}'=J_{5^{-1}\para}$]. It yields
\begin{equation}\label{c2_C_2_1}
C_{\infty}^{\para} = \Fi{}{1-m} J_{\para} f^{\flat} 
\left( I-J_{\para}'\right)\Fi{\para}{\index+\indexg}.
\end{equation}
It turns out that the situation is even better. Let~$v$ belongs to the Schwartz class~$\mathcal{S}$. 
The spectrum (support of the Fourier transform, hereafter denoted by~${{\rm spec}\,}$) of 
$(J_{\para} f) \bigl( (I-J_{\para} ') v\bigr)$ is contained in 
$$
{{\rm spec}\,} (J_{\para} f) + {{\rm spec}\,} (I-J_{\para} ') v  \subset 
B\left(0,2\para^{-1}\right) + B\left(0,10\para^{-1}\right)^{c}\subset B\left(0,3\para^{-1}\right)^{c},
$$
the exterior of the ball centered at~$0$ of radius~$3\para^{-1}$. 
Which results in 
$$
J_{\para} \bigl( (J_{\para}  f) ( (I-J_{\para} ') v )\bigr)=0,
$$
that is~$J_{\para}\bigl(f((I-J_{\para} ') v )\bigr)=J_{\para} 
\bigl((I-J_{\para})f)((I-J_{\para}')v)\bigr)$. 
By combining this identity with~\eqref{c2_C_2_1}, we are left with 
\begin{equation*}
C_{\infty}^{\para}  = \Fi{}{1-m} J_{\para} \bigl((I-J_{\para} )f\bigr)^{\flat} 
\left(I-J_{\para}'\right)\Fi{\para}{\index+\indexg}.
\end{equation*}

To estimate~$\lA  C_{\infty}^{\para} \rA_{\Fl{H^{\indexg}}{H^{\indexg}}}$ 
we prove a dual estimate for the operator adjoint. 
We write~$\bigl( C_{\infty}^{\para} \bigr)^{\star}$ 
as a product of two operators:
$$
\bigl( C_{\infty}^{\para} \bigr)^{\star} =  
\Bigl\{ \para^{-\index-\indexg} (I-J_{\para}' ) \Fi{\para}{\index+\indexg} \Bigr\} 
\Bigl\{\bigl( \overline{(I-J_{\para})f}\bigr)^{\flat} 
\left( \para^{\index+\indexg} J_{\para} \right)\Fi{}{1-m}\Bigr\},
$$
where~$\bar{z}$ denotes the complex conjugated of~$z$. 
The problem reduces to establishing that
\begin{equation}\label{c2_C_2_end2}
\bigl\lVert\bigl( \overline{(I-J_{\para})f}\bigr) ^\flat 
\bigl(\para^{\index+\indexg} J_{\para} \Fi{}{1-m}\bigr) 
\bigr\rVert_{\Fl{H^{-\indexg}}{H^{\index}}} \les 
\para^{m}\lA f\rA_\nh{\index}.
\end{equation}
Indeed, since the family~$\{\, \para^{-\index-\indexg} (I-J_{\para}' ) 
\Fi{\para}{{\index+\indexg}} \mid 0<\para \le 1 \,\}$ is 
bounded in~$\mathcal{L}(H^{\index};H^{-\indexg})$ 
(see~\eqref{c2_Lambda_at_high_frequencies}),
the estimate~\eqref{c2_C_2_end2} implies that 
$$
\bigl\lVert\bigl(C_{\infty}^{\para}\bigr)^{\star}  
\bigr\rVert_{\Fl{H^{-\indexg}}{H^{-\indexg}}} \les \para^{m}\lA f\rA_\nh{\index}.
$$ 
Which in turn implies 
$\lA C_{\infty}^{\para}  \rA_{\Fl{H^{\indexg}}{H^{\indexg}}} \les \para^{m}\lA f\rA_\nh{\index}$ 
and completes the proof.

We now have to prove~\eqref{c2_C_2_end2}. 
Let~$v$ belongs to the Schwartz class~$\mathcal{S}$. 
Since~$\index>d/2$, 
we can invoke the standard tame estimate for products, which leads to
\begin{equation}\label{c2_C_2_tame}
\begin{split}
&\bigl\lVert \overline{(I-J_{\para})f} \bigl( \para^{\index+\indexg} J_{\para} \Fi{}{1-m} v  \bigr) 
\bigr\rVert_{\nh{\index}}\les  \\
&\qquad\qquad \lA (I-J_{\para})f \rA_{L^{\infty}} 
\lA \left( \para^{\index+\indexg} J_{\para} \Fi{}{1-m} \right) v \rA_\nh{\index} + \\ 
&\qquad\qquad \lA (I-J_{\para})f \rA_\nh{\index} 
\lA \left( \para^{\index+\indexg} J_{\para} \Fi{}{1-m}\right) v \rA_{L^{\infty}} .
\end{split}
\end{equation} 
To estimate these fours terms, we use the embedding of~$H^{\index-1}(\xR^d)$ into 
$L^{\infty}(\xR^d)$ and the bounds given in~\eqref{c2_gain}, to obtain
\begin{align*}
\lA (I-J_{\para})f\rA_{L^{\infty}} &\les \lA (I-J_{\para})f\rA_\nh{\index-1}
\les\para \lA f\rA_\nh{\index},\\
\lA \left( \para^{\index+\indexg} J_{\para} \Fi{}{1-m}\right) v \rA_\nh{\index} &= 
\para^{m-1}\lA \left( \para^{\index+\indexg+1-m} J_{\para} \right) v \rA_\nh{\index+1-m}
\les\!\para^{m-1}\!\lA v \rA_\nh{-\indexg},\\
\lA (I-J_{\para})f \rA_\nh{\index}&\le \lA f \rA_\nh{\index},\\
\lA \left( \para^{\index+\indexg} J_{\para}\Fi{}{1-m} \right) v \rA_{L^{\infty}} & \les
\lA \left( \para^{\index+\indexg} J_{\para}\Fi{}{1-m} \right) v \rA_\nh{\index-1}\\
&=\para^{m}\lA \left( \para^{{\index+\indexg}-m} J_{\para}\right) v \rA_\nh{\index-m}\\
&\les\para^{m}\lA  v \rA_\nh{\index-m-({\index+\indexg}-m)}
=\para^{m}\lA  v \rA_\nh{-\indexg}.
\end{align*}

Inserting the four estimates we just proved in~\eqref{c2_C_2_tame}, we 
conclude that
\begin{equation*}
\bigl\lVert \overline{(I-J_{\para})f} \bigl( \para^{\index+\indexg} J_{\para}\Fi{}{1-m} v \bigr) 
\bigr\rVert_{\nh{\index}}\les
\para^{m} \lA f\rA_\nh{\index}\lA v \rA_\nh{-\indexg}.
\end{equation*}
This completes the proof of \eqref{c2_C_2_end2}.
%\qed
\end{proof}

\begin{proof}[Proof of Lemma~\ref{commutateur_S^1_f}] To avoid 
trivialities ($\indexg\in\emptyset$), assume that $m\le 2\index-1$.
We establish the estimate 
by using the para-differential calculus of Bony~\cite{Bony}. 
In keeping with the notations of the previous proof, $f^{\flat}$
denotes the multiplication operator $ u\mapsto fu$. 
We denote by~$T_f$ the operator of para-multiplication by~$f$. 
Rewrite the commutator~$\left[ \Qr, f^{\flat} \right]$ as
\begin{equation*}
\left[\Qr,T_{f}\right]+\Qr(f^{\flat}-T_{f})-(f^{\flat}-T_{f})\Qr.
\end{equation*}
The claim then follows from the bounds 
\begin{alignat}{2}
&\forall \indexg\in\xR,&
\lA\left[\Qr,T_f \right]\rA _{\Fl{H^{\indexg}}{H^{\indexg-m+1}}}&\le 
c_{1}(q,\indexg)\lA f\rA_\nh{\index},\label{Com_Paradif_1}\\
&\forall\indexg\in(-\index,\index-1],&
\bigl\lVert f^{\flat}-T_{f}\bigr\rVert_{\Fl{H^{\indexg}}{H^{\indexg+1}}}  
&\le c_{2}(\indexg) 
\lA f\rA_{\nh{\index}},\label{Com_Paradif_2}\\
&\forall\indexg\in\xR,&
\lA\Qr\rA_{\Fl{H^{\indexg}}{H^{\indexg-m}}}&\les\sup\nolimits_{\xi}\la\L{\xi}^{-m}q(\xi)\ra,
\label{Com_Paradif_3}
\end{alignat}
where~$c_{1}(\cdot,\indexg)\colon S^m\rightarrow \xR_+$ takes bounded
sets to bounded sets. We refer the reader to~\cite{Bony} and~\cite{Meyer} 
for the proofs of the first two inequalities (see also~\cite[Prop.~$10.2.2$]{HorL} and 
\cite[Th.~$9.6.4'$]{HorL} for the proof of~\eqref{Com_Paradif_1}; 
a detailed proof of~\eqref{Com_Paradif_2} is given in \cite[Prop.~$10.2.9$]{HorL}). 
The estimate~\eqref{Com_Paradif_3} is obvious.
 \end{proof}

\section{Energy estimates for the linearized system}\label{section:apriori}
Many results have been obtained in the past two decades concerning the 
symmetrization of the Navier Stokes equations (see, e.g.,~\cite{BRDL,DanE1,KS,KY}).
Yet, the previous works do not include the dimensionless numbers. 
Here we prove estimates valid for all $a=(\eps,\r,\pe)$ in
$\PA$, where~$A$ is defined in Notation~\ref{nota:A}. 
As already written in the introduction, our result improves earlier 
works \cite{TA,KM1,MS1} on allowing $\pe\neq 0$. Indeed, 
when~$\pe=0$, the penalization operator is skew-symmetric and hence
the perturbation terms do not appear in the $L^{2}$ estimate, 
so that the classical proof for solutions to the unperturbed equations
holds. In sharp contrast (as observed in~\cite{Maj}), when~$\pe\neq 0$ and the 
initial temperature variations are large,  the problem is 
more involved because the singular operator is no longer skew-symmetric. 
Several difficulties also specifically arise for the purpose of
proving estimates that are independent of~$\r$ and~$\pe$. 
In this regard we prove some additional smoothing 
effects for $\cn\vard$ and $\nabla\vari$.

\smallbreak

We consider the following linearized equations:
\begin{equation}\label{system:NSi}
\left\{
\begin{aligned}
&g_{1}(\phi)\bigl(\partial_{t}\tvari+\vard\cdot\nabla {\tvari}\bigr)
+\frac{1}{\eps}\cn \tvard-\frac{\pe}{\eps} \cn(\beta_{1}(\phi)\nabla\tvariii)=f_{1},\\
&g_{2}(\phi)\bigl(\partial_{t}\tvard+\vard\cdot\nabla {\tvard}\bigr)
+\frac{1}{\eps}\nabla\tvari-\r\beta_{2}(\phi)\Delta\tvard-\r\beta_{2}^{\sharp}(\phi)\nabla\cn\tvard=f_{2},\\
&g_{3}(\phi)\bigl(\partial_{t}\tvariii
+\vard\cdot\nabla{\tvariii}\bigr)
+\cn\tvard -\pe\beta_{3}(\phi)\Delta\tvariii=f_{3}.
\end{aligned}
\right.
\end{equation}
To fix matters, the 
unknown~$(\tvari,\tvard,\tvariii)$ is a function of the variables~
$(t,x)\in [0,T]\times\xD$ ($T$ is a given positive real number and
$\xD$ denotes either $\xR^{d}$ or $\xT^{d}$) 
with values in~$\xR\times\xR^{d}\times\xR$. 
The coefficients~$\phi=\phi(t,x)$ and~$\vard=\vard(t,x)$ 
take their values in~$\xR^{N}$ and~$\xR^{d}$, respectively ($N$ is a given integer). 

\smallbreak

Parallel to Assumption~\ref{assu:structural}, we make the
following hypotheses.
\begin{assu}\label{assu:linearized}
Throughout this section we require~$g_{1},g_{2},g_{3}$ to be
$C^{\infty}$ positive functions of~$\phi\in\xR^{N}$, 
without recalling this assumption explicitly in the
statements. Similarly, it is assumed that~$\beta_{1}$, $\beta_{2}$,
$\beta_{2}^{\sharp}$ and $\beta_{3}$ are $C^{\infty}$ functions of
$\phi\in\xR^{N}$ satisfying
$$
\beta_{1}>0,\quad \beta_{2}>0,\quad
\beta_{2}^{\sharp}+\beta_{2}>0, 
\quad \beta_{3}>0.
$$
Our main assumption reads
\begin{equation}\label{mainassumption}
\beta_{1}< \beta_{3}.
\end{equation}
\end{assu}
\begin{assu}\label{assu:apriori} For our
purposes, it is sufficient to prove {\em a priori\/} estimates. 
It is always assumed that the unknown
$\tvar\defn(\tvari,\tvard,\tvariii)$, the coefficients $(\vard,\phi)$
as well as the source term $f\defn(f_{1},f_{2},f_{3})$ are in $C^{0}([0,T];H^{\infty}(\xD))$.
\end{assu}

\subsection{Statements of the results}\label{subsection:reultsL2}
We establish estimate on~$\norm{\tvar}_{\Hr^{0}_{\pa}(T)}$ in
terms of the norm~$\norm{\tvar(0)}_{\Hr^{0}_{\pa,0}}$ of the data and
norm of the source term~$f$. Recall from
Definition~\ref{defi:decoupling} that
\begin{align*}
&\norm{\tvar}_{\Hr^{0}_{\pa}(T)}\defn
\sup_{t\in [0,T]}\,\bigl\{\,\norm{(\tvari,\tvard)(t)}_{H^{1}_{\eps\ncp}}
+\norm{\tvariii(t)}_{H^{1}_{\ncp}}\,\bigr\}\\
&\quad+\biggl(\int_{0}^{T}\!\!\pe\norm{\nabla\tvariii}_{H^{1}_{\ncp}}^{2}
+\r\norm{\nabla\tvard}_{H^{1}_{\eps\ncp}}^{2}
+\pe\norm{\cn\tvard}_{\nhz}^{2}+(\r+\pe)\norm{\nabla\tvari}_{\nhz}^{2} \, dt\biggr)^{\frac{1}{2}},
\end{align*}
with $\ncp\defn\sqrt{\r+\pe}$. 
Recall the following
notation we use continually in the sequel: 
given $\sigma\in\xR$ and $\varrho\ge 0$ we set 
$\norm{\cdot}_{H^{\sigma}_{\varrho}}\defn
\norm{\cdot}_{\nh{\sigma-1}}+\varrho\norm{\cdot}_{\nh{\sigma}}$. 
\begin{theorem}\label{theo:L2}
There is a smooth non-decreasing function~$C$ from~$[0,+\infty)$ 
to~$[0,+\infty)$ such that
for all~$\pa\in\PA$, all~$T\in(0,1]$, all coefficients~$\vard$ and $\phi$, 
and all~$(\tVar,f)$ satisfying~\eqref{system:NSi}, the
norm~$\norm{\tvar}_{\Hr^{0}_{\pa}(T)}$ satisfies the estimate
\begin{equation*}
\norm{\tvar}_{{\Hr^{0}_{\pa}(T)}}\le
C(R_{0})e^{TC(R)}\norm{\tvar(0)}_{\Hr^{0}_{\pa,0}}
+C(R) \int_{0}^{T}\norm{(f_{1},f_{2})}_{H^1_{\eps\ncp}}+\norm{f_{3}}_{H^{1}_{\ncp}}\,dt,
\end{equation*}
where $\ncp\defn\sqrt{\r+\pe}$ and
\begin{equation}\label{defi:Rcoef}
R_{0}\defn\norm{\phi(0)}_{L^{\infty}(\xD)},~~~ R\defn\sup_{t\in[0,T]}
\norm{(\phi,\partial_{t}\phi,\nabla\phi,\ncp\nabla^{2}\phi,\vard,\nabla\vard)(t)}_{L^{\infty}(\xD)}.
\end{equation}
\end{theorem}
\begin{remark}
By Assumption~\ref{assu:apriori}, all the functions we
encounter in the previous statement are in $C^{0}([0,T];H^{\infty}(\xD))$. Yet, one can
verify that the estimate is valid whenever its two sides are well defined.
\end{remark}

\smallbreak
Theorem~\ref{theo:L2} is not enough for the purpose of proving high frequency 
estimates independent of $\pe$. We will need the following version.
\begin{theorem}\label{theo:L2slow}
The statement of Theorem~$\ref{theo:L2}$ remains valid if 
the system~\eqref{system:NSi} is replaced with
\begin{equation}\label{system:NSi'}
\left\{
\begin{aligned}
&g_{1}(\phi)\bigl(\partial_{t}\tvari+\vard\cdot\nabla {\tvari}\bigr)
+\frac{1}{\eps}\cn \tvard-\frac{\pe}{\eps} \cn(\beta_{1}(\phi)\nabla\tvariii)=f_{1},\\
&g_{2}(\phi)\bigl(\partial_{t}\tvard+\vard\cdot\nabla {\tvard}\bigr)
+\frac{1}{\eps}\nabla\tvari-\r\beta_{2}(\phi)\Delta\tvard-\r\beta_{2}^{\sharp}(\phi)\nabla\cn\tvard=f_{2},\\
&g_{3}(\phi)\bigl(\partial_{t}\tvariii
+\vard\cdot\nabla{\tvariii}\bigr) +G(\phi,\nabla\phi)\cdot\tvard
+\cn\tvard -\pe\beta_{3}(\phi)\Delta\tvariii=f_{3},
\end{aligned}
\right.
\end{equation}
where $G$ is smooth in its arguments with values in~$\xR^{d}$.
\end{theorem}

The proof of Theorem~\ref{theo:L2} relies upon two $L^{2}$ estimates. 
Because they may be useful in other circumstances, we give separate statements in~$\S\ref{2L2esti}$. 
The proof of Theorem~\ref{theo:L2slow} follows from a close
inspection of the proof of Theorem~\ref{theo:L2}, and so will be
omitted.

%We make the following running convention: 
%%In all the statements, we will focus on estimates: 
%%it is always assumed that the functions are smooth enough so that 
%%the norms involved in the estimates are well-defined. 
%Given a normed space $X$ and function $u\colon [0,T]\rightarrow X$, we 
%denote by $\norm{u}_{X}$ the function 
%$[0,T]\ni t \mapsto \norm{u(t)}_{X}$. In particular, we do not make further reference to 
%the time variable $t$ (as it is customary). Yet it is always checked that 
%the estimates are independent of the parameter $T$.

\subsection{Example}\label{L2:example}
Some important features of the proof of Theorem~\ref{theo:L2} 
can be revealed by analyzing the following simplified system:
\begin{equation}\label{system:ex}
\left\{
\begin{aligned}
&\partial_{t}\vari+\frac{1}{\eps}\cn\vard-\frac{1}{\eps}\Delta\variii=0,\\
&\partial_{t}\vard+\frac{1}{\eps}\nabla\vari=0,\\
&\partial_{t}\variii+\cn\vard -\beta\Delta\variii=0.
\end{aligned}
\right.
\end{equation}
For the sake of notational simplicity we abandon the tildes in this subsection. 
Parallel to~\eqref{mainassumption}, we suppose
\begin{equation}\label{exam:mainassu}
\beta>1.
\end{equation}

To symmetrize the large terms in $\eps^{-1}$, 
we introduce $v_{e}\defn v-\nabla\theta$. This change of variables transforms~\eqref{system:ex} into
\begin{equation*}
\left\{
\begin{aligned}
&\partial_{t}\vari+\frac{1}{\eps}\cn\vard_{e}=0,\\
&\partial_{t}\vard_{e}+\frac{1}{\eps}\nabla\vari-\nabla\cn\vard_{e}+(\beta-1)\nabla\Delta\theta=0,\\
&\partial_{t}\theta+\cn\vard_{e} -(\beta-1)\Delta\theta=0.
\end{aligned}
\right.
\end{equation*}
We take the $L^{2}$ scalar product of the first [resp.\ the second] 
equation with $\vari$ [resp. $\vard_{e}$]. It yields:
\begin{equation}\label{iden:ex1}
\frac{1}{2}\frac{d}{dt}\norm{(\vari,\vard_{e})}_{\nhz}^{2}
+\norm{\cn\vard_{e}}_{\nhz}^{2}-(\beta-1)\scal{\Delta\theta}{\cn\vard_{e}}=0.
\end{equation}
We take the $L^{2}$ scalar product of the third equation 
with $-\eta\Delta\theta$, where $\eta$ is a positive constant to be determined later on. 
It yields
\begin{equation*}
\begin{split}
&\frac{1}{2}\frac{d}{dt}\norm{(\vari,\vard_{e},\sqrt{\eta}\nabla\theta)}_{\nhz}^{2}\\
&\quad +\norm{\cn\vard_{e}}_{\nhz}^{2}
-(\beta-1+\eta)\scal{\Delta\variii}{\cn\vard_{e}}+\eta(\beta-1)\norm{\Delta\theta}_{\nhz}^{2}=0 .
\end{split}
\end{equation*}
Set $\eta\defn\beta-1$, which is positive thanks to assumption~\eqref{exam:mainassu}. We get
$$
\frac{1}{2}\frac{d}{dt}\norm{(\vari,\vard_{e},\sqrt{\beta-1}\nabla\theta)}_{\nhz}^{2}+\norm{\cn\vard_{e}-(\beta-1)\Delta\variii}_{\nhz}^{2}=0.
$$
Let $t\ge 0$. Integrating the previous identity from $0$ to $t$, 
and using the triangle inequality to replace $\vard_{e}$ by
$\vard$, we get 
\begin{equation}\label{iden:ex2}
\norm{(\vari,\vard,\nabla\theta)(t)}_{\nhz}^{2}+\int_{0}^{t}\norm{\cn\vard-\beta\Delta\variii}_{\nhz}^{2}\,d\tau
\le K_{\beta}\norm{(\vari,\vard,\nabla\theta)(0)}_{\nhz}^{2}.
\end{equation}
Here and below, $K_{\beta}$ denotes a generic constant which depends
only on $\beta$. 

We thus have proved an $L^{2}$ estimate independent of $\eps$. 
To go beyond and obtain smoothing effect on $\cn\vard$ it is sufficient 
to estimate $\Delta\variii$ independently. The strategy is to
incorporate the troublesome term $\cn\vard$ [in the equation for
$\variii$] into a skew-symmetric operator.

To do so introduce 
$$
\zeta\defn \eps\beta\vari-\variii\quad\mbox{and}\quad v_{\eps}\defn \eps\vard.
$$
We compute
 \begin{equation*}
\left\{
\begin{aligned}
&\partial_{t}\zeta+\frac{\beta-1}{\eps}\cn\vard_{\eps}=0,\\
&\partial_{t}\vard_{\eps}+\frac{1}{\beta\eps}\nabla\zeta+\frac{1}{\beta\eps}\nabla\variii=0,\\
&\partial_{t}\variii+\frac{1}{\eps}\cn\vard_{\eps} -\beta\Delta\variii=0.
\end{aligned}
\right.
\end{equation*}
We now use the essential feature of the system, that is 
assumption~\eqref{exam:mainassu}. Multiply the first equation 
[resp. the second] by $1/(\beta-1)$ [resp. $\beta$], to put 
the penalization operator in symmetric form. The energy estimate thus reads
\begin{equation*}
\frac{1}{2}\frac{d}{dt}
\norm{(\sqrt{1/(\beta-1)}\zeta,\sqrt{\beta}\vard_{\eps},\variii)}_{\nhz}^{2}
+\beta\norm{\nabla\variii}_{\nhz}^{2}=0.
\end{equation*}
Integrate this inequality, to obtain
\begin{equation}\label{iden:ex3}
\norm{(\zeta,\vard_{\eps},\theta)(t)}_{\nhz}^{2}+
\int_{0}^{t}\norm{\nabla\theta}_{\nhz}^{2}\,d\tau\le 
K_{\beta}\norm{(\zeta,\vard_{\eps},\theta)(0)}_{\nhz}^{2}.
\end{equation}
Since the coefficients are constants, the same estimate holds true for the 
first order derivatives of $\zeta$, $\vard_{\eps}$ and 
$\variii$. Namely, one has
\begin{equation}\label{iden:ex4}
\norm{\nabla(\zeta,\vard_{\eps},\theta)(t)}_{\nhz}^{2}+
\int_{0}^{t}\norm{\nabla^{2}\theta}_{\nhz}^{2}\,d\tau\le K_{\beta}\norm{\nabla(\zeta,\vard_{\eps},\theta)(0)}_{\nhz}^{2}.
\end{equation}
On applying the triangle inequality, one can replaced $\zeta$ with $\eps\vari$ 
in the previous estimates. Finally, by
combining~\eqref{iden:ex2},~\eqref{iden:ex3} and~\eqref{iden:ex4}, we get
\begin{align}
&\norm{(\vari,\vard,\variii)(t)}_{\nhz}+
\norm{\nabla(\variii,\eps\vari,\eps\vard)(t)}_{\nhz}+
\left(\int_{0}^{t}\norm{\cn\vard}_{\nhz}^{2} +\norm{\nabla\theta}_{\nh{1}}^{2} \,d\tau\right)^{1/2}\notag \\
&\qquad\qquad\qquad\qquad  \le K_{\beta}\norm{(\vari,\vard,\variii)(0)}_{\nhz}
+K_{\beta}\norm{\nabla(\variii,\eps\vari,\eps\vard)(0)}_{\nhz}.\label{idenex:final}
\end{align}

It is worth remarking that, to establish~\eqref{idenex:final}, 
it is enough to prove~\eqref{iden:ex1} and~\eqref{iden:ex4}. 
Let us now compare~\eqref{idenex:final} with 
the estimate given in Theorem~\ref{theo:L2}. We see that the previous 
study fails to convey one feature of the proof, namely the usefulness of the additional 
smoothing effect for $\nabla\vari$. The reason is that we worked with constant coefficients. 
In the general case, an estimate for $\int_{0}^{t}\norm{\nabla\vari}_{\nhz}^{2}\,d\tau$ 
is needed to control the left hand side of~\eqref{iden:ex4} (see
Lemma~\ref{lemm:L2partialslow} below). 

Let us explain how to estimate $\int_{0}^{t}\norm{\nabla\vari}_{\nhz}^{2}\,d\tau$. 
Multiply the second equation in~\eqref{system:ex} 
by $\eps\nabla\vari$ and integrate over the strip $[0,t]\times\xD$, to obtain
$$
\int_{0}^{t}\norm{\nabla\vari}_{\nhz}^{2}\,d\tau=
-\int_{0}^{t}\scal{\partial_{t}\vard}{\eps\nabla\vari}\,d\tau.
$$
Integrating by parts both in space and time yields
\begin{align*}
\int_{0}^{t}\scal{\eps\partial_{t}\vard}{\nabla\vari}\,d\tau&=
-\int_{0}^{t}\scal{\vard}{\eps\partial_{t}\nabla\vari}\,d\tau
+\eps\bigl[\scal{\vard(\tau)}{\nabla\vari(\tau)}\bigr]_{\tau=0}^{\tau=t}\\
&=\int_{0}^{t}\scal{\cn\vard}{\eps\partial_{t}\vari}\,d\tau
+\eps\bigl[\scal{\vard(\tau)}{\nabla\vari(\tau)}\bigr]_{\tau=0}^{\tau=t}.
\end{align*}
In view of the first equation in~\eqref{system:ex}, we are left with
$$
\int_{0}^{t}\norm{\nabla\vari}_{\nhz}^{2}\,d\tau=
\int_{0}^{t}\norm{\cn\vard}_{\nhz}^{2}-\scal{\cn\vard}{\Delta\variii}\,d\tau
-\eps\bigl[\scal{\vard(\tau)}{\nabla\vari(\tau)}\bigr]_{\tau=0}^{\tau=t}.
$$
All the terms that appear in the previous identity have been estimated previously. As a 
consequence the estimate~\eqref{idenex:final} holds true if 
we include $\int_{0}^{t}\norm{\nabla\vari}_{\nhz}^{2}\,d\tau$ in its left hand side. 
By so doing, we obtain the exact analogue of~the estimate given in 
Theorem~\ref{theo:L2} (for $\r=0$ and $\pe=1$).

%The essential feature of the computation is that $\beta>1$. 
%Had we assumed instead that $\beta<1$, the corresponding inequality would have been false.

\subsection{$L^{2}$ estimates}\label{2L2esti}
Guided by the previous example, we want to prove $L^{2}$ estimates for $(\tvari,\tvard)$ 
and $(\tvariii,\eps\tvari,\eps\tvard)$. The strategy for proving both estimates is the same: 
transform the system~\eqref{system:NSi} so as to obtain 
$L^{2}$ estimates uniform in~$\eps$ by a simple integration 
by parts in which the large terms in~$1/\eps$ cancel out. Namely, we want to obtain 
systems having the form
\begin{equation}\label{normalform}
\mathcal{L}_{1}(\vard,\phi)\Ur-\mathcal{L}_{2}(\r,\pe,\phi)\Ur + \frac{1}{\eps}S(\phi)\Ur = F,
\end{equation}
where~$\mathcal{L}_{1}(\vard,\phi)-\mathcal{L}_{2}(\r,\pe,\phi)$ is a 
mixed hyperbolic/parabolic system of equations:
\begin{equation*}
L_{0}(\phi)\partial_{t} + \sum_{1\le j\le d}L_{j}(\vard,\phi)\partial_{j} 
-\sum_{1\le j,k\le d}L_{jk}(\r,\pe,\phi)\partial_{j}\partial_{k},
\end{equation*}
and the singular perturbation~$S(\phi)$ is a differential operator in the space variable which is 
skew-symmetric (with not necessarily constant coefficients).

\smallbreak
We first prove an estimate parallel to~\eqref{iden:ex3}.
\begin{proposition}\label{theo:decoupling1}
Using the same notations as in Theorem~$\ref{theo:L2}$, we have
\begin{equation}\label{theo:slowineq}
\begin{split}
&\sup_{t\in [0,T]} \norm{(\tvariii,\eps\tvari,\eps\tvard)(t)}_{\nhz} 
+\biggl(\int_{0}^{T}
\pe\norm{\nabla\tvariii}_{\nhz}^{2}+\r\norm{\eps\nabla\tvard}_{\nhz}^{2}\,dt \biggr)^{1/2} \\
&\le C(R_{0})e^{TC(R)}\norm{(\tvariii,\eps\tvari,\eps\tvard)(0)}_{\nhz}
+C(R)\biggl(\int_{0}^{T} B_{\eps}(f,\tvar)\,dt \biggr)^{1/2},
\end{split}
\end{equation}
where
\begin{equation}\label{sourceterm}
B_{\eps}(f,\tvar)\defn \norm{(f_{3},\eps f_{1},\eps
  f_{2})}_{\nhz}\norm{(\tvariii,\eps\tvari,\eps\tvard)}_{\nhz}.
\end{equation}
\end{proposition}

\begin{corollary}\label{coro:slowcomponent}
The estimate~\eqref{theo:slowineq} holds true with 
$\bigl(\int_{0}^{T} B_{\eps}(f,\tvar)\,dt\bigr)^{1/2}$ replaced by
\begin{equation*}
\int_{0}^{T}\eps \norm{(f_{1},f_{2})}_{\nhz}+\norm{f_{3}}_{\nhz}\,dt.
\end{equation*}
\end{corollary}
\begin{proof}[Proof of Corollary~$\ref{coro:slowcomponent}$ given 
Proposition~$\ref{theo:decoupling1}$] 
One has, for all $\lambda>0$,
\begin{equation*}
\biggl(\int_{0}^{T} B_{\eps}(f,\tvar)\,dt\biggr)^{1/2}\le 
\frac{1}{\lambda}\sup_{t\in [0,T]}\norm{(\tvariii,\eps\tvari,\eps\tvard)}_{\nhz} + 
\lambda\int_{0}^{T}\norm{(f_{3},\eps f_{1},\eps f_{2})}_{\nhz}\,dt.
\end{equation*}
Using~\eqref{theo:slowineq} we obtain the desired result by taking 
$\lambda$ large enough.
%\qed
\end{proof}
\begin{nota}\label{nota:C0C}Within the proofs of
  Proposition~\ref{theo:decoupling1} and~\ref{prop:fastL2}: Rather
  than writing $C(R)$ and $C(R_{0})$ in full, we will use the following abbreviations:
~$C=C(R)$ and $C_{0}=C(R_{0})$. So that $C$ denotes generic constants which depend only on 
$R\defn \norm{\bigl(\phi,\partial_{t}\phi,
\nabla\phi,\ncp\nabla^{2}\phi,\vard,\nabla\vard\bigr)}_{L^{\infty}([0,T]\times\xD)}$ and 
$C_{0}$ denotes generic constants which depend only
on~$R_{0}\defn \norm{\phi(0)}_{L^{\infty}(\xD)}$.  
As usual, the values of $C_{0}$ and $C$ may vary from relation to relation.
\end{nota}

\begin{proof}[Proof of Proposition~$\ref{theo:decoupling1}$]
All the computations 
given below are meaningful since we concentrate on regular solutions
of~\eqref{system:NSi}. Introduce 
$$
\widetilde{\zeta}\defn \eps g_{1}\beta_{3} \tvari -g_{3}\beta_{1}\tvariii,\quad 
\widetilde{v}_{\eps}\defn\eps\tvard \quad\mbox{and}\quad
\widetilde{\Ur}\defn(\widetilde{\zeta},\widetilde{v}_{\eps},\tvariii)^{t}.
$$ 
where we simply write~$g_{i}$,~$\beta_{i}$ instead of~$g_{i}(\phi)$,~$\beta_{i}(\phi)$. 
We first show that~$\widetilde{\Ur}$ solves a system of the form~\eqref{normalform}. 
Expressing~$\tvari$ in terms of~$\widetilde{\zeta}$ and~$\tvariii$, 
replacing~$\tvard$ by~$\eps^{-1}\widetilde{v}_{\eps}$ and performing a little algebra, yields
\begin{align*}
&\partial_{t}{\widetilde{\zeta}}
+\vard\cdot\nabla{\widetilde{\zeta}}+\frac{\beta_{3}
-\beta_{1}}{\eps}\cn\widetilde{v}_{\eps}=f_{1}',\\
g_{2}\bigl(&\partial_{t}{\widetilde{v}_{\eps}}
+\vard\cdot\nabla{\widetilde{v}_{\eps}}\bigr)
+ \frac{1}{\eps}\nabla \bigl(\gamma_{1}\widetilde{\zeta}\bigr) 
+\frac{1}{\eps}\nabla\bigl(\gamma_{2}\tvariii\bigr)
-\r\beta_{2}\Delta \widetilde{v}_{\eps}-\r\beta_{2}^{\sharp}\nabla\cn\widetilde{v}_{\eps}=
\eps f_{2},\notag \\
g_{3}\bigl(&\partial_{t}{\tvariii}
+\vard\cdot\nabla{\tvariii}\bigr)
+ \frac{1}{\eps}\cn \widetilde{v}_{\eps}
-\pe\beta_{3}\Delta\tvariii =  f_{3},
\end{align*}
where~${\displaystyle 
\gamma_{1}\defn \frac{1}{g_{1}\beta_{3}}}$, 
$ {\displaystyle\gamma_{2}\defn \frac{g_{3}\beta_{1}}{g_{1}\beta_{3}}}$ and
\begin{align*}
f_{1}'&\defn\eps\beta_{3}f_{1} -\beta_{1}f_{3}+\pe\beta_{3}\nabla\beta_{1}\cdot\nabla\tvariii \\
&\quad+\eps\tvari\bigl(\partial_{t}{(g_{1}\beta_{3})}+\vard\cdot\nabla{(g_{1}\beta_{3})}\bigr)-
\tvariii\bigl(\partial_{t}{(g_{3}\beta_{1})}+\vard\cdot(\nabla{g_{3}\beta_{1}})\bigr).
\end{align*}
In order to symmetrize these equations we make use of the structural 
assumption~\eqref{mainassumption}, that is 
$\beta_{3}-\beta_{1}>0$. Multiply the first equation [resp.\ the third] 
by~$\gamma_{1}(\beta_{3}-\beta_{1})^{-1}$ [resp.~$\gamma_{2}$]. 
By so doing, we obtain that~$\widetilde{\Ur}$ solves 
\begin{equation}\label{normalformUr}
\mathcal{L}_{1}(\vard,\phi)\widetilde{\Ur}-\mathcal{L}_{2}(\r,\pe,\phi)\widetilde{\Ur}
+\frac{1}{\eps}S(\phi)\widetilde{\Ur} =F,
\end{equation}
with~$F\defn\bigl(\gamma_{1}(\beta_{3}-\beta_{1})^{-1}f_{1}',\eps f_{2},\gamma_{2}f_{3}\bigr)^{T}$, 
$\mathcal{L}_{1}(\vard,\phi) \defn L_{0}(\phi)\partial_{t} + \sum L_{j}(\vard,\phi)\partial_{j}$ where
$$
L_{0}\defn
\begin{pmatrix}
\frac{\gamma_{1}}{\beta_{3}-\beta_{1}}&0&0\\
0 & g_{2}I_{d} &0 \\
0 &0 & \gamma_{2}g_{3}\end{pmatrix}\quad\mbox{and}\quad
\quad L_{j}\defn v_{j} L_{0}.
$$
Moreover, the singular and viscous perturbations are defined by
$$
S=
\begin{pmatrix} 
0 & \gamma_{1}\cn & 0 \\ 
\nabla (\gamma_{1}\,\cdot\,) & 0 & \nabla(\gamma_{2}\,\cdot\,) \\
0 & \gamma_{2}\cn & 0
\end{pmatrix}, ~~ \mathcal{L}_{2}=\begin{pmatrix} 0& 0 & 0 \\
0&\r\beta_{2}\Delta+\r\beta_{2}^{\sharp}\nabla\cn&0\\
0&0&\hspace{-5pt}\pe\gamma_{2}\beta_{3}\Delta
\end{pmatrix}.
$$

The end of the proof proceeds by multiplying by~$\widetilde{\Ur}$ and 
integration by parts. 
The point is that the large terms in~$\eps^{-1}$ cancel out for~$S$ is
skew-symmetric.  
Let $t\in (0,T]$. If we multiply~\eqref{normalformUr} by~$\widetilde{\Ur}$ 
and integrate from~$0$ to~$t$, we obtain
\begin{equation}\label{ineq:L2}
\nscal{L_{0}(\phi)\widetilde{\Ur}}{\widetilde{\Ur}}(t) 
-2\int_{0}^{t} \nscal{\mathcal{L}_{2}(\phi)\widetilde{\Ur}}{\widetilde{\Ur}}\,d\tau
=\nscal{L_{0}(\phi)\widetilde{\Ur}}{\widetilde{\Ur}}(0)+ Y(t),
\end{equation}
where~$\scal{\cdot}{\cdot}$ denotes the scalar product in~$L^{2}(\xD)$ and
$$
Y(t)\defn \int_{0}^{t} \Bscal{\Bigl\{\partial_{t}L_{0}(\phi)+
\sum_{1\le j\le d}\partial_{j}L_{j}(\vard,\phi)\Bigr\}\widetilde{\Ur}}{\widetilde{\Ur}}\,d\tau
+ 2\int_{0}^{t}\scal{F}{\widetilde{\Ur}}\,d t.
$$
Here we used the symmetry of the matrices~$L_{0}$ and~$L_{j}$. 

\smallbreak
We begin by estimating the right-hand side. One easily gathers
$$
\Bigl\lVert \partial_{t} L_{0}(\phi) + \sum_{1\le j\le d}\partial_{j}L_{j}(\vard,\phi)
\Bigr\rVert_{L^{\infty}}\le C, 
$$
(where $C$ is as in Notation~\ref{nota:C0C}) and 
\begin{equation}\label{st}
\lvert\scal{F}{\widetilde{\Ur}}\rvert \le
\norm{F}_{\nhz}\norm{\widetilde{\Ur}}_{\nhz} 
\le C \norm{(f_{1}',\eps f_{2}, f_{3})}_{\nhz}\norm{\widetilde{\Ur}}.
\end{equation}
Furthermore, since by definition $f_{1}'$ is a linear combination 
of $\eps f_{1}$, $f_{3}$, $\eps\tvari$, $\tvariii$ and
$\pe\nabla\tvariii$ whose coefficients are estimated in $L^{\infty}$ norm 
by a constant depending only on $R$, we get
\begin{equation}\label{st:2}
\norm{f_{1}'}\le C 
\norm{(f_{3},\eps f_{1})}_{\nhz} 
+C\norm{(\tvariii,\eps\tvari)}_{\nhz}+C\norm{\pe\nabla\tvariii}_{\nhz}.
\end{equation}
Moreover, observe that 
\begin{equation}\label{ineqL2:111}
(\tvariii,\eps\tvari,\eps\tvard)^{t}=M(\phi)\widetilde{\Ur},
\end{equation}
where~$M(\phi)$ is an~$n\times n$ invertible matrix. As a consequence,
there exists a constant $C$ such that 
$C^{-1}\norm{\widetilde{\Ur}}_{\nhz}\le\norm{(\tvariii,\eps\tvari,\eps\vard)}_{\nhz}\le C\norm{\widetilde{\Ur}}_{\nhz}$. 
Hence, the estimates~\eqref{st} and~\eqref{st:2} result in
\begin{equation}
\lvert\scal{F}{\widetilde{\Ur}}\rvert \le C B_{\eps}(f,\widetilde{U}) + 
C\norm{\widetilde{\Ur}}_{\nhz}^{2}+
C\norm{\pe\nabla\tvariii}_{\nhz}\norm{\widetilde{\Ur}}_{\nhz},
\end{equation}
where~$B_{\eps}(f,\tvar)$ is given by~\eqref{sourceterm}. 
For all~$\lambda\ge 1$, the last term in the right 
hand-side of the previous estimate is bounded by
\begin{equation*}
\lambda C\norm{\widetilde{\Ur}}_{\nhz}^{2}+\frac{\pe^{2}}{\lambda}\norm{\nabla\tvariii}_{\nhz}^{2}.
\end{equation*}
Consequently the right-hand side of~\eqref{ineq:L2} is less than
\begin{equation}\label{slow:RHS}
C_{0}\norm{\widetilde{\Ur}(0)}_{\nhz}^{2}+\lambda C\int_{0}^{t}
\norm{\widetilde{\Ur}}_{\nhz}^{2}\,d\tau
+\frac{\pe^{2}}{\lambda}\int_{0}^{t}\norm{\nabla\tvariii}_{\nhz}^{2}\,d\tau 
+C\int_{0}^{t}B_{\eps}(f,\tvar)\,d\tau.
\end{equation}

\smallbreak
We want to estimate now the left-hand side of~\eqref{ineq:L2}. In this regard, 
the most direct estimates show that
\begin{equation*}
\nscal{L_{0}(\phi)\widetilde{\Ur}}{\widetilde{\Ur}}\ge
\norm{L_{0}(\phi)^{-1}}_{L^{\infty}}^{-1}\norm{\widetilde{\Ur}}_{\nhz}^{2}.
\end{equation*}
Similarly, using Lemma~\ref{lemm:Raphael} below, one has
\begin{equation*}
-\nscal{\mathcal{L}_{2}(\phi)\widetilde{\Ur}}{\widetilde{\Ur}} 
\ge  K \r m \norm{\nabla \widetilde{v}_{\eps}}^{2}
+K \pe m\norm{\nabla\tvariii}^{2}
-C\norm{\widetilde{\Ur}}_{\nhz}^{2},
\end{equation*}
where $m\defn
\min\bigl\{\norm{\beta_{2}^{-1}}_{L^{\infty}}^{-1},\norm{(\beta_{2}+\beta_{2}^{\sharp})^{-1}}_{L^{\infty}}^{-1},\norm{\beta_{3}^{-1}}_{L^{\infty}}^{-1}\bigr\}$ 
and $K$ is a generic constant depending only on the dimension.

To estimate the~$L^{\infty}$ norm of~$L_{0}^{-1}$, 
$\beta_{2}^{-1}$, $(\beta_{2}+\beta_{2}^{\sharp})^{-1}$
and~$\beta_{3}^{-1}$, write 
\begin{equation}\label{boundLinfty}
\begin{split}
\norm{F(\phi(t))}_{L^{\infty}} &\le \norm{F(\phi(0))}_{L^{\infty}}+\int_{0}^{t}
\norm{\partial_{t}F(\phi(\tau))}_{L^{\infty}}\,d\tau\\ 
&\le C_{0}+t\sup_{\tau\in [0,t]}\norm{F'(\phi(\tau))}_{L^{\infty}}
\norm{\partial_{t}\phi(\tau)}_{L^{\infty}}\\
&\le C_{0}+TC\le C_{0}e^{TC}.
\end{split}
\end{equation}
Consequently, the left hand-side of~\eqref{ineq:L2} is greater than
\begin{equation}\label{slow:LHS}
C_{0}e^{-TC} \Bigl\{ \norm{\widetilde{\Ur}(t)}_{\nhz}^{2} + 
\int_{0}^{t} \r\norm{\nabla \widetilde{v}_{\eps}}^{2} + \pe\norm{\nabla\tvariii}^{2} \,d\tau \Bigr\}
-C\int_{0}^{t}\norm{\widetilde{\Ur}}_{\nhz}^{2}\,d\tau.
\end{equation}
From this and the assumption~$\pe\le 1$, we see that one can absorb the third 
term in~\eqref{slow:RHS} by 
taking $\lambda$ large enough. Then, Gronwall's Lemma implies
\begin{equation*}
\norm{\widetilde{\Ur}(t)}_{\nhz}^{2}
+\int_{0}^{t}\r\norm{\nabla\widetilde{v}_{\eps}}^{2}+
\pe\norm{\nabla\tvariii}^{2} \,d\tau
\le C_{0}e^{TC}\norm{\widetilde{\Ur}(0)}_{\nhz}^{2}+C\int_{0}^{t}B_{\eps}(f,\tvar)\,d\tau.
\end{equation*}

We next claim that the previous estimate 
holds true with~$\widetilde{\Ur}$ replaced by~$(\tvariii,\eps\tvari,\eps\tvard)$. Indeed, 
using the identity~\eqref{ineqL2:111} and the bound~\eqref{boundLinfty}, one has 
$\norm{(\tvariii,\eps\tvari,\eps\tvard)}_{\nhz}\le C_{0}e^{TC}\norm{\widetilde{\Ur}}_{\nhz}$. 
To complete the proof, take the square root of the inequality thus obtained and 
take the supremum over $t\in [0,T]$. \end{proof}

Let us briefly recall a standard estimate used in the previous proof.
\begin{lemma}\label{lemm:Raphael}
There exists two constants $K_{i}=K_{i}(d)$ such that for all Lipschitz functions
$\zeta=\zeta(x)$ and $\eta=\eta(x)$ such that $\zeta>0$ and
$\eta+\zeta>0$, and for all vector field $u\in H^{1}(\xD)$, there holds
\begin{equation}
-\scal{\zeta\Delta u+ \eta\nabla\cn u}{u}_{H^{-1}\times H^{1}}\ge K_{1} m 
\norm{\nabla u}_{L^{2}}^{2} - \frac{K_{2}M^{2}}{m} \norm{u}_{L^{2}}^{2},
\end{equation}
where $m=\inf_{x\in\xD}\{ \zeta(x),\zeta(x)+\eta(x)\}$ and 
$M=\norm{\nabla\zeta}_{L^{\infty}}+\norm{\nabla\eta}_{L^{\infty}}$.
\end{lemma}
\begin{proof} Decompose $u$ as $\underline{u}+\nabla\phi\defn \Qr u
  +(\id-\Qr)u$ where $\Qr$ is the Leray projector onto divergence free 
vector field \footnote{This trick is due to R. Danchin.}, so that
$$
\zeta\Delta u+ \eta\nabla\cn u
=\zeta\Delta \underline{u}+(\zeta+\eta)\nabla\Delta\phi.
$$
It results from $\cn \underline{u}=0$ that
$$
-\scal{\zeta\Delta u+ \eta\nabla\cn u}{u}_{H^{-1}\times H^{1}}
=\scal{\zeta\nabla
  \underline{u}}{\nabla\underline{u}}+\scal{(\zeta+\eta)\Delta\phi}{\Delta
  \phi} +R,
$$
where $R$ is such that $\la R\ra\le K_{3} M
\norm{\nabla u}_{\nhz}\norm{u}_{L^{2}}$. 

To infer the
desired bound, we use two usual inequalities. For all $\lambda\ge 1$, one has 
$\la R\ra \le (m/\lambda)\norm{\nabla u}_{\nhz}^{2} + (\lambda
K_{3}^{2}M^{2}/m)\norm{u}_{\nhz}^{2}$. We next recall 
the most simplest of all Calder\'on--Zygmund estimates: 
$\norm{\nabla^{2}\phi}_{L^{2}}\le K_{4}\norm{\Delta\phi}_{L^{2}}$.
\end{proof}

We next prove an estimate parallel to~\eqref{iden:ex1}.
\begin{proposition}\label{prop:fastL2}
Using the same notations as in Theorem~$\ref{theo:L2}$, we have
\begin{align*}
&\sup_{t\in[0,T]} \norm{(\tvari,\tvard)(t)}_{\nhz} 
+\biggl(\int_{0}^{T}\r\norm{\nabla\tvard}_{\nhz}^{2}
+\pe\norm{\cn\tvard}_{\nhz}^{2}\, dt\biggr)^{1/2}\\
&\qquad\le C(R_{0})e^{TC(R)}\norm{(\tvari,\tvard)(0)}_{\nhz}
+ C(R)\int_{0}^{T}\norm{(f_{1},f_{2})}_{\nhz} +\pe\norm{f_{3}}_{\nh{1}}\,dt\\[0.5ex]
&\qquad+ C(R_{0})e^{TC(R)}\biggl\{\sup_{t\in [0,T]} 
\norm{\tvariii(t)}_{H^{1}_{\ncp}}+\biggl(\int_{0}^{T}\pe
\norm{\nabla\tvariii}_{H^{1}_{\ncp}}^{2}\,dt \biggr)^{1/2}\biggr\},
\end{align*}
where~$\ncp\defn\sqrt{\r+\pe}$ and $\norm{\cdot}_{H^{1}_{\ncp}}=
\norm{\cdot}_{\nhz}+\ncp\norm{\cdot}_{\nh{1}}$.
\end{proposition}
\begin{remark}
From now on, we make intensive and implicit uses of the assumptions
$\r\le 1$ and $\pe\le 1$. In particular we freely use obvious estimates
like $\r\le\sqrt{\r}\le \cpe$. Similarly, we use the estimate
$\cpe\norm{\nabla u}_{\nhz}\le \norm{u}_{\nhsc{1}}$ without mentioning it
explicitly in the proofs.
\end{remark}
\begin{proof}
%The proof of Proposition~\ref{prop:fastL2} is in two steps. 
Introduce $\tvard_{e}\defn\tvard-\pe\beta_{1}\nabla\tvariii$. 
Performing a little algebra we find that~$\widetilde{\mathcal{V}}\defn
(\tvari,\tvard_{e})^{t}$ satisfies
\begin{equation*}
\begin{pmatrix}
g_{1}&0\\
0&g_{2}I_{d}
\end{pmatrix}\bigl(\partial_{t}\widetilde{\mathcal{V}}+\vard\cdot
\nabla\widetilde{\mathcal{V}}\bigr)
+\frac{1}{\eps}\begin{pmatrix}0&\cn\\ \nabla&0\end{pmatrix}
\widetilde{\mathcal{V}}
-
\begin{pmatrix}
0 & 0 \\ 0 & \mathcal{L}_{22}
\end{pmatrix}\widetilde{\mathcal{V}}=F,
\end{equation*}
where
$$
\mathcal{L}_{22}\defn \pe g_{2}\beta_{1}\nabla(g_{3}^{-1}\cn\,\cdot\,)+\r\beta_{2}\Delta+\r\beta_{2}^{\sharp}\nabla\cn,
$$
and $F=(f_{1},f_{2}+f_{2}')^{t}$ with
\begin{equation}\label{L2esti:fast.source}
\begin{split}
f_{2}'&\defn \r\pe\beta_{2}\Delta(\beta_{1}\nabla\tvariii) 
+\r\pe\beta_{2}^{\sharp}\nabla\cn(\beta_{1}\nabla\tvariii)
-\pe g_{2}\bigl(\partial_{t}\beta_{1}+\vard\cdot\nabla {\beta_{1}}\bigr)\nabla\tvariii
\\
&\quad + \pe g_{2}\beta_{1}\nabla\bigl(\pe g_{3}^{-1}\cn(\beta_{1}\nabla\tvariii)-
\pe g_{3}^{-1}\beta_{3}\Delta\tvariii-g_{3}^{-1}f_{3}\bigr)
+\pe g_{2}\beta_{1}(\nabla\vard)^{t}\nabla\tvariii.
\end{split}
\end{equation}
The thing of greatest interest here is that~$-\mathcal{L}_{22}$ is a 
differential operator whose leading symbol is greater than 
$C\pe \xi^{t}\cdot\xi + C \r \la\xi\ra^{2}$.

As before, the proof proceeds by multiplying
by~$\widetilde{\mathcal{V}}$ and integrating on the strip $[0,t]\times\xD$.
Then the analysis establishing that the right [resp.\ left] hand side
of~\eqref{ineq:L2} is smaller than~\eqref{slow:RHS} [resp.\ greater
than~\eqref{slow:LHS}] also gives:
\begin{equation}\label{L2esti:fast.ineq}
\begin{split}
&\norm{\widetilde{\mathcal{V}}(t)}_{\nhz}^{2}
+\int_{0}^{t}\r\norm{\nabla \tvard_{e}}_{\nhz}^{2}+\pe\norm{\cn \tvard_{e}}_{\nhz}^{2} \, d\tau\\
&\qquad\le
C_{0}e^{TC}\norm{\widetilde{\mathcal{V}}(0)}_{\nhz}^{2}+C\int_{0}^{t}
\norm{\widetilde{\mathcal{V}}}_{\nhz}^{2}\,d\tau+
C\int_{0}^{t}\lvert\scal{F}{\widetilde{\mathcal{V}}}\rvert\,d\tau.
\end{split}
\end{equation}
Let us estimate the last term in the right-hand side of~\eqref{L2esti:fast.ineq}. 
Using~\eqref{L2esti:fast.source} one can decompose~$f_{2}'$ as~$f_{2,1}'+\sqrt{\pe}\nabla f_{2,2}'$ where 
\begin{align*}
\norm{f_{2,1}'}_{\nhz}&\le C \pe\norm{f_{3}}_{\nh{1}}
+C\pe\norm{\nabla\tvariii}_{\nhz}+C\sqrt{\pe}(\r+\pe)\norm{\nabla^{2}\tvariii}_{\nhz},\\
\norm{f_{2,2}'}_{\nhz}&\le C\pe\norm{\nabla\tvariii}_{\nhz}+
C\sqrt{\pe}(\r+\pe)\norm{\nabla^{2}\tvariii}_{\nhz}.
\end{align*}
By notation, $\scal{F}{\widetilde{\mathcal{V}}}=\scal{f_{1}}{\tvari}
+\scal{f_{2}}{\tvard_{e}}+\scal{f_{2,1}'}{\tvard_{e}}-\scal{f_{2,2}'}{\sqrt{\pe}\cn \tvard_{e}}$. 
We thus get, for all $\lambda\ge 1$,
\begin{equation}\label{L2esti:fast.ineqR}
\begin{aligned}
\int_{0}^{t}\bigl\lvert\scal{F}{\widetilde{\mathcal{V}}}\bigr\rvert \,d\tau & \le\lambda C 
\biggl(\int_{0}^{t}\norm{(f_{1},f_{2})}_{\nhz}+\pe\norm{f_{3}}_{\nh{1}}\,d\tau\biggr)^{2} \\
&\quad+ \frac{1}{\lambda}\int_{0}^{t}\pe\norm{\cn \tvard_{e}}_{\nhz}^{2}\,d\tau +\frac{1}{\lambda}
\sup_{\tau\in [0,t]}\norm{\widetilde{\mathcal{V}}(t)}_{\nhz}^{2}\\
&\quad + \lambda C \int_{0}^{t}
\norm{(\pe\nabla\tvariii,\sqrt{\pe}(\r+\pe)\nabla^{2}\tvariii}_{\nhz}^{2}\,d\tau.
\end{aligned}
\end{equation}

Replacing $\tvard_{e}$ with $\tvard$ in~\eqref{L2esti:fast.ineq}--\eqref{L2esti:fast.ineqR}, 
taking $\lambda$ large enough and applying Gronwall's lemma leads to the expected result.
%\qed
\end{proof}

\subsection{End of the proof of Theorem~\ref{theo:L2}.}
From now on, we consider a time $0<T\le 1$, a fixed triple of parameter
$\pa=(\eps,\r,\pe)\in\PA$ and a solution $\tvar=(\tvari,\tvard,\tvariii)$ of the
system~\eqref{system:NSi}. We denote by $R_{0}$ and $R$ the norms
defined in the statement of Theorem~\ref{theo:L2}
(see~\eqref{defi:Rcoef}). We also set $\cpe\defn\sqrt{\r+\pe}$.

Introduce the functions $N,w,y,Y,Z\colon [0,T]\rightarrow \xR_{+}$
given by 
\begin{equation}\label{defi:N}
N(t)\defn \norm{\tvar}_{\Hr^{0}_{\pa}(t)},
\end{equation}
and
\begin{equation*}
w(t)\defn\sup_{\tau\in[0,t]} \norm{(\tvari,\tvard)(\tau)}_{\nhz} 
+\biggl(\int_{0}^{t}\r\norm{\nabla\tvard}_{\nhz}^{2}
+\pe\norm{\cn\tvard}_{\nhz}^{2}\, d\tau\biggr)^{1/2},
\end{equation*}

\begin{equation*}
y(t)\defn 
\sup_{\tau\in [0,t]}\norm{(\tvariii,\eps\tvari,\eps\tvard)(\tau)}_{\nhz} 
+\biggl(\int_{0}^{t}
\pe\norm{\nabla\tvariii}_{\nhz}^{2}+\r\norm{\eps\nabla\tvard}_{\nhz}^{2}\,d\tau \biggr)^{1/2},
\end{equation*}

\begin{equation*}
Y(t)\defn 
\sup_{\tau\in [0,t]}\norm{\nabla(\tvariii,\eps\tvari,\eps\tvard)(\tau)}_{\nhz} 
+\biggl(\int_{0}^{t}
\pe\norm{\nabla^{2}\tvariii}_{\nhz}^{2}+\r\norm{\eps\nabla^{2}\tvard}_{\nhz}^{2}\,d\tau \biggr)^{1/2},
\end{equation*}
\begin{equation*}
z(t)\defn \left(\int_{0}^{t}\norm{\nabla\tvari}_{\nhz}\,d\tau\right)^{1/2}.
\end{equation*}

Note that 
$N(T)\le w(T)+\cpe z(T)+y(T)+\cpe Y(T)$. 
Since we have estimated $w(T)$ and $y(T)$, it remains only to
estimate $\cpe Y(T)$ and $\cpe z(T)$. Parallel to~\eqref{iden:ex4}, we begin by
establishing an estimate for $Y$.
\begin{lemma}\label{lemm:L2partialslow}
There exists a constant $C_{0}$ depending only on
$R_{0}$ and a constant $C$ depending only on $R$ such that for all 
$t\in [0,T]$ and all $\lambda\ge 1$,
\begin{equation}\label{1}
\begin{split}
Y(t)^{2} &\le C_{0} e^{T C} Y(0)^{2}
+\lambda C\int_{0}^{t} Y(\tau)^{2}\,d\tau +\frac{C}{\lambda}\int_{0}^{t}\norm{(\cn\tvard,\nabla\tvari)}_{\nhz}^{2}\,d\tau\\
&\quad 
+ \biggl(C\int_{0}^{t}\norm{f_{3}}_{\nh{1}}+\norm{(f_{1},f_{2})}_{H^{1}_{\eps}}\,d\tau\biggr)^{2}.
\end{split}
\end{equation}
\end{lemma}
\begin{proof}
The proof, if tedious, is elementary. Indeed, the strategy for
proving~\eqref{1} consists in differentiating the
system~\eqref{system:NSi} so as to apply Proposition~\ref{theo:decoupling1} with
$(\tvariii,\eps\tvari,\eps\tvard)$ replaced by 
$\nabla(\tvariii,\eps\tvari,\eps\tvard)$. 

Let $1\le j\le d$ and set 
$$
\tvar_{j}\defn(\tvari_{j},\tvard_{j},\tvariii_{j})\defn
(\partial_{j}\tvari,\partial_{j}\tvard,\partial_{j}\tvariii).
$$ 
We commute\footnote{Which means that we commute $\partial_{j}$ with the
  $i{\rm th}$ equation ($i=1,2,3$) premultiplied by $g_{i}^{-1}$ and we
  next multiply the result by $g_{i}$.} the $i{\rm th}$ equation in~\eqref{system:NSi} with 
$g_{i}(\partial_{j} g_{i}^{-1}\cdot)$. It yields
\begin{equation*}
\left\{
\begin{aligned}
&g_{1}(\phi)\bigl(\partial_{t}{\tvari_{j}}+\vard\cdot\nabla {\tvari_{j}}\bigr)
+\frac{1}{\eps}\cn \tvard_{j}-\frac{\pe}{\eps} \cn(\beta_{1}(\phi)\nabla\tvariii_{j})=F_{1},\\
&g_{2}(\phi)\bigl(\partial_{t}{\tvard_{j}}+\vard\cdot\nabla {\tvard_{j}}\bigr)
+\frac{1}{\eps}\nabla\tvari_{j}-\r\beta_{2}(\phi)\Delta\tvard_{j}-\r\beta_{2}^{\sharp}(\phi)\nabla\cn\tvard_{j}=F_{2},\\
&g_{3}(\phi)\bigl(\partial_{t}{\tvariii_{j}}
+\vard\cdot\nabla{\tvariii_{j}}\bigr)
+\cn\tvard_{j} -\pe\beta_{3}(\phi)\Delta\tvariii_{j}=F_{3},
\end{aligned}
\right.
\end{equation*}
where, for $i=1,2,3$, the source term $F_{i}$ is given by
$$
F_{i}\defn F_{i}^{1}+F_{i}^{2}\defn g_{i}\partial_{j}(g_{i}^{-1}f_{i})+g_{i}\tilde{F}_{i}^{2},
$$
with
\begin{alignat*}{2}
\tilde{F}_{1}^{2}&\defn && -\partial_{j}v\cdot\nabla\tvari-\frac{1}{\eps}
\partial_{j}g_{1}^{-1}\cn\tvard + \frac{\pe}{\eps}\partial_{j}g_{1}^{-1}\cn(\beta_{1}\nabla\tvariii)
+\frac{\pe}{\eps}g_{1}^{-1}\cn(\partial_{j}\beta_{1}\nabla\tvariii),\\
\tilde{F}_{2}^{2}&\defn && -\partial_{j}v\cdot\nabla\tvard-\frac{1}{\eps}
\partial_{j}g_{2}^{-1}\nabla\tvari+\r\partial_{j}(g_{2}^{-1}\beta_{2})\Delta\tvard
+\r\partial_{j}(g_{2}^{-1}\beta_{2}^{\sharp})\nabla\cn\tvard,\\
\tilde{F}_{3}^{2}&\defn && -\partial_{j}v\cdot\nabla\tvariii-\partial_{j}g_{3}^{-1}\cn\tvard+\pe\partial_{j}(g_{3}^{-1}\beta_{3})\Delta\tvariii.
\end{alignat*}

Proposition~\ref{theo:decoupling1} implies that 
\begin{equation}\label{esti:partialY}
Y(t)^{2} \le C_{0} e^{T C}\norm{\widetilde{\mathcal{U}}(0)}_{\nhz}^{2}
+C\sum_{1\le j\le d}\int_{0}^{t}B_{\eps}(F,\tvar_{j})\,d\tau.%+
\end{equation}
where $\widetilde{\mathcal{U}}\defn(\nabla\tvariii,\eps\nabla\tvari,\eps\nabla\tvard)$ and $B_{\eps}$ is defined by~\eqref{sourceterm}. 
We now have to estimate $B_{\eps}(F,\tvar_{j})$. The source terms are directly estimated by
\begin{align*}
\norm{\bigl(F^{1}_{3},\eps F^{1}_{2},\eps F^{1}_{3}\bigr)}_{\nhz}&\le 
C\norm{f_{3}}_{\nh{1}}+C\norm{(f_{1},f_{2})}_{H^{1}_{\eps}},\\
\norm{\bigl(\tilde{F}^{2}_{3},\eps\tilde{F}^{2}_{2},\eps\tilde{F}^{2}_{3}\bigr)}_{\nhz}&\le 
C \norm{\widetilde{\mathcal{U}}}_{\nhz}+C\norm{(\cn\tvard,\nabla\tvari)}_{\nhz}+
C\norm{(\eps\r\nabla^{2}\tvard,\pe\nabla^{2}\tvariii)}_{\nhz}.
\end{align*}
The last estimate implies that, for all positive $\lambda$ and $\lambda'$, one has
\begin{align*}
B_{\eps}(F^{2},\tvar_{j}) 
&\le (1+\lambda+\lambda') C \norm{\widetilde{\mathcal{U}}}_{\nhz}^{2}\\
&\quad +\frac{C}{\lambda}\norm{(\cn\tvard,\nabla\tvari)}_{\nhz}^{2}+\frac{C}{\lambda'}
\norm{(\eps\r\nabla^{2}\tvard,\pe\nabla^{2}\tvariii)}_{\nhz}^{2}.
\end{align*}
Moreover, the term $\int_{0}^{t}\lvert B_{\eps}(F^{1},\tvar_{j})\rvert\,d\tau$ 
is estimated as in the proof of Corollary~\ref{coro:slowcomponent}. 
With these estimates in hands, we find that the last term in the right-side 
of~\eqref{esti:partialY} is less then
\begin{equation}\label{esti:partialYRHS}
\begin{split}
&\int_{0}^{t}(1+\lambda+\lambda') C \norm{\widetilde{\mathcal{U}}}_{\nhz}^{2}\,d\tau +
\frac{C}{\lambda''}\sup_{\tau\in [0,t]}\norm{\widetilde{\mathcal{U}}}_{\nhz}^{2}\\
&+\int_{0}^{t}\frac{C}{\lambda}\norm{(\cn\tvard,\nabla\tvari)}_{\nhz}^{2}\,d\tau
+\frac{C}{\lambda'}\int_{0}^{t}\norm{(\eps\r\nabla^{2}\tvard,\pe\nabla^{2}\tvariii)}_{\nhz}^{2}\,d\tau\\
&+\lambda''  \biggl(C\int_{0}^{t}\norm{f_{3}}_{\nh{1}}+\norm{(f_{1},f_{2})}_{H^{1}_{\eps}}\,d\tau\biggr)^{2}.
\end{split}
\end{equation}
Yet, by definition, $\sup_{\tau\in [0,t]}\norm{\widetilde{\mathcal{U}}}_{\nhz}^{2}+
\int_{0}^{t}\norm{(\eps\r\nabla^{2}\tvard,\pe\nabla^{2}\tvariii)}_{\nhz}^{2}\,d\tau\le
Y(t)^{2}$. 
Hence, taking $\lambda'=\lambda''$ large enough, one can absorb the second and fourth terms in~\eqref{esti:partialYRHS} 
in the left hand side of~\eqref{esti:partialY}, thereby obtaining the
desired estimate~\eqref{1}.
\end{proof}

To obtain a closed set of inequalities it remains to estimate~$z(t)$.
\begin{lemma}
There  exists a constant $C_{0}$ depending only on
$R_{0}$ and a constant $C$ depending only on $R$ such that for all
$t\in [0,T]$,
\begin{equation}\label{esti:zL2}
\cpe^{2} z(t)^{2}\le  C_{0}e^{TC}\bigl\{ \cpe^{2}
Y(t)^{2}+w(t)^{2}\bigr\}
+C\biggr(\int_{0}^{T}\norm{(f_{1},f_{2})}_{\nhz} \,dt\biggr)^{2}.
\end{equation}
\end{lemma}
\begin{proof}
Set $\dom\defn [0,t]\times\xD$ and denote by~$\dscal{}{}$ the scalar product in~$L^{2}\bigl(\dom\bigr)$. 
Multiplying the second equation in~\eqref{system:NSi} by 
$\eps\ncp^{2}g_{2}^{-1}\nabla\tvari$, and integrating over~$\dom$ yields
\begin{equation}\label{112}
\begin{split}
&\ncp^{2}\dscal{g_{2}^{-1}\nabla\tvari}{\nabla\tvari} = \\
&\qquad\qquad- \eps\ncp^{2}\dscal{(\partial_{t}+\vard\cdot\nabla)\tvard}{\nabla\tvari}
+\eps\r\ncp^{2}\dscal{g_{2}^{-1}\beta_{2}\Delta\tvard}{\nabla\tvari}\\
&\qquad\qquad+\eps\r\ncp^{2}\dscal{g_{2}^{-1}\beta_{2}^{\sharp}\nabla\cn\tvard}{\nabla\tvari}+
\eps\ncp^{2}\dscal{g_{2}^{-1}f_2}{\nabla\tvari}.
\end{split}
\end{equation}
The most direct estimates show that
\begin{align*}
\dscal{g_{2}^{-1}\nabla\tvari}{\nabla\tvari} &\ge \norm{g_{2}}_{L^{\infty}(\dom)}^{-1} 
\norm{\nabla\tvari}_{L^{2}(\dom)}^{2},\\
\eps\r\ncp^{2}\dscal{g_{2}^{-1}\beta_{2}\Delta\tvard}{\nabla\tvari}&\le 
\norm{g_{2}^{-1}\beta_{2}}_{L^{\infty}(\dom)}
\norm{\eps\r\ncp\nabla^{2}\tvard}_{L^{2}(\dom)}\norm{\ncp\nabla\tvari}_{L^{2}(\dom)},\\
\eps\r\ncp^{2}\dscal{g_{2}^{-1}\beta_{2}^{\sharp}\nabla\cn\tvard}{\nabla\tvari}&\le 
\norm{g_{2}^{-1}\beta_{2}^{\sharp}}_{L^{\infty}(\dom)}
\norm{\eps\r\ncp\nabla^{2}\tvard}_{L^{2}(\dom)}\norm{\ncp\nabla\tvari}_{L^{2}(\dom)},\\
\eps\ncp^{2}\dscal{g_{2}^{-1}f_2}{\nabla\tvari}&\le 
\norm{g_{2}^{-1}}_{L^{\infty}(\dom)}
\norm{f_{2}}_{L^{1}_{t}L^{2}}\norm{\eps\ncp^{2}\nabla\tvari}_{C^{0}_{t}L^{2}},
\end{align*}
where we used the shorthand notations
$$
\norm{\,\cdot\,}_{C^{0}_{t}X}\defn \sup_{\tau \in [0,t]}\norm{\,\cdot\,}_{X}\quad\mbox{and}\quad
\norm{\,\cdot\,}_{L^{p}_{t}X}\defn \biggl(\int_{0}^{t}\norm{\,\cdot\,}_{X}^{p}\,d\tau\biggr)^{1/p}.
$$
With the bound \eqref{boundLinfty}, 
the previous four estimates and~\eqref{112} imply
\begin{equation}\label{ll}
\begin{aligned}
\norm{\ncp \nabla\tvari}_{L^{2}(\dom)}^{2}
&\le 
C_{0}e^{TC}\la\eps\nu^{2}\dscal{(\partial_{t}+\vard\cdot\nabla)\tvard}{\nabla\tvari}\ra  \\
&\quad +C_{0}e^{TC} 
\Bigl(\norm{\eps\r\ncp\nabla^{2}\tvard}_{L^{2}(\dom)} + 
\norm{f_{2}}_{L^{1}_{t}L^{2}}+\norm{\eps\ncp^{2}\nabla\tvari}_{C^{0}_{t}L^{2}}\Bigr)^{2}.
\end{aligned}
\end{equation}
Note that the second term in the right-hand side of the previous
estimate is bounded by the right-hand side of~\eqref{esti:zL2}. 
It thus remains to estimate the first term. To do so integrate by parts in both the space and time 
variables to obtain
\begin{equation}\label{prop:apriori3.2}
\begin{aligned}
&\eps\nu^{2}\ndscal{(\partial_{t}+\vard\cdot\nabla)\tvard}{\nabla\tvari} \\
&\qquad = \eps\nu^{2}\ndscal{\cn\tvard}{(\partial_{t}+\vard\cdot\nabla)\tvari} \\
&\qquad\quad - \eps\nu^{2}\nscal{\cn\tvard(t)}{\tvari(t)}+\eps\nu^{2}\nscal{\cn\tvard(0)}{\tvari(0)}\\
&\qquad\quad + \eps\nu^{2}\dscal{\nabla\vard}{\nabla\tvari\otimes\tvard}
-\eps\nu^{2}\dscal{\cn\vard}{\tvard\cdot\nabla\tvari}.
\end{aligned}
\end{equation}
The last two terms in the right-hand side of~\eqref{prop:apriori3.2} are estimated by
$$
Kt\norm{\nabla\vard}_{C^{0}_{t}L^{\infty}}^{2}\norm{\tvard}_{C^{0}_{t}L^{2}}^{2}
+\norm{\eps\ncp\nabla\tvari}_{C^{0}_{t}L^{2}}^{2}\le TC w(t)^{2}
+\cpe^{2}Y(t)^{2},
$$
and the second and third terms are estimated by
$$
\norm{\eps\cpe\nabla\tvard}_{C^{0}_{t}L^{2}}^{2}+\norm{\tvari}_{C^{0}_{t}L^{2}}^{2}\le
\cpe^{2} Y(t)^{2}+w(t)^{2}.
$$
In particular, the sum of the last four terms in the right-hand side
of~\eqref{prop:apriori3.2} is estimated by $e^{TC}\{w(t)^{2}+\cpe^{2}Y(t)^{2}\}$. 
This brings us to estimate the first term. 
Using the first equation in~\eqref{system:NSi}, we get
\begin{align*}
&\eps\ncp^{2}\dscal{\cn\tvard}{(\partial_{t}+\vard\cdot\nabla)\tvari}= -\ncp^{2}\ndscal{\cn\tvard}{g_{1}^{-1}\cn \tvard} \\
&\quad+\ncp^{2}\pe\dscal{\cn\tvard}{g_{1}^{-1}\cn(\beta_{1}\nabla\tvariii)}
+ \eps\ncp^{2}\dscal{\cn\tvard}{g_{1}^{-1}f_{1}}.
\end{align*}
Then, the analysis establishing~\eqref{ll} also gives
\begin{align*}
&\eps\ncp^{2}\la\dscal{\cn\tvard}{(\partial_{t}+\vard\cdot\nabla)\tvari}\ra\\
&\quad \le  C_{0}e^{TC}\Bigl\{\ncp^{2}\norm{(\cn\tvard,\pe\nabla^{2}\tvariii)}_{L^{2}(\dom)}^{2}
+\ncp^{2}\pe\norm{(\nabla\tvariii,\eps\nabla\tvard)}_{C^{0}_{t}L^{2}}^{2}
+\norm{f_{1}}_{L^{1}_{t}L^{2}}^{2}\Bigr\},
\end{align*}
which in turn is bounded by the right-hand side
of~\eqref{esti:zL2}. The proof is complete.
%\qed
\end{proof}

We thus have proved a closed system of estimates. Indeed, by taking
appropriate combinations of the previous estimates, we see that there
exists a constant $C_{0}$ depending only on $R_{0}$ and a 
constant $C$ depending only on $R$ such that for all $t\in [0,T]$ and
all $\lambda\ge 1$, the norm $N(t)$ (as defined in~\eqref{defi:N})
satisfies the estimate
\begin{equation*}
N(t)^{2}\le C_{0}e^{TC}N(0)^{2} + \lambda C\int_{0}^{t}N(\tau)^{2}\,d\tau
+\frac{C}{\lambda} N(t)^{2} +C F(T)^{2},
\end{equation*}
where $F(T)\defn
\int_{0}^{t}\norm{(f_{1},f_{2})}_{H^{1}_{\eps\cpe}}+\norm{f_{3}}_{\nhsc{1}}\,dt$. 

Again, take $\lambda$ large enough and apply the
Gronwall's lemma.  The proof of Theorem~\ref{theo:L2} is complete
since $\norm{\tvar}_{{\Hr^{0}_{\pa}(T)}}\le N(T)$.

%%%%%%%%%%%%%%%%%%%%%%%%%%%%%%%%%%%%%%%%%%%%%%%%%%%%%%%%%%%%%%%%%%%%%%%%%%%%%%%%%%%%%%%%%%%%%%%%%%%%%%%%%%%%%%
%%%%%%%%%%%%%%%%%%%%%%%%%%%%%%%%%%%%%%%%%%%%%%%%%%%%%%%%%%%%%%%%%%%%%%%%%%%%%%%%%%%%%%%%%%%%%%%%%%%%%%%%%%%%%%
%%%%%%%%%%%%%%%%%%%%%%%%%%%%%%%%%%%%%%%%%%%%%%%%%%%%%%%%%%%%%%%%%%%%%%%%%%%%%%%%%%%%%%%%%%%%%%%%%%%%%%%%%%%%%%
%%%%%%%%%%%%%%%%%%%%%%%%%%%%%%%%%%%%%%%%%%%%%%%%%%%%%%%%%%%%%%%%%%%%%%%%%%%%%%%%%%%%%%%%%%%%%%%%%%%%%%%%%%%%%%
%%%%%%%%%%%%%%%%%%%%%%%%%%%%%%%%%%%%%%%%%%%%%%%%%%%%%%%%%%%%%%%%%%%%%%%%%%%%%%%%%%%%%%%%%%%%%%%%%%%%%%%%%%%%%%

\section{High frequency regime}\label{section:HFR}
To establish Theorem~\ref{theo:uniform}, 
the crucial part consists in obtaining {\em a priori\/} estimates in Sobolev norms 
independent of~$\pa\in\PA$. 
Theorem~\ref{theo:L2} provides the basic $L^{2}$ estimates. However, the 
estimates of the derivatives also require a careful analysis. Indeed, 
the classical approach, which consists in differentiating the
equations, certainly fails since it reveals unbounded terms in
$\eps^{-1}$. However, one can follow this strategy in the high
frequency regime where the parabolic behavior prevails.

\smallbreak
More precisely, given a smooth solution
$\var=(\vari,\vard,\variii)$ to~\eqref{system:NS2}, we will estimate 
the $\Hr^{s}_{\pa}(T)$ norm of $(\id-J_{\para})\var$, 
where $\Hr^{s}_{\pa}(T)$ is as defined in
Definition~\ref{defi:decoupling} and, in keeping with the notations of \S\ref{Nonlinear}, $\{ \,J_{\para}\,\arrowvert\,\para\in [0,1]\,\}$ is a
Friedrichs' mollifier and $\Fi{}{s}\defn(\id-\Delta)^{s/2}$. 
In order to make our energy estimates applicable, the main difficulty is to verify that the commutator of
$(\id-J_{\para})\Fi{}{s}$ and the equations~\eqref{system:NS2} can
be seen as a source term. To do so we first note that for
$\para=\Or(\eps)$, one can gain an extra factor $\eps$ in the
commutator estimates [see \eqref{commutator:estimate2} below]. Yet, this costs a derivative. To compensate this loss of
derivative, we use in an essential way 
the parabolic behavior of the equations. Consequently, we search $\para$
under the form $c(\r,\pe)\eps$.  
Since for our purposes the main smoothing effect concerns the
penalized terms $\cn\vard$ and $\nabla\vari$, we take $c(\r,\pe)=\sqrt{\r+\pe}$.

\begin{proposition}\label{prop:HFuniform}
Given an integer~$s>1+d/2$, there exists a continuous non-decreasing 
function $C$ such that for all 
$\pa=(\eps,\r,\pe)\in\PA$, all $T\in [0,1]$ and all 
$U=(\vari,\vard,\variii)\in C^{1}([0,T];H^{\infty}(\xD))$ 
satisfying~\eqref{system:NS2},
\begin{equation}\label{prop:slowcomponent:1}
\norm{(\id-J_{\eps\cpe})\var}_{\Hr_{\pa}^{s}(T)}\le C(\Omega_{0})e^{\sqrt{T}C(\Omega)},
\end{equation}
where $\cpe\defn\sqrt{\r+\pe}$,
$\Omega_{0}\defn\norm{\var(0)}_{\Hr_{\pa,0}^{s}}$, 
$\Omega\defn \norm{\var}_{\Hr_{\pa}^{s}(T)}$  (see Definition~$\ref{defi:decoupling}$).
\end{proposition}
\begin{remark}
Note that we establish estimates for the exact solutions of~\eqref{system:NS2} and we do not estimate the solutions of the
linearized system~\eqref{system:NSi}.
\end{remark}

%We make two running conventions: In all the statements, we will focus on estimates: 
%it is always assumed that the functions are smooth enough so that 
%the norms involved in the estimates are well-defined. 
%Given a normed space $X$ and function $u\colon [0,T]\rightarrow X$, we 
%denote by $\norm{u}_{X}$ the function 
%$[0,T]\ni t \mapsto \norm{u(t)}_{X}$. In particular, we do not make further reference to 
%the time variable $t$ (as it is customary). Yet it is always checked that 
%the estimates are independent of the parameter $T$.

\subsection{Preliminaries} 
To avoid interruptions of the proofs later on, we now collect a few
nonlinear estimates we use throughout this section. 

Recall that $\norm{\cdot}_{H^{\sigma}_{\varrho}}\defn\norm{\cdot}_{H^{\sigma-1}}+\varrho\norm{\cdot}_{\nh{\sigma}}$.

\begin{lemma}
Let $s>1+d/2$. 
There exists a constant $K$ such that for all $\para\in [0,1]$, all
$\varrho\ge 0$, all 
$f\in H^{s+1}(\xD)$ and all $u\in H^{s}(\xD)$, one has
\begin{equation}\label{commutator:estimate2}
\norm{\bigl[ f,(\id-J_{\para})\Fi{}{s}\bigr]u}_{H^{1}_{\varrho}} 
\le (\varrho+\para) K \bigl\{\norm{f}_{\Lip}\norm{u}_{\nh{s}}+\norm{f}_{\nh{s+1}} \norm{u}_{L^{\infty}}\bigr\},
\end{equation}
where $\norm{f}_{\Lip}\defn\norm{f}_{L^{\infty}}+\norm{\nabla f}_{L^{\infty}}$.
\end{lemma}
\begin{remark}
This inequality provides a way to gain an extra factor $\para$ [with
$\varrho=\para$] or a derivative [with $\varrho=1$]. In this respect
it is like the estimate~\eqref{Friedrichs}, to which we will return
in a moment (see Lemma~\ref{bignorms:commutator}).
\end{remark}

\begin{proof}
To prove this result we need a tame estimate version
of~\eqref{commutator:estimateusual}. We use the following result,
which is Proposition~$3.6.\mbox{A}$ in~\cite{Taylor} (see also~\cite{David}). 
Let $m>0$. For all Fourier multiplier $\Qr$ with symbol $q\in S^{m}$ and for all $\sigma\ge 0$, we have
\begin{equation}\label{commutator:estimate}
\lA \bigl[f,\Qr\bigr]u \rA_{H^{\indexg}}\le K\lA f\rA_{\Lip}\norm{u}_{\nh{\indexg+m-1}}+K\norm{f}_{\nh{m+\indexg}}\norm{u}_{L^{\infty}},
\end{equation}
where the constant $K$ depends only on $m$, $s$, $\sigma$, $d$ and a finite number of semi-norms of~$q$ in~$S^{m}$. 

Since the support of the symbol~$1-\jmath_{\para}$ is limited by $|\xi|\ge 1/\para$, 
the most direct estimates show that $\{\,\para^{-1}(1-\jmath(\para\xi))\L{\xi}^{s}\,\arrowvert \para\in [0,1]\,\}$ 
is uniformly bounded in the symbol class $S^{s+1}$. Thus, the commutator 
estimate~\eqref{commutator:estimate}, applied with $(m,\sigma)\defn (s+1,0)$, implies
\begin{equation}\label{commutator:estimate.HL2}
\norm{\bigl[ f,(\id-J_{\para})\Fi{}{s}\bigr]u}_{\nhz} 
\les \para \bigl(\norm{f}_{\Lip}\norm{u}_{\nh{s}}+\norm{f}_{\nh{s+1}} \norm{u}_{L^{\infty}}\bigr).
\end{equation}
On the other hand, the family $\{\,(1-\jmath(\para\xi))\L{\xi}^{s}\,\arrowvert \para\in [0,1]\,\}$ 
is uniformly bounded in $S^{s}$, so that~\eqref{commutator:estimate}
applied with $(m,\sigma)=(s,1)$ implies
\begin{equation}\label{commutator:estimate.HL3}
\norm{\bigl[ f,(\id-J_{\para})\Fi{}{s}\bigr]u}_{\nh{1}} 
\les\norm{f}_{\Lip}\norm{u}_{\nh{s}}+\norm{f}_{\nh{s+1}}\norm{u}_{L^{\infty}}.
\end{equation}
Combining~\eqref{commutator:estimate.HL2} with~\eqref{commutator:estimate.HL3} multiplied by $\varrho$, yields~\eqref{commutator:estimate2}.
%\qed
\end{proof}

We next prove two Moser-type estimates for the norms~$\norm{\cdot}_{H_{\varrho}^{\indexg}}$. 
\begin{lemma}\label{lemm:weightedprod}
Let $s>1+d/2$. There exists a constant $K$ such that for all
$\varrho\ge 0$ and for all $u_{i}\in H^{s}_{\varrho}(\xD)$,
\begin{equation}
\norm{u_{1}u_{2}}_{H_{\varrho}^{s}}\le K \norm{u_{1}}_{H_{\varrho}^{s}}\norm{u_{2}}_{H_{\varrho}^{s}}.\label{esti:weightedprod}
\end{equation}
This result extends to vector valued functions. 
\end{lemma}
\begin{proof}
Using the standard tame estimate for products at orders 
$s-1$ and~$s$, we get
\begin{align*}
\norm{u_{1}u_{2}}_{\nh{s-1}}&\les 
\norm{u_{1}}_{L^{\infty}}\norm{u_{2}}_{\nh{s-1}}
+\norm{u_{1}}_{\nh{s-1}}\norm{u_{2}}_{L^{\infty}},\\
\norm{u_{1}u_{2}}_{\nh{s}}&\les 
\norm{u_{1}}_{L^{\infty}}\norm{u_{2}}_{\nh{s}}
+\norm{u_{1}}_{\nh{s}}\norm{u_{2}}_{L^{\infty}}.
\end{align*}
Using the definition $\norm{\cdot}_{H_{\varrho}^{s}}\defn \norm{\cdot}_{\nh{s-1}}
+\varrho\norm{\cdot}_{\nh{s}}$, and putting 
the parameter~$\varrho$ in appropriate spots yields
$$
\norm{u_{1}u_{2}}_{H_{\varrho}^{s}}\les\norm{u_{1}}_{L^{\infty}}
\norm{u_{2}}_{H_{\varrho}^{s}}+\norm{u_{1}}_{H_{\varrho}^{s}}\norm{u_{2}}_{L^{\infty}}.
$$
Since $s-1>d/2$, the Sobolev Theorem implies 
$\norm{u}_{L^{\infty}}\les \norm{u}_{\nh{s-1}}\les \norm{u}_{H_{\varrho}^{s}}$.
Which completes the proof.
%\qed
\end{proof}

\begin{lemma}\label{lemm:weightedcomp}
Let~$s>1+d/2$ and~$F\colon \xR^{n} \rightarrow \xC$ be a $C^{\infty}$
function such that~$F(0)=0$. 
Then, for all $\varrho\ge 0$ and 
all~$u\in H^{s}_{\varrho}(\xD)$ with values in $\xR^{n}$,
\begin{equation} \label{esti:weightedcomp}
\norm{F(u)}_{H_{\varrho}^{s}}\le C\bigl(\norm{u}_{H_{\varrho}^{s}}\bigr),
\end{equation}
where $C(\cdot)$ depends only on a finite number of
semi-norms of $F$ in $C^{\infty}$.
\end{lemma}
\begin{proof}
Recall that, for all $\index> d/2$ and all $C^{\infty}$ function $F$ satisfying $F(0)=0$,
\begin{equation} \label{comp}
\lA F(u)\rA_\nh{\index}\leqslant C(\lA u\rA_{L^{\infty}})\norm{ u}_{\nh{\index}}.
\end{equation}
Using this estimate at orders $s-1$ and $s$, one has 
\begin{equation}\label{esti:weightcomppre}
\norm{F(u)}_{H_{\varrho}^{s}}\le C\bigl(\norm{u}_{L^{\infty}}\bigr)\norm{u}_{H_{\varrho}^{s}}.
\end{equation}
This in turn implies the desired estimate.
%\qed
\end{proof}

\begin{lemma}
Let $s>1+d/2$,~$F\colon \xR^{n} \rightarrow \xC$ be a $C^{\infty}$
function such that~$F(0)=0$ and $\Qr$ be a Fourier multiplier with symbol $q\in S^{s}$. 
For all vector-valued function $u\in H^{s}(\xD)$ one has
\begin{equation}\label{paralinearization}
\norm{\Qr\bigl(F(u)\bigr)-F'(u)\Qr u}_{\nh{1}}\le C\bigl(\norm{u}_{\nh{s}}\bigr),
\end{equation}
where $F'$ is the differential of $F$ and $C(\cdot)$ is a smooth non-decreasing function depending only on $s$, $d$, a finite number of semi-norms of $q$ in $S^{s}$ 
and a finite number of semi-norms of $F$ in $C^{\infty}$. 
\end{lemma}
\begin{proof}
To establish~\eqref{paralinearization} we use the para-differential calculus of Bony~\cite{Bony}. 
Denote by~$T_f$ the operator of para-multiplication by~$f$. Starting from
$$
F(u) = T_{F'(u)}u+R(u;x,D_{x})u,
$$
where $R(u;x,D_{x})$ is a smoothing operator (see~\eqref{Com_Paradif_32}), we obtain
$$
\Qr \bigl(F(u)\bigr)-F'(u)\Qr u = \bigl(T_{F'(u)}-F'(u)\bigr)\Qr
u+\bigl[\Qr,T_{F'(u)}\bigr] u + \Qr R(u;x,D_{x})u.
$$
The claim then follows from the
bounds~\eqref{Com_Paradif_1}--\eqref{Com_Paradif_3} and the estimate
\begin{equation}
\lA R(u;x,D_{x})u\rA_{\nh{s+1}}\le C\bigl(\norm{u}_{\nh{s}}\bigr).
\label{Com_Paradif_32}
\end{equation}
See~\cite[Th.~$10.3.1$]{HorL} and~\cite[Corollary~$9.3.6$]{HorL} for the proof of~\eqref{Com_Paradif_32}. 
%\qed
\end{proof}

\subsection{Localization in the high frequency region}
To simplify the presentation, we fix a real number~$s$ strictly
greater than~$1+d/2$. 

To proceed further, we need some more terminology. 
\begin{nota} For all~$\para\in [0,2]$, define
$$
\Qr_{\para}\defn (\id-J_{\para})\Fi{}{s}.
$$
Hereafter, the parameter $\para$ is a product $\eps\cpe$ with
$(\eps,\cpe)\in [0,1]\times [0,2]$.
\end{nota}

Introduce next the commutator of the equations~\eqref{system:NS2}
and $\Qr_{\eps\cpe}$.
\begin{nota}\label{nota:sourceHF}Given~$\pa=(\eps,\r,\pe)\in\PA$, $\cpe\in [0,2]$ and 
$\var=(\vari,\vard,\variii)$, set
\begin{alignat*}{4}
f^{\pa,\cpe}_{1,\HF}(\var)&\defn\bigl[g_{1}(\phi),\Qr_{\eps\cpe}\bigr]\partial_{t}\vari
&&+\bigl[g_{1}(\phi)\vard,\Qr_{\eps\cpe}\bigr]\cdot\nabla\vari  
&&- \frac{\pe}{\eps}&&\bigl[B_{1}(\phi),\Qr_{\eps\cpe}\bigr]\variii,\\
f^{\pa,\cpe}_{2,\HF}(\var)&\defn\bigl[g_{2}(\phi),\Qr_{\eps\cpe}\bigr]\partial_{t}\vard
&&+\bigl[g_{2}(\phi)\vard,\Qr_{\eps\cpe}\bigr]\cdot\nabla\vard
&&-\r&&\bigl[B_{2}(\phi),\Qr_{\eps\cpe}\bigr]\vard,\\[0.5ex]
f^{\pa,\cpe}_{3,\HF}(\var)&\defn \bigl[g_{3}(\phi),\Qr_{\eps\cpe}\bigr]\partial_{t}\variii 
&&+\bigl\{ g_{3}(\phi)\vard\cdot\nabla\variii;\Qr_{\eps\cpe} \bigr\}
&&-\pe&&\bigl[B_{3}(\phi),\Qr_{\eps\cpe}\bigr]\variii,
\end{alignat*}
where $\phi$ is a shorthand notation for $(\variii,\eps\vari)$ and 
$$
\bigl\{ g_{3}(\phi)\vard\cdot\nabla\variii;\Qr_{\eps\cpe} \bigr\}\defn 
g_{3}(\phi)\vard\cdot\nabla\Qr_{\eps\cpe}\variii+g_{3}(\phi)\nabla\variii\cdot\Qr_{\eps\cpe}\vard
-\Qr_{\eps\cpe}\bigl(g_{3}(\phi)\vard\cdot\nabla\variii\bigr).
$$
\end{nota}
\begin{remark}
For the purpose of proving estimates independent of $\pe$, we do not consider the exact commutator of the
third equation in~\eqref{system:NS2} with~$\Qr_{\eps\cpe}$.
\end{remark}

Let us recall that
\begin{align*}
B_{1}(\phi)&\defn\chi_{1}(\eps\vari)\cn(\beta(\variii)\nabla\,\cdot\,), \\[0.5ex]
B_{2}(\phi)&\defn\chi_{2}(\eps\vari)\cn(2\zeta(\variii)D\cdot)
+\chi_{2}(\eps\vari)\nabla(\eta(\variii)\cn\cdot),\\[0.5ex]
B_{3}(\phi)&\defn\chi_{3}(\eps\vari)\cn(\beta(\variii)\nabla\cdot).
\end{align*}

The following lemma shows that $f^{\pa,\cpe}_{\HF}$ can be seen
as a source term.

\begin{lemma}\label{lemm:commHF}
There exists a smooth non-decreasing function $C=C(\cdot)$ such that for all $\pa\in\PA$, 
all $T\ge 0$, all $\cpe\in [0,2]$ and all vector valued function 
$\var\defn(\vari,\vard,\variii)\in C^{1}([0,T];H^{\infty}(\xD))$,
\begin{alignat*}{5}
& \bigl\lVert{f^{\pa,\cpe}_{1,\HF}(\var)\bigr\rVert}_{H^{1}_{\eps\cpe}} &&\le &&C(R) 
\bigr\{&&1+ \norm{\eps\partial_{t}\vari}_{H^{s}_{\cpe}}
&&+\pe\norm{\variii}_{\nh{s+2}}\bigr\},\\
&\bigl\lVert{f^{\pa,\cpe}_{2,\HF}(\var)\bigr\rVert}_{H^{1}_{\eps\cpe}} &&\le &&C (R)
\bigr\{&&1+\norm{\eps\partial_{t}\vard}_{H^{s}_{\cpe}}
&&+\r\norm{\eps\vard}_{\nh{s+2}}\bigr\},\\
&\bigl\lVert{f^{\pa,\cpe}_{3,\HF}(\var)\bigr\rVert}_{\nhsc{1}}
&&\le &&C (R)\bigr\{&&1+\norm{\partial_{t}\variii}_{H^{s}_{\cpe}}&&+\pe\norm{\variii}_{\nh{s+2}}\bigr\},
\end{alignat*}
where $R\defn \norm{(\vari,\vard)}_{H^{s+1}_{\eps\cpe}}+\norm{\variii}_{\nhsc{s+1}}$.
\end{lemma}
\begin{remark}
We will apply this lemma with $\cpe\defn\sqrt{\r+\pe}$. Yet, we prove
estimates independent of $\cpe\in [0,2]$ to explain why the frequency
space is cut around $\la\xi\ra \approx 1/(\eps\sqrt{\r+\pe})$.
\end{remark}
\begin{proof}
We make intensive use of the following obvious observations. 
Firstly, using  the Sobolev embedding Theorem and the very
definition of the norms $\norm{\cdot}_{H^{\sigma}_{\varrho}}$, one has
\begin{equation}\label{esti:HFobservation}
\norm{u}_{L^{\infty}}\les \norm{u}_{\nh{s-1}}\le \norm{u}_{\nhsc{s}}
\quad\mbox{and}\quad \norm{u}_{\Lip}\les \norm{u}_{\nh{s}}\le \norm{u}_{\nhsc{s+1}}.
\end{equation}
Since $\eps\le 1$ and $\cpe\le 2$, directly from the definition of
$\norm{\cdot}_{\nhsc{\sigma}}$, one has
\begin{equation}\label{esti:HFobservation2}
\cpe\norm{u}_{\nh{\sigma}}\le \norm{u}_{\nhsc{\sigma}}\le 3
\norm{u}_{\nh{\sigma}}\quad\mbox{and}\quad 
\norm{\eps u}_{\nhsc{s+1}}\le \norm{u}_{H^{s+1}_{\eps\cpe}}.
\end{equation}
Note the following corollary of the second inequality:
$
\norm{(\eps\vari,\eps\vard)}_{\nhsc{s+1}}\le R.
$
This in turn implies
$\norm{\phi}_{\nhsc{s+1}}\defn\norm{(\variii,\eps\vari)}_{\nhsc{s+1}}\le R$.
\smallbreak

\noindent
{\sc{Step 1:}} Estimate for $f^{\pa,\cpe}_{1,\HF}$(\var). {\textbf{a)}} We begin by proving that
\begin{equation}\label{commHF:01}
\norm{\bigl[g_{1}(\phi),\Qr_{\eps\cpe}\bigr]\partial_{t}\vari}_{H^{1}_{\eps\cpe}}\le
C(R) \norm{\eps\partial_{t}\vari}_{H^{s}_{\cpe}}.
\end{equation}
Applying the commutator estimate~\eqref{commutator:estimate2}
with~$\para=\varrho=\eps\cpe$, we have
\begin{align*}%\label{commHF:019}
\norm{\bigl[g_{1}(\phi),\Qr_{\eps\cpe}\bigr]\partial_{t}\vari}_{H^{1}_{\eps\cpe}}
&\les
\eps\cpe\norm{\tilde{g}_{1}(\phi)}_{\Lip}\norm{\partial_{t}\vari}_{\nh{s}}
+\eps\cpe\norm{\tilde{g}_{1}(\phi)}_{\nh{s+1}}\norm{\partial_{t}\vari}_{L^{\infty}}\\
&\les
\norm{\tilde{g}_{1}(\phi)}_{\Lip}\norm{\eps\partial_{t}\vari}_{\nhsc{s}}
+\norm{\tilde{g}_{1}(\phi)}_{\nhsc{s+1}}\norm{\eps\partial_{t}\vari}_{L^{\infty}},
\end{align*}
where $\tilde{g}_{1}$ is defined by $\tilde{g}_{1}=g_{1}-g_{1}(0)$. 
The estimates~\eqref{esti:HFobservation} imply that the right side is
less than
$K\norm{\tilde{g}_{1}(\phi)}_{\nhsc{s+1}}\norm{\eps\partial_{t}\vari}_{\nhsc{s}}$. 
Using Lemma~\ref{lemm:weightedcomp}, we obtain $\norm{\tilde{g}_{1}(\phi)}_{\nhsc{s+1}}\le
C(\norm{\phi}_{\nhsc{s+1}})\le C(R)$. 
This proves~\eqref{commHF:01}.

\smallbreak
\noindent {\textbf{b)}} Next, we prove that
\begin{equation}\label{commHF:02}
\norm{\bigl[g_{1}(\phi)\vard,\Qr_{\eps\cpe}\bigr]\cdot\nabla\vari}_{H^{1}_{\eps\cpe}}\le
C(R) .
\end{equation}
Again, this follows from the commutator
estimate~\eqref{commutator:estimate2} applied with $\para=\varrho=\eps\cpe$. 
Indeed, it yields
\begin{equation*}
\norm{\bigl[g_{1}(\phi)\vard,\Qr_{\eps\cpe}\bigr]\cdot\nabla\vari}_{H^{1}_{\eps\cpe}}\les 
\norm{\tilde{g}_{1}(\phi)\vard}_{\Lip}\norm{\eps\nabla\vari}_{\nhsc{s}}
+\norm{\eps\tilde{g}_{1}(\phi)\vard}_{\nhsc{s+1}}\norm{\nabla\vari}_{L^{\infty}}.
\end{equation*}
The first term in the right-hand side is estimated by $C(R)$ since
$$
\norm{\tilde{g}_{1}(\phi)\vard}_{\Lip}\le 
C\bigl(\norm{\phi}_{\Lip},\norm{\vard}_{\Lip}\bigr)\le
C(\norm{\phi}_{\nh{s}},\norm{\vard}_{\nh{s}}) \le C(R),
$$
and $\norm{\eps\nabla\vari}_{\nhsc{s}}\le
\norm{\vari}_{H^{s+1}_{\eps\cpe}}\le R$. 

It thus remains to estimate $\norm{\eps\tilde{g}_{1}(\phi)\vard}_{\nhsc{s+1}}$. 
To do so use Lemma~\ref{lemm:weightedprod} and~\ref{lemm:weightedcomp}, to obtain  
$$
\norm{\eps\tilde{g}_{1}(\phi)\vard}_{\nhsc{s+1}}\les
\norm{\tilde{g}_{1}(\phi)}_{\nhsc{s+1}}\norm{\eps\vard}_{\nhsc{s+1}}\le C(R)R.
$$

\smallbreak
\noindent {\textbf{c)}} Let us prove that
\begin{equation}\label{commHF:03}
\frac{\pe}{\eps}\norm{\bigl[B_{1}(\phi),\Qr_{\eps\cpe}\bigr]\variii}_{H^{1}_{\eps\cpe}}\le
C(R)\norm{\pe\variii}_{\nh{s+2}}.
\end{equation}
By definition $B_{1}(\phi)=\chi_{1}(\eps\vari)\cn(\beta(\variii)\nabla\,\cdot\,)$, so that one 
can decompose the commutator 
$\bigl[B_{1}(\phi),\Qr_{\eps\cpe}\bigr]\variii$ as
\begin{equation}\label{split:HFB1kl}
\bigl[\Er_{1}(\phi),\Qr_{\eps\cpe}\bigr]\Delta\variii
+\bigl[\Er_{2}(\phi,\nabla\phi),\Qr_{\eps\cpe}\bigr]\cdot\nabla\variii,
\end{equation}
where $\Er_{1}(\phi)=\chi_{1}(\eps\vari)\beta(\variii)$ and $\Er_{2}(\phi,\nabla\phi)=
\chi_{1}(\eps\vari)\nabla\beta(\variii)=\chi_{1}(\eps\vari)\beta'(\variii)\nabla\variii$. 
The commutator estimate~\eqref{commutator:estimate2} [applied with
$\para=\varrho=\eps\cpe$] implies
\begin{align*}
\norm{\bigl[\Er_{1},\Qr_{\eps\cpe}\bigr]\Delta\variii}_{H^{1}_{\eps\cpe}}&\les\eps
\norm{\tilde{\Er}_{1}}_{\Lip}\norm{\Delta\variii}_{\nhsc{s}}
+\eps\norm{\tilde{\Er}_{1}}_{\nhsc{s+1}}\norm{\Delta\variii}_{L^{\infty}},\\
\norm{\bigl[\Er_{2},\Qr_{\eps\cpe}\bigr]\cdot\nabla\variii}_{H^{1}_{\eps\cpe}}&\les\eps
\norm{\cpe\Er_{2}}_{\Lip}\norm{\nabla\variii}_{\nh{s}}
+\eps\norm{\Er_{2}}_{\nhsc{s+1}}\norm{\nabla\variii}_{L^{\infty}},
\end{align*}
where $\tilde{\Er}_{1}\defn\Er_{1}-\Er_{1}(0)$. We gather
\begin{align*}
&\norm{\tilde{\Er}_{1}(\phi)}_{\Lip}\le C(\norm{\phi}_{\nh{s}})\le C(R),
&&\norm{\Delta\variii}_{\nhsc{s}}\les \norm{\variii}_{\nh{s+2}},\\
&\norm{\tilde{\Er}_{1}(\phi)}_{\nhsc{s+1}}\le C(\norm{\phi}_{\nhsc{s+1}})\le C(R),
&&\norm{\Delta\variii}_{L^{\infty}}\les \norm{\variii}_{\nh{s+2}},\\
&\norm{\cpe\Er_{2}(\phi,\nabla\phi)}_{\Lip}\le C(\norm{\phi}_{\nhsc{s+1}})\le C(R),
&&\norm{\nabla\variii}_{\nh{s}}\le \norm{\variii}_{\nh{s+2}},\\
&\norm{\Er_{2}(\phi,\nabla\phi)}_{\nhsc{s+1}}\le C(R)\norm{\variii}_{\nh{s+2}},
&&\norm{\nabla\variii}_{L^{\infty}}\les \norm{\variii}_{\nh{s}}\le R.
\end{align*}
The first six inequalities follows from Lemma~\ref{lemm:weightedcomp}, 
the Sobolev Theorem and the estimates \eqref{esti:HFobservation}--\eqref{esti:HFobservation2}. 
To estimate $\norm{\tilde{\Er}_{2}(\phi,\nabla\phi)}_{\nhsc{s+1}}$,  
we first use Lemma~\ref{lemm:weightedprod}, to obtain
\begin{equation}\label{need:beta1variii}
\norm{\Er_{2}(\phi,\nabla\phi)}_{\nhsc{s+1}}\les
(1+\norm{\tilde{\chi_{1}}(\eps\vari)}_{\nhsc{s+1}})(1+\norm{\tilde{\beta'}(\variii)}_{\nhsc{s+1}})\norm{\nabla\variii}_{\nhsc{s+1}},
\end{equation}
and next use Lemma~\ref{lemm:weightedcomp} to bound the
$H^{s+1}_{\cpe}$-norms of $\tilde{\chi_{1}}(\eps\variii)$ and
$\tilde{\beta'}(\variii)$ by $C(R)$. This yields the desired bound since
$\norm{\nabla\variii}_{\nhsc{s+1}}
\les\norm{\variii}_{\nh{s+2}}$.

The claim~\eqref{commHF:03} then easily follows from these estimates. 

\smallbreak
\noindent{\sc{Step 2:}}{ Estimate for $f^{\pa,\cpe}_{2,\HF}(\var)$}. 
Note that one can obtain $f^{\pa,\cpe}_{2,\HF}(\var)$ from
$f^{\pa,\cpe}_{1,\HF}(\var)$ by replacing $\vari$ by $\vard$ and $\variii$ by
$\eps\vard$. Therefore, we are back in the situation of the previous step, and hence
conclude that
\begin{align*}
\norm{\bigl[g_{2}(\phi),\Qr_{\eps\cpe}\bigr]\partial_{t}\vard}_{H^{1}_{\eps\cpe}}&\le
C(R) \norm{\eps\partial_{t}\vard}_{H^{s}_{\cpe}},\\
\norm{\bigl[g_{2}(\phi)\vard,\Qr_{\eps\cpe}\bigr]\cdot\nabla\vard}_{H^{1}_{\eps\cpe}}&\le
C(R),\\
\r\norm{\bigl[B_{2}(\phi),\Qr_{\eps\cpe}\bigr]\vard}_{H^{1}_{\eps\cpe}}&\le
C(R)\norm{\r\eps\vard}_{\nh{s+2}}.
\end{align*}

\smallbreak
\noindent{\sc{Step 3:}} Estimate for $f^{\pa,\cpe}_{3,\HF}(\var)$. 
The estimates for the first and the last terms in
$f^{\pa,\cpe}_{3,\HF}(\var)$ can be deduced by following the previous
analysis. We only indicate the point at which the argument differs: 
we use the commutator estimate~\eqref{commutator:estimate.HL3} 
with $\para=\eps\cpe$ and $\varrho=\cpe$ (instead of $\varrho=\eps\cpe$). 

Let us concentrate on the second term. We claim that
\begin{equation*}%\label{esti:goodlinea}
\norm{\bigl\{ g_{3}(\phi)\vard\cdot\nabla\variii;\Qr_{\eps\cpe}
  \bigr\}}_{\nhsc{1}}\le C(R).
\end{equation*}
We first decompose $\bigl\{
g_{3}(\phi)\vard\cdot\nabla\variii;\Qr_{\eps\cpe} \bigr\}$ as 
\begin{equation}\label{defi:ZHF}
\bigl[g_{3}(\phi),\Qr_{\eps\cpe}\bigr]\vard\cdot\nabla\variii 
+g_{3}(\phi) Z
\end{equation}
where $Z\defn 
\vard\cdot\nabla\Qr_{\eps\cpe}\variii
+\nabla\variii\cdot\Qr_{\eps\cpe}\vard
-\Qr_{\eps\cpe}(\vard\cdot\nabla\variii)$.

Applying the commutator estimate~\eqref{commutator:estimate2}
with~$\para=\eps\cpe$ and $\varrho=\cpe$, and using the
estimates~\eqref{esti:HFobservation}, we get
$$
\norm{\bigl[g_{3}(\phi),\Qr_{\eps\cpe}\bigr]\vard\cdot\nabla\variii}_{\nhsc{1}}
\les \norm{\tilde{g}_{3}(\phi)}_{\nhsc{s+1}}\norm{\vard\cdot\nabla\variii}_{\nhsc{s}}.
$$
Using Lemma~\ref{lemm:weightedprod} and~\ref{lemm:weightedcomp},
we find that the right-hand side of the previous estimate is dominated
by $C(\norm{\phi}_{\nhsc{s+1}})\norm{\vard}_{\nhsc{s}}\norm{\nabla\variii}_{\nhsc{s}}\le
C(R)$.

\smallbreak
Since $\norm{g_{3}(\phi) Z}_{\nhsc{1}}\le
\norm{g_{3}(\phi)}_{W^{1,\infty}}\norm{Z}_{\nhsc{1}}\le
C(R)\norm{Z}_{\nhsc{1}}$, 
to control the $H^{1}_{\cpe}$-norm of the second term 
in~\eqref{defi:ZHF}, only the estimate
of~$\norm{Z}_{\nhsc{1}}$ is missing. We split $\norm{Z}_{\nhsc{1}}$ 
as $\norm{Z}_{\nhz}+\cpe\norm{Z}_{\nh{1}}$, and we prove that 
$\norm{Z}_{\nhz}$ and $\cpe\norm{Z}_{\nh{1}}$ are both estimated by
$C(R)$. To estimate $\norm{Z}_{\nhz}$, note that
\begin{equation}\label{split:ZHFL2}
Z=\nabla\variii\cdot\Qr_{\eps\cpe}\vard+ [\vard,\Qr_{\eps\cpe}]\cdot\nabla\variii.
\end{equation}
The second term in the right-hand side is estimated by way of the 
commutator estimate~\eqref{commutator:estimateusual}. Indeed, 
applying Lemma~\ref{commutateur_S^1_f} with $\sigma_{0}=m=s$ 
and $\sigma=s-1$, we get
$$
\norm{[\vard,\Qr_{\eps\cpe}]\cdot\nabla\variii}_{\nhz}\les
\norm{\vard}_{\nh{s}}\norm{\nabla\variii}_{\nh{s-1}}\le 
\norm{\vard}_{\nh{s}}\norm{\variii}_{\nh{s}}\le R^{2}.
$$
As regards the first term in the right-hand side of~\eqref{split:ZHFL2}, write
$$
\norm{\nabla\variii\cdot\Qr_{\eps\cpe}\vard}_{\nhz}\le 
\norm{\nabla\variii}_{L^{\infty}}\norm{\Qr_{\eps\cpe}\vard}_{\nhz}
\les \norm{\variii}_{\nh{s}}\norm{\vard}_{\nh{s}}\le R^{2}.
$$

Moving to the estimate of $\cpe\norm{Z}_{\nh{1}}$, remark that
$$
\cpe Z= F'(u)\Qr_{\eps\cpe}u-\Qr_{\eps\cpe}(F(u)),
$$
with $u=(u_{1},u_{2})\defn(\vard,\cpe\nabla\variii)$ and
$F(u)=u_{1}u_{2}$. Hence, \eqref{paralinearization}
yields 
$$
\cpe\norm{Z}_{\nh{1}}\le
C\bigl(\norm{(\vard,\cpe\nabla\variii)}_{\nh{s}}\bigr)\le
C\bigl(\norm{\vard}_{\nh{s}}+\norm{\variii}_{\nhsc{s+1}}\bigr)\le C(R).
$$
The previous estimates imply that 
$\norm{Z}_{\nhsc{1}}$ is
controlled by $C(R)$. 
This completes the proof of Lemma~\ref{lemm:commHF}.
%\qed
\end{proof}
\begin{remark}
Let us explain the reason why we assume that $\beta$ depends only on $\variii$. 
Had we worked instead with general coefficient $\beta$ depending 
also on $\eps\vari$, 
the corresponding inequality~\eqref{need:beta1variii} would 
have involved $\norm{\eps\nabla\vari}_{\nhsc{s+1}}$. The problem
presents itself: $\norm{\eps\vari}_{L^{1}(0,T;H^{s+2}_{\cpe})}$ is not 
controlled by the
norm~$\norm{(\vari,\vard,\variii)}_{\Hr^{s}_{\pa}(T)}$. It is possible
to get around
the previous problem, yet we do not address this question.
\end{remark}

In view of Lemma~\ref{lemm:commHF} we are led to estimate 
$\partial_{t}\ter$ where $\ter\defn (\variii,\eps\vari,\eps\vard)$. 

\begin{lemma}\label{lemm:estdescoef123}
There exists a continuous non-decreasing 
function $C(\cdot)$ such that for all 
$\pa\in\PA$, all $T\ge 0$, all $\cpe\in [0,2]$ and all 
$(\vari,\vard,\variii)\in \mathcal{H}^{s}_{\pa}(T)$ 
solving~\eqref{system:NS2}, the function 
$\ter\defn (\variii,\eps\vari,\eps\vard)$ satisfies
\begin{equation}\label{mlemm:estdescoef}
\norm{\partial_{t}\ter}_{\nhsc{s}}\le
C(R)\bigl\{1+R'\bigr\},
\end{equation}
where
\begin{equation}\label{defi:RR'}
\begin{split}
R&\defn\norm{(\vari,\vard)}_{H^{s+1}_{\eps{\cpe}}}
+\norm{\variii}_{H^{s+1}_{\cpe}},\\
R'&\defn \cpe\norm{(\cn\vard,\nabla\vari)}_{\nh{s}}
+\sqrt{\r}\norm{\nabla\vard}_{H^{s+1}_{\eps{\cpe}}}
+\sqrt{\pe}\norm{\nabla\variii}_{H^{s+1}_{{\cpe}}}.
\end{split}
\end{equation}
\end{lemma} 
\begin{remark}
This estimate plays a key role for the purpose 
of proving estimates independent of $\r$ and $\pe$. 
Indeed, this is the only step in which we use 
the additional smoothing effect for $\cn\vard$ and $\nabla\vari$. 
More precisely, the fact that the estimate~\eqref{mlemm:estdescoef}
is tame [linear in $R'$] allows us to control 
$\norm{\partial_{t}\ter}_{L^{2}(0,T;\nhsc{s})}$ 
by $C\bigl(\norm{(\vari,\vard,\variii)}_{\Hr^{s}_{\pa}(T)}\bigr)$ for
all $\cpe \le \sqrt{\r+\pe}$.
\end{remark}
\begin{proof}
In order to establish~\eqref{mlemm:estdescoef}, observe that 
\begin{equation}\label{Eq:terest}
G(\phi)\bigl(\partial_{t}\psi+\vard\cdot\nabla\psi)+L\var +B_{\r,\pe}(\phi)\psi=0,
\end{equation}
where $\psi$ and $\var$ are identified with
$(\eps\vari,\eps\vard,\variii)^{t}$ and $(\vari,\vard,\variii)^{t}$,
respectively, and 
\begin{equation*}
G\defn\begin{pmatrix} g_{1}& 0 &0 \\ 
0& g_{2}I_{d}&0 \\
0&0 & g_{3}
\end{pmatrix},~
L\defn \begin{pmatrix}
0 & \cn & 0 \\
\nabla & 0 & 0 \\
0 & \cn & 0
\end{pmatrix},~
B_{\r,\pe}\defn 
\begin{pmatrix}
0 & 0 & \pe B_{1} \\
0 &\r B_{2} & 0 \\
0 & 0 & \pe B_{3}
\end{pmatrix}.
\end{equation*}
Since~$\max \{ \r,\pe\} \le 1$, one can easily verify that there
exists a family $\{F_{\pa}\,\arrowvert\,\pa\in\PA\}$ 
uniformly bounded in~$C^{\infty}$, such that $F_{\pa}(0)=0$ and
\begin{equation*}%\label{repr:1}
\partial_{t}\ter = F_{\pa}(\Xi),
\end{equation*}
with $\Xi\defn\bigl(\vard,\ter,\nabla\ter,\cn\vard,\nabla\vari,
\eps\r\nabla^{2}\vard,\pe\nabla^{2}\variii\bigr)$. 

Using the Moser-type estimate~\eqref{esti:weightcomppre}, we find
\begin{equation*}
\norm{\partial_{t}\ter}_{\nhsc{s}} \le C(\norm{\Xi}_{\nh{s-1}}) 
\{1+\cpe\norm{\Xi}_{\nh{s}}\}.
\end{equation*}
To complete the proof note that, directly from the definitions~\eqref{defi:RR'} and the
assumption $(\r,\pe)\in [0,1]^{2}$, one has 
$\norm{\Xi}_{\nh{s-1}}\le R$ and $\cpe\norm{\Xi}_{\nh{s}} \le R+R'$.
%%\qed
\end{proof}

We are now prepared to prove Proposition~\ref{prop:HFuniform}. 
\subsection{Proof of Proposition~\ref{prop:HFuniform}}\label{subsection:proofHF}
Set $\cpe\defn\sqrt{\r+\pe}\le \sqrt{2}$ and 
$(\tvari,\tvard,\tvariii)\defn
\bigl(\Qr_{\eps\cpe}\vari,\Qr_{\eps\cpe}\vard,\Qr_{\eps\cpe}\variii\bigr)$. 
We first show that~$\tvar\defn(\tvari,\tvard,\tvariii)$
satisfies~\eqref{system:NSi'} for suitable source terms 
$f_1$, $f_{2}$, $f_{3}$ and coefficients $\beta_{1}$, $\beta_{2}$,
$\beta_{2}^{\sharp}$, $\beta_{3}$ and $G$. 
It readily follows from Notation~\ref{nota:sourceHF} that $(\tvari,\tvard,\tvariii)$ solves
\begin{equation*}
\left\{
\begin{aligned}
&g_{1}(\phi)\bigl(\partial_{t}\tvari+\vard\cdot\nabla {\tvari}\bigr)
+\frac{1}{\eps}\cn \tvard-\frac{\pe}{\eps} B_{1}(\phi)\tvariii=f^{\pa,\cpe}_{1,\HF}(\var),\\
&g_{2}(\phi)\bigl(\partial_{t}\tvard+\vard\cdot\nabla {\tvard}\bigr)
+\frac{1}{\eps}\nabla\tvari-\r B_{2}(\phi)\tvard=f^{\pa,\cpe}_{2,\HF}(\var),\\
&g_{3}(\phi)\bigl(\partial_{t}\tvariii
+\vard\cdot\nabla{\tvariii}+\tvard\cdot\nabla\variii\bigr)
+\cn\tvard -\pe B_{3}(\phi)\tvariii=f^{\pa,\cpe}_{3,\HF}(\var),
\end{aligned}
\right.
\end{equation*}
where $\phi=(\variii,\eps\vari)$. Set $\beta_{1}\defn\chi_{1}\beta$,
$\beta_{2}\defn\chi_{2}\zeta$, $\beta_{2}^{\sharp}\defn\chi_{2}\eta$, 
$\beta_{3}\defn\chi_{3}\beta$ and $G(\phi,\nabla\phi)=g_{3}(\phi)\nabla\variii$, to obtain 
that $(\tvari,\tvard,\tvariii)$ satisfies~\eqref{system:NSi'}
where
\begin{align*}
f_{1}&\defn f^{\pa,\cpe}_{1,\HF}(\var)+
f_{1,\HF}'+\Qr_{\eps\cpe}\Upsilon_{1} \quad\mbox{with}~~ 
f_{1,\HF}'\defn -\frac{\pe}{\eps}\nabla\chi_{1}(\eps\vari)\cdot\bigl(\beta(\variii)\nabla\tvariii\bigr),\\
f_{2}&\defn f^{\pa,\cpe}_{2,\HF}(\var) +f_{2,\HF}'\quad\mbox{with}~~ f_{2,\HF}'\defn
\r\chi_{2}(\eps\vari)\bigl\{2D\tvard \nabla\zeta(\variii)+\cn\tvard\nabla\eta(\variii)\bigr\},\\
f_{3}&\defn f^{\pa,\cpe}_{3,\HF}(\var)+
f_{3,\HF}'+\eps\Qr_{\eps\cpe}\Upsilon_{3} \quad\mbox{with}~~
f_{3,\HF}'\defn
\pe\chi_{3}(\eps\vari)\nabla\beta(\variii)\cdot\nabla\tvariii,
\end{align*}
where $\Upsilon_{i}$ are as in system~\eqref{system:NS2}.

By definition of the
norm~$\norm{\cdot}_{\Hr^{s}_{\pa}(T)}$ (see
Definition~\ref{defi:decoupling}), the Sobolev embedding Theorem and
Lemma~\ref{lemm:estdescoef123} give
$$
\norm{(\phi,\partial_{t}\phi,\nabla\phi,\cpe\nabla^{2}\phi,\vard,\nabla\vard)}
_{L^{\infty}([0,T]\times\xD)}\le C(\Omega),
$$
where $\Omega\defn\norm{\var}_{\Hr^{s}_{\pa}(T)}$. Similarly one has~$\norm{\phi(0)}_{L^{\infty}(\xD)}\le \Omega_{0}\defn
\norm{\var(0)}_{\Hr_{\pa,0}^{s}}$. 

Observe that Assumption~\ref{assu:structural} implies that the
conditions in Assumption~\ref{assu:linearized} are satisfied.

With these preliminary remarks in hands, Theorem~\ref{theo:L2slow} yields
\begin{equation*}
\norm{\tVar}_{\Hr^{0}_{\pa}(T)}\le 
C(\Omega_{0})e^{TC(\Omega)}\norm{\tVar(0)}_{\Hr^{0}_{\pa,0}}+C(\Omega)
\bigl\{\mathfrak{F}(T)+\mathfrak{F}'(T)+\mathfrak{F}''(T)\bigr\},
\end{equation*}
with
\begin{align*}
\mathfrak{F}(T)&\defn
\bigl\lVert (f^{\pa,\cpe}_{1,\HF}(\var), f^{\pa,\cpe}_{2,\HF}(\var))\bigr\rVert_{L^{1}(0,T;H^1_{\eps\cpe})}
+\bigl\lVert f^{\pa,\cpe}_{3,\HF}(\var)\bigr\rVert_{L^{1}(0,T;H^{1}_{\cpe})},\\
\mathfrak{F}'(T)&\defn
\bigl\lVert (f'_{1,\HF}, f'_{2,\HF})\bigr\rVert_{L^{1}(0,T;H^1_{\eps\cpe})}
+\bigl\lVert f'_{3,\HF}\bigr\rVert_{L^{1}(0,T;H^{1}_{\cpe})},\\
\mathfrak{F}''(T)&\defn 
\bigl\lVert \Qr_{\eps\cpe}\Upsilon_{1}\bigr\rVert_{L^{1}(0,T;H^1_{\eps\cpe})}
+\bigl\lVert \eps\Qr_{\eps\cpe}\Upsilon_{3}\rVert_{L^{1}(0,T;H^{1}_{\cpe})}.
\end{align*}

One has $\norm{\tvar}_{\Hr^{0}_{\pa}(T)}=
\norm{(\id-J_{\eps\cpe})\var}_{\Hr^{s}_{\pa}(T)}$ and similarly 
 $\norm{\tvar(0)}_{\Hr^{0}_{\pa,0}}=
\norm{(\id-J_{\eps\cpe})\var(0)}_{\Hr^{s}_{\pa,0}}\le \Omega_{0}$. 
Therefore, in view of the elementary inequality 
$x+y\le 2xe^{y}$ (for all $x\ge 1$ and $y\ge 0$), the proof of
Proposition~\ref{prop:HFuniform} reduces to establishing that 
$\mathfrak{F}(T)+\mathfrak{F}'(T)\le\sqrt{T}C(\Omega)$. 

The estimate $\mathfrak{F}(T)\le\sqrt{T}C(\Omega)$ immediately follows 
from the preliminaries. Indeed, define $R$ and $R'$ as
in~\eqref{defi:RR'}. By definition of $\cpe\defn\sqrt{\r+\pe}$, one has
$\r\le\sqrt{\r}\cpe$ and $\pe\le\sqrt{\pe}\cpe$. 
Hence, it is easily verified that 
Lemma~\ref{lemm:commHF} and Lemma~\ref{lemm:estdescoef123} imply
\begin{equation}\label{esti:sourceHFproof}
\bigl\lVert{f^{\pa,\cpe}_{1,\HF}(\var)\bigr\rVert}_{H^{1}_{\eps\cpe}} 
+\bigl\lVert{f^{\pa,\cpe}_{2,\HF}(\var)\bigr\rVert}_{H^{1}_{\eps\cpe}}
+\bigl\lVert{f^{\pa,\cpe}_{3,\HF}(\var)\bigr\rVert}_{H^{1}_{\cpe}}
\le C(R)\bigl\{1+R'\bigr\}.
\end{equation}
Since $C(\cdot)$ is non-decreasing, integrating and using the
Cauchy--Schwarz estimate, there results
\begin{align*}
\mathfrak{F}(T)
\le \sqrt{T}C\bigl(\norm{R}_{L^{\infty}(0,T)}\bigr)\bigl\{1+
\norm{R'}_{L^{2}(0,T)}\bigr\}.
\end{align*}
This in turn implies the desired estimate $\mathfrak{F}(T)\le\sqrt{T}C(\Omega)$ since, 
by definition, one has 
$\Omega\approx \norm{R}_{L^{\infty}(0,T)}+\norm{R'}_{L^{2}(0,T)}$.

Let us prove that, similarly,
$\mathfrak{F}'(T)\le\sqrt{T}C(\Omega)$. To do that it is sufficient to
prove that~\eqref{esti:sourceHFproof} holds true with
$f^{\pa,\cpe}_{i,\HF}(\var)$ replaced by $f'_{i,\HF}$. This in turn
follows from direct estimates. Indeed, observe that
$$
f'_{1,\HF}=\pe\Fr_{1}(\phi,\nabla\vari)\nabla\tvariii,\quad
f'_{2,\HF}=\r\Fr_{2}(\phi,\nabla\variii)\nabla\tvard,\quad
f'_{3,\HF}=\pe\Fr_{3}(\phi,\nabla\variii)\nabla\tvariii,
$$
for some $C^{\infty}$ functions $\Fr_{i}$ vanishing at the
origin. As already seen, one can give estimates for the coefficients $\Fr_{i}$
by combining the Sobolev embedding Theorem with the Moser
estimate~\eqref{esti:weightedcomp}. It is found that
\begin{equation}\label{esti:HFF'LHS}
\norm{(f'_{1,\HF},f'_{2,\HF})}_{H^{1}_{\eps\cpe}}+\norm{f'_{3,\HF}}_{H^{1}_{\cpe}}\le
C(R)\{\norm{\sqrt{\pe}\nabla\tvariii}_{\nhsc{1}}+\norm{\sqrt{\r}\nabla\tvard}_{H^{1}_{\eps\cpe}}\}.
\end{equation}
We next use $\Qr_{\eps\cpe}\le\Fi{}{s}$, to obtain
$$
\norm{\sqrt{\pe}\nabla\tvariii}_{\nhsc{1}}+\norm{\sqrt{\r}\nabla\tvard}_{H^{1}_{\eps\cpe}}\le
\sqrt{\pe}\norm{\variii}_{\nhsc{s+2}}+\sqrt{\r}\norm{\vard}_{H^{s+2}_{\eps\cpe}}\le R'.
$$
Consequently, the left-hand side of~\eqref{esti:HFF'LHS} is controlled
by $C(R)R'$.

To conclude the proof it remains to show that $\mathfrak{F}''(T)\le
\sqrt{T}C(\Omega)$. This is nothing new in that it follows from the estimates:
\begin{align*}
\mathfrak{F}''(T)&\les
\sqrt{T}\norm{(\Upsilon_{1},\Upsilon_{3})}_{L^{2}(0,T;H^{s+1}_{\eps\cpe})}\\
&\le
\sqrt{T}C\bigl(\norm{(\phi,\sqrt{\r}\nabla\vard)}_{L^{\infty}([0,T]\times\xD)}\bigr)
\bigl\{1+\sqrt{\r}\norm{\nabla\vard}_{L^{2}(0,T;H^{s+1}_{\eps\cpe})}\bigr\}\\
&\le
\sqrt{T}C(\norm{R}_{L^{\infty}(0,T)})\{1+\norm{R'}_{L^{2}(0,T)}\}\\
&\le \sqrt{T}C(\Omega). 
\end{align*}

\smallbreak
We have proved Proposition~\ref{prop:HFuniform}. 
This completes the analysis of the high frequency regime.

\section{Low frequency regime}\label{section:BF}
This section is devoted to the proof of {\em a priori\/} estimates 
in the low frequency region, which is the most delicate part.

\begin{proposition}\label{prop:LFuniform}
Given an integer~$s>1+d/2$, there exists a continuous non-decreasing 
function $C$ such that for all 
$\pa=(\eps,\r,\pe)\in\PA$, all $T\in [0,1]$ and all 
$U=(\vari,\vard,\variii)\in C^{1}([0,T];H^{\infty}(\xD))$ 
satisfying~\eqref{system:NS2},
\begin{equation}\label{prop:fastcomponent:1}
\norm{J_{\eps\cpe}\var}_{\Hr_{\pa}^{s}(T)}\le C(\Omega_{0})e^{(\sqrt{T}+\eps)C(\Omega)},
\end{equation}
where $\cpe\defn\sqrt{\r+\pe}$,
$\Omega_{0}\defn\norm{\var(0)}_{\Hr_{\pa,0}^{s}}$, 
$\Omega\defn \norm{\var}_{\Hr_{\pa}^{s}(T)}$ (see Definition~$\ref{defi:decoupling}$).
\end{proposition}

As alluded to previously, the nonlinear energy estimates
cannot be obtained from the $L^{2}$ estimates by an elementary
argument using differentiation of the equations with respect to 
spatial derivatives. For such problems a general strategy can 
be used. One first applies to the equations some 
operators based on $\edtt$. Next, one uses the special structure 
of the equations to estimate the spatial derivatives. 

This basic strategy has many roots, at least for hyperbolic problems
(see, e.g., \cite{TA,Iso,Se,Secchi}). 
For our purposes, the key point is that the 
hyperbolic behavior prevails in the low frequency regime. 
Yet, in sharp contrast with the Euler equations ($\r=\pe=0$), 
the form of the equations~\eqref{system:NS2} shows that the time derivative and 
the spatial derivatives have not the same weight. 

In particular, our analysis requires some preparation. 
We begin our discussion in
\S\ref{subsection:nonisotropic} by establishing some 
estimates which allows us to commute $J_{\eps\cpe}\edtt^{\ik}$ with
the equations. The latter task is achieved in~\S\ref{subsection:localizationLF}. 

With these preliminaries established, we can proceed to give an
estimate for $J_{\eps\cpe}\edtt^{s}\var$. 
The fast components $J_{\eps\cpe}(\cn\vard,\nabla\vari)$ are estimated 
next by using an
induction argument. 
To conclude, we give the estimates for the slow components $\variii$ and
$\curl \vard$.

\smallbreak
For the sake of notational clarity, in this section we deliberately omit the
terms $\Upsilon_{1}$ and $\eps\Upsilon_{3}$ in the
system~\eqref{system:NS2}. Nothing is changed in the statements of the
results, nor in their proofs. 

\subsection{Non-isotropic estimates}\label{subsection:nonisotropic}
The fact that the time derivative and the spatial derivatives 
have not the same weight is made precise by the
following lemma (whose easy proof is left to the reader).

\begin{lemma}
There is family $\bigl\{\,B_{\pa,\alpha}\,\arrowvert\,\pa\in\PA,~\alpha\in\xN^{d},~1\le\la\alpha\ra\le 2 \bigr\}$ 
uniformly bounded in $C^{\infty}\bigl(\xR^{N};\xR^{N\times N}\bigr)$ (where $N=(d+2)^{2}$) 
such that for all $\pa\in\PA$ and all smooth solution $(\vari,\vard,\variii)$ of~\eqref{system:NS2}, 
the function $\TER$ defined by
\begin{equation}\label{defi:intter}
\TER\defn \bigl(\ter,\partial_{t}\ter,\nabla\ter\bigr)\quad\mbox{where}\quad \ter\defn (\variii,\eps\vari,\eps\vard),
\end{equation}
solves
\begin{equation}\label{eq:TER}
\eps\partial_{t}\TER = \sum_{1\le j\le d}B_{\pa,j}(\TER)\partial_{j}\TER 
+\eps(\r+\pe)\sum_{1\le j,k\le d}\partial_{j}\bigl(B_{\pa,jk}(\TER)\partial_{k}\TER\bigr).
\end{equation}
\end{lemma}

We want to introduce an operator based on $\edtt$ which has the weight of a spatial
derivative. The previous result 
suggests introducing the following family of operators.
\begin{definition}For all~$\eps\ge 0$, $\cpe\ge 0$ and~$\ell\in\xN$, define the operators
$$
\Ze{\ell}\defn \Fi{\eps\cpe}{-\ell}(\eps\partial_{t})^{\ell},
$$
where, recall from~\eqref{defi:localizor}, that $\Fi{\eps\cpe}{-\ell}\defn\bigl(\id-(\eps\cpe)^{2}\Delta\bigr)^{-\ell/2}$.
\end{definition}

We will need the following technical ingredient.
\begin{lemma}\label{bignorms:comp}
Given~$F\in C^{\infty}(\xR^{n})$ satisfying $F(0)=0$ and
$\index>d/2$, there exists a function~$C(\cdot)$ such that 
for all $\eps\in [0,1]$, all $\cpe\in [0,2]$, all $T>0$, all vector-valued
function $U\in C^{\infty}([0,T];H^{\infty}(\xD))$ and all $\xN\ni\ik\le\index$, 
\begin{equation}\label{esti:Zecomp}
\bigl\lVert\Ze{\ik} F(U)\bigr\rVert_{\nh{\index-\ik}}
\le C\Bigl(\, \sum_{0\le\ell\le\ik}\bigl\lVert\Ze{\ell}U\bigr\rVert_{\nh{\index-\ell}} \,\Bigr) .
\end{equation}
\end{lemma}
\begin{proof}
To prove this claim, observe that $\edtt^{\ik}F(U)$ is a sum of terms of the form
$$
f(U)\edtt^{\ell_{1}} u_{1} \cdots \edtt^{\ell_{p}} u_{p},
$$
with~$p\le\ik$ and~$\sum_{1\le i\le p} \ell_{i} = \ik$. 
In this formula,~$f$ is a $C^{\infty}$ function and 
$u_{1},\ldots,u_{p}$ denote coefficients of~$U$. 

In order to estimate these terms, we use the following result 
(whose proof follows from Proposition~\ref{prop:Product} by
induction): 

Let 
$\alpha=(\alpha_{0},,\ldots,\alpha_{p}) \in\xN^{p+1}$ be such that 
$\sum_{i=0}^{p}\alpha_i =: \la\alpha\ra \le \index$, then
$$
\Bigl\lVert  \Fi{\eps\cpe}{-\la\alpha\ra} \prod_{i=0}^{p} V_i \Bigr\rVert_{\nh{\index-\la\alpha\ra}} 
\les \prod_{i=0}^{p} \lA \Fi{\eps\cpe}{-\alpha_i} V_i \rA_{\nh{\index-\alpha_{i}}},
$$
where the implicit constant depends  only on~$d$,$\index$ and~$\alpha$.

Set $\widetilde{f}\defn f-f(0)$. We apply the previous result with~$\alpha_{0}=0$,~$\alpha_{i}=\ell_{i}$ ($i\ge 1$), 
$V_{0} = \widetilde{f}(U)$ and~$V_{i} = \edtt^{\ell_{i}}u_{i}$ ($i\ge 1$). 
This yields
\begin{align*}
&\bigl\lVert\Fi{\eps\cpe}{-\ik}\bigl(f(U)\edtt^{\ell_{1}} u_{1} 
\cdots \edtt^{\ell_{p}} u_{p}\bigr) \bigr\rVert_{\nh{\index-\ik}} 
\\
&\qquad\qquad
\les \bigl(1+\bigl\lVert \widetilde{f}(U) \bigr\rVert_{\nh{\index}}\bigr) \bigl\lVert \Ze{\ell_{1}}u_{1} 
\bigr\rVert_{\nh{\index-\ell_{1}}} \cdots 
\bigl\lVert \Ze{\ell_{p}}u_{p} \bigr\rVert_{\nh{\index-\ell_{p}}}.
\end{align*}
Since $\Ze{0}=\id$, the estimate~\eqref{comp} implies
$\norm{\widetilde{f}(U)}_{\nh{\index}}\le
C\bigl(\norm{\Ze{0}U}_{\nh{\index}}\bigr)$. Which completes the proof.
%\qed
\end{proof}

Now we are in position to prove that $\Ze{1}$ has the weight
of a spatial derivative. The following result states that, for
$\ik\in\xN$, $\Ze{\ik}\TER$
satisfies the same estimates as $\Fi{}{\ik}F(\TER)$ does (where $F$ is
a given function).

\begin{proposition}\label{prop:notrouble}
Let $s>1+d/2$ be an integer. There exists a function~$C(\cdot)$ such
that for all $\pa=(\eps,\r,\pe)\in\PA$, all $T >0$ 
and all smooth solution $(\vari,\vard,\variii)\in
C^{\infty}([0,T];H^{\infty}(\xD))$ of~\eqref{system:NS2}
, if $\cpe \in [\r+\pe,2]$ then 
the function $\TER$ defined by~\eqref{defi:intter} satisfies 
\begin{gather}
\sum_{0\le \ell\le s}\norm{\Ze{\ell}\TER}_{\nh{s-\ell-1}}\le C \bigl(\norm{\TER}_{H^{s-1}}\bigr),\label{induction:TER1}\\
\sum_{0\le \ell\le s}\norm{\Ze{\ell}\TER}_{\nhsc{s-\ell}}\le C \bigl(\norm{\TER}_{H^{s-1}}\bigr)\norm{\TER}_{\nhsc{s}}.\label{induction:TER2}
\end{gather}
\end{proposition}
\begin{proof}
We prove by induction on $\ik\in\{0,\ldots,s\}$ that 
\begin{gather}
\sum_{0\le \ell\le \ik}\norm{\Ze{\ell}\TER}_{\nh{s-\ell-1}}\le C \bigl(\norm{\TER}_{H^{s-1}}\bigr),\label{induction:2TER1}\\
\sum_{0\le \ell\le \ik}\norm{\Ze{\ell}\TER}_{\nh{s-\ell}}\le C \bigl(\norm{\TER}_{H^{s-1}}\bigr)\norm{\TER}_{\nh{s}}.\label{induction:2TER2}
\end{gather}
Note that these results are obvious with $\ik=0$. 

Assume the results \eqref{induction:2TER1}--\eqref{induction:2TER1} at order $\ik <s$. 
By definition, one has $\Ze{\ik+1}=\Fi{\eps\cpe}{-1}\Ze{\ik}\edtt$. It thus follows from~\eqref{eq:TER} that
$$
\Ze{\ik+1}\TER=
\sum_{j,k}\Fi{\eps\cpe}{-1}\Ze{\ik}\bigl(B_{\pa,j}(\TER)\partial_{j}\TER\bigr)
+\bigl(\eps(\r+\pe)\partial_{j}\Fi{\eps\cpe}{-1}\bigr)\Ze{\ik}\bigl(B_{\pa,jk}(\TER)\partial_{k}\TER\bigr).
$$
Since $\cpe\ge\r+\pe$ one has
$\eps(\r+\pe)\partial_{j}\Fi{\eps\cpe}{-1}\les\id$ (see~\eqref{esti:handled}). 
Similarly one has $\Fi{\eps\cpe}{-1}\les\id$. 
The proof thus reduces to estimating terms having the form 
$\Ze{\ik} \bigl(B(\TER)\partial_{j}\TER\bigr)$. More precisely, it is
sufficient to prove that, for all integer $\ik< s$ and for all
smooth function $B$, one has
\begin{align}
&\norm{\Ze{\ik}\bigl(B(\TER)\partial_{j}\TER\bigr)}_{H^{s-\ik-2}} \le 
C\bigl(\norm{\TER}_{H^{s-1}}\bigr) ,\label{TER:todo1}\\
&\norm{\Ze{\ik}\bigl(B(\TER)\partial_{j}\TER\bigr)}_{H^{s-\ik-1}} \le 
C\bigl(\norm{\TER}_{H^{s-1}}\bigr)\norm{\TER}_{\nh{s}},\label{TER:todo2}
\end{align}
where $C(\cdot)$ depends only on a finite number of semi-norms of $B$ in $C^{\infty}$. 

Firstly, writing $B(\TER)\partial_{j}\TER$ as $\partial_{j}F(\TER)$ 
for some $C^{\infty}$ function $F$ such that $F(0)=0$, 
and using Lemma~\ref{bignorms:comp} with $\index=s-1$, we find
$$
\norm{\Ze{\ik}\bigl(B(\TER)\partial_{j}\TER\bigr)}_{H^{s-\ik-2}}
\les \norm{\Ze{\ik}F(\TER)}_{H^{s-1-\ik}}
\le 
C\Bigl(\sum_{0\le \ell\le\ik}\norm{\Ze{\ell} \TER}_{\nh{s-1-\ell}}\Bigr).
$$
As a consequence, the estimate~\eqref{TER:todo1} follows from the
induction hypothesis~\eqref{induction:2TER1}. Moving to the proof of~\eqref{TER:todo2}, we begin with
the Leibniz rule
$$
\edtt^{\ik}\bigl(B(\TER)\partial_{j}\TER\bigr)=\sum_{0\le\ell\le\ik}\binom{\ik}{\ell}\edtt^{\ik-\ell}B(\TER) \edtt^{\ell}\partial_{j}\TER.
$$
Let $\xN\ni\ell\le\ik$. 
Since~$s>1+d/2$ and $\ik\le s-1$, Proposition~\ref{prop:Product} applies with 
$\index=s-1$, $\indexg_{1}=m_{1}=\ik-\ell$ and $\indexg_{2}=m_{2}=\ell$. It yields
\begin{align*}
\norm{\Ze{\ik}\bigl(B(\TER)\partial_{j}\TER\bigr)}_{\nh{s-\ik-1}}\les\sum_{0\le\ell\le\ik}\norm{\Ze{\ik-\ell}B(\TER)}_{\nh{s-1-(\ik-\ell)}} 
\norm{\Ze{\ell}\partial_{j}\TER}_{\nh{s-1-\ell}}.
\end{align*}
Using Lemma~\ref{bignorms:comp} to estimate the first term in the summand, we get
\begin{align*}
\norm{\Ze{\ik}\bigl(B(\TER)\partial_{j}\TER\bigr)}_{\nh{s-\ik-1}}\le 
C\Bigl(\sum_{0\le\ell\le\ik}\norm{\Ze{\ell}\TER}_{\nh{s-\ell-1}}\Bigr)
\sum_{0\le\ell\le \ik}\norm{\Ze{\ell}\TER}_{\nh{s-\ell}}.
\end{align*}
Using the induction
hypotheses~\eqref{induction:2TER1}--\eqref{induction:2TER1}, we prove~\eqref{TER:todo2}.
%\qed
\end{proof}

We next prove a commutator estimate with gain of a factor~$\eps$. 
This technical ingredient is an analogue of the commutator
estimate~\eqref{commutator:estimate2} [with $\varrho=\para$] we used in the analysis of the high frequency regime. 

\begin{lemma}\label{bignorms:commutator}
Given $s>1+d/2$, there exists a constant $K$ such that for all
$\eps\in [0,1]$, all $\cpe\in [0,2]$, all $T>0$, all $\ik\in\xN$ such that $1\le\ik\le s$ and all 
$f,u\in C^{\infty}([0,T];H^{\infty}(\xD))$,
\begin{multline}\label{6.11}
\norm{\bigl[f,J_{\eps\cpe}\edtt^{\ik}]u}_{H^{s-\ik+1}_{\eps\cpe}}\\
\le K\eps
\Bigl\{\norm{f}_{H^{s}}+\sum_{\ell=0}^{\ik-1}\norm{\Ze{\ell}\partial_{t}
  f}_{H^{s-1-\ell}}\Bigr\}\Bigl\{\norm{\Ze{\ik}u}_{\nhsc{s-\ik}}+\sum_{\ell=0}^{\ik-1}\norm{\Ze{\ell}u }_{\nh{s-1-\ell}}\Bigr\}.
\end{multline}
\end{lemma}
\begin{proof}
The commutator $\bigl[f,J_{\eps\cpe}\edtt^{\ik}]u$ is expanded to
\begin{equation}\label{comest0}
\bigl[f,J_{\eps\cpe}]\edtt^{\ik}u+J_{\eps\cpe}\bigl[f,\edtt^{\ik}]u.
\end{equation}
The proof of~\eqref{6.11} is based on the tools we developed in
Section~\ref{Nonlinear}.

\textbf{a)} We first claim that
\begin{equation}\label{comest1}
\norm{\bigl[f,J_{\eps\cpe}]\edtt^{\ik} u}_{H^{s-\ik+1}_{\eps\cpe}}
\les\eps\cpe\norm{f}_{H^{s}}
\norm{\Ze{\ik}u}_{\nh{s-\ik}}.
\end{equation}
Since $s>1+d/2$ and $1\le\ik\le s$, Proposition~\ref{Friedrichs} applies with $(\para,m,\index,\indexg)$\footnote{Here, $m$ refers to
  the index used in the statement of Lemma~\ref{commutateur_S^1_f}.}
replaced by $(\eps\cpe,1,s,s-\ik)$. It yields
$$
\norm{\bigl[f,J_{\eps\cpe}]\edtt^{\ik} u}_{\nh{s-\ik}}
\les\eps\cpe\norm{f}_{H^{s}}
\norm{\Fi{\eps\cpe}{-(2s-\ik)}\edtt^{\ik}u}_{\nh{s-\ik}}.
$$
Since $\ik\le s$, there holds
$\Fi{\eps\cpe}{-(2s-\ik)}\le \Fi{\eps\cpe}{-\ik}$. 
Hence, we have the first half of~\eqref{comest1}, namely:
\begin{equation}\label{comest15}
\norm{\bigl[f,J_{\eps\cpe}]\edtt^{\ik} u}_{\nh{s-\ik}}
\les\eps\cpe\norm{f}_{H^{s}}
\norm{\Ze{\ik}u}_{\nh{s-\ik}}.
\end{equation}
The technique for obtaining the second half is similar. We first apply
Proposition~\ref{Friedrichs} with $m=0$, to obtain
\begin{equation}\label{comest2}
\norm{\bigl[f,J_{\eps\cpe}]\edtt^{\ik} u}_{\nh{s-\ik+1}}
\les\norm{f}_{H^{s}}\norm{\Ze{\ik}u}_{\nh{s-\ik}},
\end{equation}
and we next multiply~\eqref{comest2} by~$\eps\cpe$.

\noindent
\textbf{b)} Moving to the second term in~\eqref{comest0}, we claim that
\begin{multline}\label{com:Leibniz}
\norm{J_{\eps\cpe}\bigl[f,\edtt^{\ik}] u}_{H^{s-\ik+1}_{\eps\cpe}}\\
\les\eps\sum_{\ell=0}^{\ik-1}\norm{\Ze{\ell}\partial_{t} f}_{H^{s-1-\ell}}
\sum_{\ell=0}^{\ik-1}\norm{\Ze{\ell}u}_{\nh{s-1-\ell}}.
\end{multline}
Starting from the Leibniz rule, we get
\begin{equation}
\bigl[f,\edtt^{\ik}]u
= \eps \sum_{\ell=0}^{\ik-1}\binom{\ik}{\ell+1} 
\bigl(\edtt^{\ell}\partial_{t}f\bigr) \,\edtt^{\ik-1-\ell} u \label{Leibniz}.
\end{equation}
Let~$0\le\ell\le\ik-1$. 
Since~$s>1+d/2$ and~$1\le\ik\le s$, Proposition~\ref{prop:Product} applies with 
$\index=s-1$,~$\indexg_{1}=m_{1}=\ell$ and~$\indexg_{2}=m_{2}=\ik-1-\ell$. 
It yields
\begin{align*}
&\norm{\Fi{\eps\cpe}{-\ik+1}\bigl(\edtt^{\ell}\partial_{t}f)\,\edtt^{\ik-1-\ell}u\bigr)}_{\nh{s-\ik}} \\
&\qquad 
\les \norm{\Fi{\eps\cpe}{-\ell}\edtt^{\ell}\partial_{t}f}_{\nh{s-1-\ell}}
\norm{\Fi{\eps\cpe}{-(\ik-1-\ell)}\edtt^{\ik-1-\ell} u }_{\nh{s-1-(\ik-1-\ell)}}.
\end{align*}
By summing over all $0\le\ell\le\ik-1$, we obtain from the definition 
of the operators~$\Ze{\cdot}$ that 
$\norm{\Fi{\eps\cpe}{-\ik+1}\bigl[f,\edtt^{\ik}] u}_{H^{s-\ik}}$ 
is estimated by the right side of~\eqref{com:Leibniz}. To complete the
proof of~\eqref{com:Leibniz}, use 
the elementary estimate
\begin{equation}\label{JvsLambda}
\norm{J_{\eps\cpe} u}_{H^{s-\ik+1}_{\eps\cpe}} \les
\norm{J_{\eps\cpe}u}_{\nh{s-\ik}}\les \norm{\Fi{\eps\cpe}{-\ik+1} u}_{H^{s-\ik}}.
\end{equation}
%\qed
\end{proof}

\subsection{Notations}\label{subsection:LFnotations}
Let us pause here to set a few notations we use continually in the
sequel, as well as to make some running conventions.

\smallbreak
\noindent {\em Notations.}
\smallbreak
From now on, we consider a time $0<T\le 1$, a fixed triple of parameter
$\pa=(\eps,\r,\pe)\in\PA$ and a smooth solution $\var=(\vari,\vard,\variii)\in
C^{\infty}([0,T];H^{\infty}(\xD))$ of the
system~\eqref{system:NS2}. The notations $\phi$, $\ter$ and $\TER$ are shorthand notations for
\begin{equation}\label{defi:phiterTER}
\phi\defn(\variii,\eps\vari),\quad\ter\defn(\variii,\eps\vari,\eps\vard)\quad
\mbox{and}\quad\TER\defn (\ter,\partial_{t}\ter,\nabla\ter).
\end{equation}
Here and below, $s$ always denotes an
integer $s>1+d/2$. We set 
$$
\cpe\defn\sqrt{\r+\pe},\quad \Omega_{0}\defn
\norm{\var(0)}_{\Hr^{s}_{\pa,0}} \quad\mbox{and}\quad 
\Omega\defn \norm{\var}_{\Hr^{s}_{\pa}(T)},
$$
where the norms are defined in Definition~\ref{defi:decoupling}. As in~\eqref{defi:RR'}, set
\begin{equation}\label{defi:RR'2}
\begin{split}
R&\defn\norm{(\vari,\vard)}_{H^{s+1}_{\eps{\cpe}}}
+\norm{\variii}_{H^{s+1}_{\cpe}},\\
R'&\defn \sqrt{\r+\pe}\norm{(\cn\vard,\nabla\vari)}_{\nh{s}}
+\sqrt{\r}\norm{\nabla\vard}_{H^{s+1}_{\eps{\cpe}}}
+\sqrt{\pe}\norm{\nabla\variii}_{H^{s+1}_{{\cpe}}}.
\end{split}
\end{equation}
With these notations, one has $\Omega\approx\norm{R}_{L^{\infty}(0,T)}+\norm{R'}_{L^{2}(0,T)}$.

\smallbreak
\noindent {\em Conventions.}
\smallbreak

To say that a smooth non-decreasing function $C\colon
[0,+\infty)\rightarrow [1,+\infty)$ is {\em generic\/} is to say that 
$C(\cdot)$ is independent of $T$, $\pa$ and $(\vari,\vard,\variii)$. 
Given a generic function $C(\cdot)$, we denote by $\widetilde{C}$
the positive constant
\begin{equation}\label{defi:widetildeC}
\widetilde{C}\defn C(\Omega_{0})e^{(\sqrt{T}+\eps)C(\Omega)}.
\end{equation}
The factor $\eps$ is of no consequence but makes some arguments work 
more smoothly. We will often use the simple observation that 
for all generic function $C(\cdot)$, one has
$C(\Omega_{0})+(\sqrt{T}+\eps)C(\Omega)\le \widetilde{C}$.

To clarify matters, with these conventions,
Proposition~\ref{prop:LFuniform} is formulated concisely in the
following way: there exists a generic function $C(\cdot)$ such that 
$\norm{J_{\eps\cpe}U}_{\Hr^{s}_{\pa}(T)}\le
\widetilde{C}$.

\subsection{Localization in the low frequency region}\label{subsection:localizationLF}
\begin{definition}
Given $\ik\in\xN$, set
\begin{equation}\label{defi:oplocal}
\mathcal{X}^{\ik}\defn J_{\eps\cpe}\edtt^{\ik}.
\end{equation}
\end{definition}
\begin{nota}\label{nota:sourceLF}
Denote by $f^{\ik}_{\LF}=(f^{\ik}_{1,\LF},f^{\ik}_{2,\LF},f^{\ik}_{3,\LF})$ the commutator of the
equations~\eqref{system:NS2} and the operator $\mathcal{X}^{\ik}$:
\begin{alignat*}{4}
f^{\ik}_{1,\LF}&\defn\bigl[g_{1}(\phi),\mathcal{X}^{\ik}\bigr]\partial_{t}\vari
&&+\bigl[g_{1}(\phi)\vard,\mathcal{X}^{\ik}\bigr]\cdot\nabla\vari  
&&- \frac{\pe}{\eps}&&\bigl[B_{1}(\phi),\mathcal{X}^{\ik}\bigr]\variii,\\
f^{\ik}_{2,\LF}&\defn\bigl[g_{2}(\phi),\mathcal{X}^{\ik}\bigr]\partial_{t}\vard
&&+\bigl[g_{2}(\phi)\vard,\mathcal{X}^{\ik}\bigr]\cdot\nabla\vard
&&-\r&&\bigl[B_{2}(\phi),\mathcal{X}^{\ik}\bigr]\vard,\\[0.5ex]
f^{\ik}_{3,\LF}&\defn \bigl[g_{3}(\phi),\mathcal{X}^{\ik}\bigr]\partial_{t}\variii 
&&+\bigl[ g_{3}(\phi)\vard,\mathcal{X}^{\ik}\bigr]\cdot\nabla\variii
&&-\pe&&\bigl[B_{3}(\phi),\mathcal{X}^{\ik}\bigr]\variii.
\end{alignat*}
\end{nota}

The aim of this paragraph is to estimate $f^{\ik}_{\LF}$. To begin
with, consider the case when $m\ge 1$. 
\begin{lemma}\label{lemm:localizationLF}
There exists a generic function $C(\cdot)$ such that for all
$\ik\in\xN$ such that $1\le \ik\le s$,
\begin{equation*}
\norm{f^{\ik}_{\LF}}_{\Hr^{s-\ik}_{\pa,0}}\defn\norm{\bigl(f^{\ik}_{1,\LF},f^{\ik}_{2,\LF}\bigr)}_{H^{s-\ik+1}_{\eps\cpe}}+
\norm{f^{\ik}_{3,\LF}}_{H^{s-\ik+1}_{\cpe}}\le C(R)\{1+R'\},
\end{equation*}
where $R$ and $R'$ are as defined in~\eqref{defi:RR'2}.
\end{lemma}
\begin{proof}
To emphasize the role of the slow component $\TER$ (as defined in
\eqref{defi:phiterTER}), many bounds are given in terms of:
$$
\gamma\defn \norm{\TER}_{\nh{s-1}} \quad\mbox{and}\quad \Gamma\defn\norm{\TER}_{\nhsc{s}}.
$$
Directly from Lemma~\ref{lemm:estdescoef123}, one has 
$
\norm{\partial_{t}\ter}_{\nh{s-1}}\le C(R)$ and 
$\norm{\partial_{t}\ter}_{\nhsc{s}}\le C(R)\{1+R'\}$. 
Furthermore, directly from the definition of $\ter$ (see
\eqref{defi:phiterTER}), one has 
$\norm{(\ter,\nabla\ter)}_{\nhsc{s}}\le
\norm{(\vari,\vard)}_{H^{s+1}_{\eps\cpe}}+\norm{\variii}_{\nhsc{s+1}}\le
R$. Hence, 
\begin{equation}\label{obse:TERR}
\gamma=\norm{\TER}_{\nh{s-1}}\le C(R)\quad\mbox{and}\quad
\Gamma=\norm{\TER}_{\nhsc{s}}\le C(R)\{1+R'\}.
\end{equation}
Note that Proposition~\ref{prop:notrouble} applies 
since the condition $\cpe\ge \r+\pe$ is fulfilled.

\smallbreak 
\noindent {\sc{Step 1:}} Estimate for $f^{\ik}_{1,\LF}$. {\textbf{a)}}
We begin by proving 
\begin{equation}\label{esti:LLF1}
\norm{\bigl[g_{1}(\phi),\mathcal{X}^{\ik}\bigr]\partial_{t}\vari}_{H^{s-\ik+1}_{\eps\cpe}}
\le C(\gamma)\Gamma.
\end{equation}
Starting from Lemma~\ref{bignorms:commutator}, we find that the left-hand side is bounded by
\begin{equation}\label{A:C11}
\Bigl\{\norm{\tilde{g}_{1}(\phi)}_{H^{s}} + \sum_{\ell=0}^{\ik-1}\norm{\Ze{\ell}\partial_{t} g_{1}(\phi)}_{H^{s-1-\ell}}\Bigr\}
\Bigr\{\sum_{\ell=0}^{\ik}\norm{\Ze{\ell}(\eps\partial_{t}\vari)}_{\nhsc{s-\ell}}\Bigr\}\cdot
\end{equation}
Since $\TER=(\ldots,\edtt\vari,\ldots)$, the
estimate~\eqref{induction:TER2} implies that the second factor 
in~\eqref{A:C11} is estimated by $C(\gamma)\Gamma$. 
Moving to the first factor², note that 
$\partial_{t}g_{1}(\phi)=F(\phi,\partial_{t}\phi)$ 
for some $C^{\infty}$ function $F$ such that $F(0)=0$. 
As a consequence, Lemma~\ref{bignorms:comp} implies
\begin{equation}\label{A:C15}
\sum_{\ell=0}^{\ik-1}\norm{\Ze{\ell}\partial_{t}g_{1}(\phi)}_{\nh{s-1-\ell}}
\le C\Bigl(\sum_{\ell=0}^{\ik-1}\norm{\Ze{\ell}(\phi,\partial_{t}\phi)}_{\nh{s-1-\ell}}\Bigr).
\end{equation}
Again, since $\TER=(\phi,\ldots,\partial_{t}\phi,\ldots)$, the estimate~\eqref{induction:TER1} implies that the right-hand side 
in the previous estimate is bounded by $C(\gamma)$. 
It remains to estimate $\norm{\tilde{g}_{1}(\phi)}_{H^{s}}$. To do so write
\begin{align*}
\norm{\tilde{g}_{1}(\phi)}_{H^{s}}\le C\bigl(\norm{\phi}_{H^{s}}\bigr)\le C\bigl(\norm{(\ter,\nabla\ter)}_{H^{s-1}}\bigr)
\le C(\gamma).
\end{align*}
\noindent {\textbf{b)}} We next prove that 
\begin{equation}\label{esti:LLF2}
\norm{\bigl[g_{1}(\phi)\vard,\mathcal{X}^{\ik}\bigr]\cdot\nabla\vari}_{H^{s-\ik+1}_{\eps\cpe}}
\le C(\gamma)\bigl\{\norm{\nabla\vari}_{\nh{s-1}}+\Gamma\bigr\}.
\end{equation}
Lemma~\ref{bignorms:commutator} implies that the left-hand side is bounded by
\begin{equation}\label{A:C21}
\begin{aligned}
&\Bigl\{\norm{\eps g_{1}(\phi)\vard}_{H^{s}} 
+ \sum_{\ell=0}^{\ik-1}\norm{\Ze{\ell}\partial_{t} (\eps g_{1}(\phi)\vard)}_{H^{s-1-\ell}}\Bigr\}\\
&\qquad\qquad\qquad\times\Bigr\{\norm{\nabla\vari}_{\nh{s-1}}
+\sum_{\ell=1}^{\ik}\norm{\Ze{\ell}\nabla\vari}_{\nhsc{s-\ell}}\Bigr\}\cdot
\end{aligned}
\end{equation}
By definition $\ter=(\variii,\eps\vari,\eps\vard)$. Hence, 
one can rewrite $\eps g_{1}(\phi)\vard$ as 
$G(\ter)$ for 
some $C^{\infty}$ function $G$. Consequently, the first term in the above written product 
is the exact analogue of the first term in~\eqref{A:C11} with $g_{1}(\phi)$ replaced by $G(\ter)$. 
We thus obtain that this term is estimated by $C(\gamma)$.

Next, using the very definitions of $\Ze{\ell}$, write
\begin{equation}\label{A:C2split}
\sum_{\ell=1}^{\ik}\norm{\Ze{\ell}\nabla\vari}_{\nhsc{s-\ell}}=
\sum_{\ell=1}^{\ik}\norm{\Ze{\ell-1}
\Fi{\eps\cpe}{-1}\nabla(\eps\partial_{t}\vari)}_{\nhsc{s-\ell}}.
\end{equation}
Since $\Fi{\eps\cpe}{-1}\les\id$, this in turn implies
\begin{equation*}
\sum_{\ell=1}^{\ik}\norm{\Ze{\ell}\nabla\vari}_{\nhsc{s-\ell}}\le\sum_{\ell=0}^{\ik-1}\norm{\Ze{\ell}
(\eps\partial_{t}\vari)}_{\nhsc{s-\ell}}.
\end{equation*}
Which, as in the previous step \textbf{a)}, is estimated by $C(\gamma)\Gamma$.

\noindent {\textbf{c)}} 
To complete the estimate of $f^{\ik}_{1,\LF}$, we establish
\begin{equation}\label{esti:LLF3}
\frac{\pe}{\eps}\norm{\bigl[B_{1}(\phi),\mathcal{X}^{\ik}\bigr]\variii}_{H^{s-\ik+1}_{\eps\cpe}}
\le C(\gamma,R)\bigl\{\norm{\pe\variii}_{\nhsc{s+2}}+\Gamma\bigr\}.
\end{equation}
Parallel
to~\eqref{split:HFB1kl}, we 
decompose the commutator $\pe\eps^{-1}\bigl[B_{1}(\phi),\mathcal{X}^{\ik}\bigr]\variii$ as
\begin{equation}\label{A:C3split1}
\frac{\pe}{\eps}\bigl[ \Er_{1}(\phi),\mathcal{X}^{\ik}\bigr]\Delta\variii
+\frac{\pe}{\eps}\bigl[\Er_{2}(\phi,\nabla\phi),\mathcal{X}^{\ik}\bigr]\cdot\nabla\variii,
\end{equation}
where $\Er_{1}(\phi)\defn \chi_{1}(\eps\vari)\beta(\variii)$ and 
$\Er_{2}(\phi,\nabla\phi)\defn\chi_{1}(\eps\vari)\beta'(\variii)\nabla\variii$.

Replacing $\partial_{t}\vari$ with $\pe\eps^{-1}\Delta\variii$, the arguments given in~\textbf{a)} yield
$$
\frac{\pe}{\eps}
\norm{\bigl[\Er_{1}(\phi),\mathcal{X}^{\ik}\bigr]\Delta\variii}_{H^{s-\ik+1}_{\eps\cpe}}
\le C(\gamma)\sum_{\ell=0}^{\ik}\norm{\Ze{\ell}(\pe\Delta\variii)}_{\nhsc{s-\ell}}.
$$
We split the sum in the previous right-hand side as
$$
\norm{\pe\Delta\variii}_{\nhsc{s}}
+\sum_{\ell=1}^{\ik}\norm{\Ze{\ell}(\pe\Delta\variii)}_{\nhsc{s-\ell}}.
$$
The first term is bounded by $\norm{\pe\variii}_{\nhsc{s+2}}$. 
As regards the second term, we proceed as in~\textbf{b)}. Namely, 
write
\begin{equation*}
\sum_{\ell=1}^{\ik}\pe\norm{\Ze{\ell}\Delta\variii}_{\nhsc{s-\ell}}
=\sum_{\ell=1}^{\ik}\norm{\Ze{\ell-1}
\bigl(\eps\pe\Fi{\eps\cpe}{-1}\Delta\bigr)\partial_{t}\variii}_{\nhsc{s-\ell}}\le
\sum_{\ell=0}^{\ik-1}\norm{\Ze{\ell}\partial_{t}\variii}_{\nhsc{s-\ell}},
\end{equation*}
where we used $\eps\pe\Fi{\eps\cpe}{-1}\Delta\les\Fi{}{1}$ 
(which stems from $\pe\le\sqrt{\pe}\le\cpe$). 
Note that $\TER=(\ldots,\partial_{t}\variii,\ldots)$. Therefore, 
the estimate~\eqref{induction:TER2} implies 
that the right-hand side of the previous inequality is controlled 
by $C(\gamma)\Gamma$.

We now have to estimate the $H^{s-\ik+1}_{\eps\cpe}$-norm of the second term in~\eqref{A:C3split1}. 
Since $\pe\le\cpe$ and $\pe\le 1$, Lemma~\ref{bignorms:commutator} implies that this term is bounded by
\begin{equation*}
\Bigl\{\norm{\cpe\Er_{2}}_{\nh{s}}+\sum_{\ell=0}^{\ik-1}\norm{\cpe\Ze{\ell}\partial_{t}\Er_{2}}_{\nh{s-1-\ell}}\Bigr\}
\Bigl\{\norm{\nabla\variii}_{\nh{s-1}}+\sum_{\ell=1}^{\ik}\norm{\cpe\Ze{\ell}\nabla\variii}_{\nh{s-\ell}}\Bigr\}.
\end{equation*}
This in turn is bounded by $C(\gamma)\{C(R)+\Gamma\}$ since
\begin{align}
&\norm{\cpe\Er_{2}(\phi,\nabla\phi)}_{\nh{s}}\le C(R),
&\cpe\norm{\Ze{\ell}\partial_{t}
  \Er_{2}(\phi,\nabla\phi)}_{\nh{s-1-\ell}}\le C(\gamma)\Gamma,\label{6.1c1} \\
&\norm{\nabla\variii}_{\nh{s-1}}\le \gamma, &
\forall \ell\in [1,\ik], \quad \norm{\cpe\Ze{\ell}\nabla\variii}_{\nh{s-\ell}}\le C(\gamma).\label{esti:toimprove}
\end{align}
The first inequality in~\eqref{6.1c1} follows from~\eqref{esti:weightedcomp}. 
The second inequality in~\eqref{6.1c1} follows from~\eqref{TER:todo2}
since the term $\partial_{t}\Er_{2}(\phi,\nabla\phi)$ can be written $F(\TER)\nabla\TER$ for some
$C^{\infty}$ function $F$. The first inequality in~\eqref{esti:toimprove} is obvious. In order to prove the last one, 
as already seen, one can write
\begin{align*}
\norm{\cpe\Ze{\ell}\nabla\variii}_{\nh{s-\ell}}&\les 
\norm{\Ze{\ell-1}(\eps\cpe\Fi{\eps\cpe}{-1}\nabla)\partial_{t}\variii}_{\nh{s-\ell}}\\
&\les \norm{\Ze{\ell-1}\TER}_{\nh{s-1-(\ell-1)}}\le C(\gamma).
\end{align*}

\smallbreak
\noindent {\textbf{d)}} 
By combining~\eqref{esti:LLF1} with~\eqref{esti:LLF2}
and~\eqref{esti:LLF3}, we obtain
$$
\norm{f_{1,\LF}^{\ik}}_{H^{s-\ik+1}_{\eps\cpe}}\le
C(\gamma,R)\{\norm{\pe\variii}_{\nhsc{s+2}}+\Gamma\}, 
$$
so that the desired estimate 
$\| f_{1,\LF}^{\ik}\|_{H^{s-\ik+1}_{\eps\cpe}}\le C(R)\{1+R'\}$
follows from the observation~\eqref{obse:TERR}.

\noindent {\sc{Step 2:}} Estimate for
$f^{\ik}_{2,\LF}$. 
Note that one can obtain $f^{\ik}_{2,\LF}$ from
$f^{\ik}_{1,\LF}$ by replacing $\vari$ by $\vard$ and $\variii$ by
$\eps\vard$. Therefore, we are in the situation of the previous step and hence
conclude that
\begin{align}
\norm{\bigl[g_{2}(\phi),\mathcal{X}^{\ik}\bigr]\partial_{t}\vard}_{H^{s-\ik+1}_{\eps\cpe}}
&\le C(\gamma)\Gamma,\label{A:C4}\\ 
\norm{\bigl[g_{2}(\phi)\vard,\mathcal{X}^{\ik}\bigr]\cdot\nabla\vard}_{H^{s-\ik+1}_{\eps\cpe}}
&\le C(\gamma)\bigl\{\norm{\nabla\vard}_{\nh{s-1}}+\Gamma\bigr\},\label{A:C5}\\
\r\norm{\bigl[B_{2}(\phi),\mathcal{X}^{\ik}\bigr]\vard}_{H^{s-\ik+1}_{\eps\cpe}}
&\le C(\gamma,R)\bigl\{\norm{\eps\r\vard}_{\nhsc{s+2}}+\Gamma\bigr\}.\label{A:C6}
\end{align}

\noindent {\sc{Step 3:}} Estimate for
$f^{\ik}_{3,\LF}$. Note that, for all $u\in H^{\sigma}_{\cpe}$, 
\begin{equation}\label{esti:obvious17}
\norm{u}_{\nhsc{\sigma}}\le
\eps^{-1}\norm{u}_{H^{\sigma}_{\eps\cpe}}. 
\end{equation}
Hence,
\begin{multline*}
\norm{\bigl[g_{3}(\phi),\mathcal{X}^{\ik}\bigr]\partial_{t}\variii}_{H^{s-\ik+1}_{\cpe}}
+\pe\norm{\bigl[B_{3}(\phi),\mathcal{X}^{\ik}\bigr]\variii}_{H^{s-\ik+1}_{\cpe}}
\\
\le \frac{1}{\eps}\norm{\bigl[g_{3}(\phi),\mathcal{X}^{\ik}\bigr]\partial_{t}\variii}_{H^{s-\ik+1}_{\eps\cpe}}
+\frac{\pe}{\eps}\norm{\bigl[B_{3}(\phi),\mathcal{X}^{\ik}\bigr]\variii}_{H^{s-\ik+1}_{\eps\cpe}}.
\end{multline*}
The first term in the right-hand side is estimated as in \textbf{a)}
above [replacing $\partial_{t}\vari$ by
$\eps^{\-1}\partial_{t}\variii$] and the second 
term has been estimated in \textbf{c)}.

For technical reasons, the estimate for
$[g_{3}(\phi)\vard,\mathcal{X}^{m}]\cdot\nabla\variii$ is somewhat
more complicated. We argue as in the proof of
Lemma~\ref{bignorms:commutator}. Set $f\defn g_{3}(\phi)\vard$ and
$u=\nabla\variii$. 
We begin by
splitting the commutator $[f,\mathcal{X}^{m}]\cdot u$ as
$$
P + Q \defn J_{\eps\cpe}[f,\edtt^{m}]\cdot u + [f,J_{\eps\cpe}]\cdot\edtt^{m}u.
$$
By combining~\eqref{comest15} with \eqref{comest2} multiplied by $\cpe$,
we find
$$
\norm{Q}_{\nhsc{s-\ik+1}}\les
\norm{f}_{\nh{s}}\norm{\Ze{\ik}u}_{\nhsc{s-\ik}}\le 
C(\gamma,R)\Gamma,
$$
where we used~\eqref{esti:toimprove}. 

Our next task is to show a similar estimate for $P$. To do so 
we decompose $P$ into two parts:
\begin{equation}\label{split:PP1P2}
P_{1}+P_{2}\defn \Bigl\{J_{\eps\cpe}[f,\edtt^{m}]\cdot u
-J_{\eps\cpe}\bigl(u\cdot\edtt^{m}f\bigr)\Bigr\}
+J_{\eps\cpe}\bigl(u\cdot\edtt^{m}f\bigr).
\end{equation}

Let us prove that $\norm{P_{1}}_{H^{s-\ik+1}_{\cpe}}\le C(\gamma)$. 
In light of~\eqref{esti:obvious17}, all we need to prove is that
\begin{equation}\label{Eq:indirect}
\norm{P_{1}}_{H^{s-\ik+1}_{\eps\cpe}}\le \eps C(\gamma).
\end{equation}
Repeat the proof of Lemma~\ref{bignorms:commutator}, to obtain
\begin{equation*}
\norm{P_{1} }_{H^{s-\ik+1}_{\eps\cpe}}
\les \eps\sum_{\ell=1}^{\ik-1}\norm{\Ze{\ell}\partial_{t} f}_{H^{s-1-\ell}}\norm{\Ze{\ell}u }_{\nh{s-1-\ell}}.
\end{equation*}
The sum differs from the one that appears in
Lemma~\ref{bignorms:commutator} in that it 
is indexed by $\ell\ge 1$ instead of $\ell\ge 0$. This fact allows us to
write 
\begin{equation*}
\norm{P_{1} }_{H^{s-\ik+1}_{\eps\cpe}}\les\eps\sum_{\ell=1}^{\ik-1}\norm{\Ze{\ell}\edtt f}_{H^{s-1-\ell}}\norm{\Ze{\ell-1}\partial_{t}u }_{\nh{s-1-\ell}}.
\end{equation*}
Let $1\le \ell\le\ik-1$. 
Write $\edtt f$ as $F(\TER)$ for some $C^{\infty}$ function $F$ such
that $F(0)=0$. By combining Lemma~\ref{bignorms:comp} 
and the estimate~\eqref{induction:TER1}, we obtain 
$$
\norm{\Ze{\ell}\edtt f}_{H^{s-1-\ell}}\le C(\gamma).
$$
Moving to the estimate of $\Ze{\ell-1}\partial_{t}u$, use the very definitions
of $u=\nabla\variii$ and
$\TER=(\ldots,\partial_{t}\variii,\ldots)$, to obtain thanks to~\eqref{induction:TER1}
$$
\norm{\Ze{\ell-1}\partial_{t}u }_{\nh{s-1-\ell}}\le
\norm{\Ze{\ell-1}\TER}_{\nh{s-\ell}}
=\norm{\Ze{\ell-1}\TER}_{\nh{s-(\ell-1)-1}}\le C(\gamma).
$$

We have proved~\eqref{Eq:indirect}, so to conclude it remains only to estimate the second term $P_{2}$
in~\eqref{split:PP1P2}. This is accomplished using 
\begin{equation*}
\norm{P_{2}}_{\nhsc{s-\ik+1}}\le \norm{u}_{\nhsc{s}}\norm{\Ze{\ik-1}\edtt f}_{\nhsc{s-\ik+1}},
\end{equation*}
as the reader can verify, yielding the bound 
$\norm{P_{2}}_{\nhsc{s-\ik+1}}\le C(\gamma,R)\Gamma$. 

This completes the proof of Lemma~\ref{lemm:localizationLF}.
%\qed
\end{proof}

Note that in the case when $\ik=0$, the previous method does not work as it
stands, as the estimate~\eqref{comest1} is no longer correct. 

Now we give an estimate valid for all $0\le\ik\le s-1$. 
\begin{lemma}\label{lemm:localizationLF2}
There exists a generic function $C(\cdot)$ such that for all
$\ik\in\xN$ such that $0\le \ik\le s-1$,
\begin{equation}\label{esti:localizationLF2}
\norm{f^{\ik}_{\LF}}_{\Hr^{s-1-\ik}_{\pa,0}}\defn\norm{\bigl(f^{\ik}_{1,\LF},f^{\ik}_{2,\LF}\bigr)}_{H^{s-\ik}_{\eps\cpe}}+
\norm{f^{\ik}_{3,\LF}}_{H^{s-\ik}_{\cpe}}\le C(R).
\end{equation}
\end{lemma}
%We now give an estimate for $f^{0}_{\LF}$.
%\begin{lemma}\label{lemm:localizationLF2}
%There exists a generic function $C(\cdot)$ such that
%\begin{equation}\label{esti:localizationLF2}
%\norm{f^{0}_{\LF}}_{H^{s}}\le C(R)\{1+R'\}.
%\end{equation}
%\end{lemma}
%Note that this estimate is weaker than the estimate obtained by
%formally setting $m=0$ in Lemma~\ref{lemm:localizationLF}. 
%Note that in the case when $\ik=0$, the previous method does not work as it
%stands, as the estimate~\eqref{comest1} is no longer correct.
The estimate~\eqref{esti:localizationLF2} is nothing new in that it can be deduced by following
the proof of Lemma~\ref{lemm:localizationLF}. One only has to use the following analogue 
of the calculus inequality in Lemma~\ref{bignorms:commutator}:
\begin{multline*}
\norm{\bigl[f,J_{\eps\cpe}\edtt^{\ik}]u}_{H^{s-\ik}_{\eps\cpe}} \le \\
K\eps
\Bigl\{\norm{f}_{H^{s}}+\sum_{\ell=0}^{\ik-1}\norm{\Ze{\ell}\partial_{t}
  f}_{H^{s-2-\ell}}\Bigr\}\Bigl\{\norm{\Ze{\ik}u}_{\nh{s-1-\ik}}+\sum_{\ell=0}^{\ik-1}\norm{\Ze{\ell}u }_{\nh{s-1-\ell}}\Bigr\}.
\end{multline*}
In the case when $\ik=0$ the sum $\sum_{0}^{-1}$ is interpreted as
$0$. The index $s-2-\ell$ in the second term of the first set of
parentheses is not a typographical error. It is of use to us for the
estimate of $\partial_{t}\Er_{2}(\phi,\nabla\phi)$ where $\Er_{2}$ is
as in~\eqref{6.1c1}.

%%%%%%%%%%%%%%%%%%%%%%%%%%%%%%%%%%%%%%%%%%%%%%%%%%%%%%%%%%%%%%%%%
%%%%%%%%%%%%%%%%%%%%%%%%%%%%%%%%%%%%%%%%%%%%%%%%%%%%%%%%%%%%%%%%%
%%%%%%%%%%%%%%%%%%%%%%%%%%%%%%%%%%%%%%%%%%%%%%%%%%%%%%%%%%%%%%%%%

\subsection{The fast components}\label{subsection:induction}
We give here the estimates for the fast components $\cn\vard$ and
$\nabla\vari$. We make use of the notations introduced in~\S\ref{subsection:LFnotations}.
\begin{nota}
For all integer $\ik\le s$, set $\var_{\ik}\defn \mathcal{X}^{\ik}\var$
where recall $\mathcal{X}^{m}=J_{\eps\cpe}\edtt^{\ik}$.
\end{nota}
As a preliminary step towards the estimate of
$(\cn\vard,\nabla\vari)$ we estimate the $\Hr^{0}_{\pa}(T)$-norm of
$\var_{\ik}$.
{\nobreak
\begin{lemma}\label{lemm:induct1}
For all integer $\ik\le s$, there is a generic function $C$ such that
$\norm{\var_{\ik}}_{\Hr^{0}_{\pa}(T)} \le \widetilde{C}$, where 
$\widetilde{C}$ is as defined in~\eqref{defi:widetildeC}.
\end{lemma}
Having estimated the commutators $f^{\ik}_{\LF}$, 
this result can be deduced by following the end of proof of 
Proposition~\ref{prop:HFuniform} (see~\S\ref{subsection:proofHF}). 
We therefore only indicate the points at which the argument is
slightly different.
\begin{proof}
It readily follows from Notation~\ref{nota:sourceLF} that 
$(\tvari,\tvard,\tvariii)\defn(\vari_{\ik},\vard_{\ik},\variii_{\ik})$ satisfies the linearized system~\eqref{system:NSi} with 
\begin{align*}
f_{1}&\defn f^{\ik}_{1,\LF}+
f_{1,\LF}'\quad\mbox{where}\quad 
f_{1,\LF}'\defn -\frac{\pe}{\eps}\nabla\chi_{1}(\eps\vari)\cdot\bigl(\beta(\variii)\nabla\variii_{\ik}\bigr),\\
f_{2}&\defn f^{\ik}_{2,\LF} +f_{2,\LF}'\quad\mbox{where}\quad f_{2,\LF}'\defn
\r\chi_{2}(\eps\vari)\bigl\{2D\vard_{\ik} \nabla\zeta(\variii)+\cn\vard_{\ik}\nabla\eta(\variii)\bigr\},\\
f_{3}&\defn f^{\ik}_{3,\LF}+
f_{3,\LF}'\quad\mbox{where}\quad
f_{3,\LF}'\defn \pe\chi_{3}(\eps\vari)\nabla\beta(\variii)\cdot\nabla\variii_{\ik}.
\end{align*}
recalling that we deliberately omit the
terms $\Upsilon$ and $\eps\Upsilon$ in the system~\eqref{system:NS2}. 

Applying Theorem~\ref{theo:L2}, we get
\begin{equation*}
\norm{\var_{\ik}}_{\Hr^{0}_{\pa}(T)}\le 
C(\Omega_{0})e^{TC(\Omega)}\norm{\var_{\ik}(0)}_{\Hr^{0}_{\pa,0}}+C(\Omega)\mathfrak{F}(T)+C(\Omega)\mathfrak{F}'(T),
\end{equation*}
with
\begin{align*}
\mathfrak{F}(T)&\defn
\bigl\lVert(f^{\ik}_{1,\LF},f^{\ik}_{2,\LF})\bigr\rVert_{L^{1}(0,T;H^1_{\eps\cpe})}
+\bigl\lVert f^{\ik}_{3,\LF}\bigr\rVert_{L^{1}(0,T;H^{1}_{\cpe})},\\
\mathfrak{F}'(T)&\defn
\bigl\lVert (f'_{1,\LF},f'_{2,\LF})\bigr\rVert_{L^{1}(0,T;H^1_{\eps\cpe})}
+\bigl\lVert f'_{3,\LF}\bigr\rVert_{L^{1}(0,T;H^{1}_{\cpe})}.
\end{align*}

The proof thus reduces to establishing that 
$\norm{\var_{\ik}(0)}_{\Hr^{0}_{\pa,0}}\le C(\Omega_{0})$ and $\mathfrak{F}(T)+\mathfrak{F}'(T)\le\sqrt{T}C(\Omega)$ . 
To fix matters, we concentrate on the hardest case when $\ik=s$. 
Note that the conclusion of Lemma~\ref{lemm:localizationLF} is the exact
analogue of~\eqref{esti:sourceHFproof}. Hence, 
one has $\mathfrak{F}(T)\le
\sqrt{T}C(\Omega)$. 
As in~\S\ref{subsection:proofHF}, all that has to be done in order to
prove $\mathfrak{F}'(T)\le\sqrt{T}C(\Omega)$ is to check 
$\sqrt{\pe}\|\nabla\variii_{s}\|_{\nhsc{1}}+\sqrt{\r}\|\nabla\vard_{s}\|_{H^{1}_{\eps\cpe}}\le
C(R)\{1+R'\}$. To do so, 
using the definition of $\variii_{s}=J_{\eps\cpe}\edtt^{s}\variii$, we first rewrite 
$\nabla\variii_{s}$ as $\eps
J_{\eps\cpe}\edtt^{s-1}\nabla\partial_{t}\variii$. 
Next, by combining the estimate $\norm{\eps J_{\eps\cpe}u}_{\nhsc{1}}\le
\norm{J_{\eps\cpe}u}_{H^{1}_{\eps\cpe}}\les \norm{J_{\eps\cpe}u}_{\nhz}$ with the inequality $\sqrt{\pe}\le \cpe$ and the
definition of $\TER=(\ldots,\partial_{t}\variii,\ldots)$, we obtain
$$
\sqrt{\pe}\norm{\nabla\variii_{s}}_{\nhsc{1}}\les
\cpe\norm{\Ze{s-1}\nabla\TER}_{\nhz}\le\norm{\Ze{s-1}\TER}_{\nhsc{1}}.
$$
Analogous computations lead to
\begin{equation*}
\sqrt{\r}\|\nabla\vard_{s}\|_{H^{1}_{\eps\cpe}}\les
\norm{\vard_{s}}_{\nhsc{1}}
\les \norm{\Ze{s-1}\edtt\vard}_{\nhsc{1}}
\le \norm{\Ze{s-1}\TER}_{\nhsc{1}}.
\end{equation*}
Hence, the desired bound follows from~\eqref{induction:TER2} and
~\eqref{obse:TERR}. 

The technique for estimating the initial data is similar. 
Indeed, one has $\norm{\var_{s}}_{\Hr^{0}_{\pa,0}}\les
\norm{\Ze{s-1}\TER}_{\nhz}$, as the reader can verify, yielding the
bound $\norm{\var_{s}(0)}_{\Hr^{0}_{\pa,0}}\le C(R(0))=C(\Omega_{0})$.
\end{proof}

We now come to the main estimates of this part. %the induction argument. 
\begin{nota}
Define
$\norm{u}_{\mathcal{K}^{\sigma}_{\cpe}(T)}\defn
\norm{u}_{L^{\infty}(0,T;H^{\sigma-1})}
+\cpe\norm{u}_{L^{2}(0,T;H^{\sigma})}$.
\end{nota}
\begin{lemma}
Let $\widetilde{\var}\defn (\tvari,\tvard,\tvariii)$
be a 
solution of the system:
\begin{equation}\label{system:NS5}
\left\{
\begin{aligned}
&g_{1}(\phi)\partial_{t}\tvari
+\frac{1}{\eps}\cn \tvard-\frac{\pe}{\eps}B_{1}(\phi)\tvariii=f_{1},\\
&g_{2}(\phi)\partial_{t}\tvard
+\frac{1}{\eps}\nabla\tvari-\r B_{2}(\phi)\tvard=f_{2},\\
&g_{3}(\phi)\partial_{t}\tvariii
+\cn\tvard -\pe B_{3}(\phi)\tvariii=f_{3}.
\end{aligned}
\right.
\end{equation}
If the Fourier transform of $\tvar$ is supported in the ball
$\{\la\xi\ra\le 2/\eps\cpe\}$, then there exists a generic function
$C(\cdot)$ such that for all $\sigma \in [1,s]$,
\begin{equation}\label{esti:LFinduction}
\begin{split}
&\norm{\tvari}_{\Kr_{\cpe}^{\sigma+1}(T)} 
+\norm{\cn\tvard}_{\Kr_{\cpe}^{\sigma}(T)}\\[0.5ex]
&\quad\le \widetilde{C}\norm{\edtt\tvari}_{\Kr_{\cpe}^{\sigma}(T)} 
+\widetilde{C}\norm{\edtt\cn\tvard}_{\Kr_{\cpe}^{\sigma-1}(T)}\\[0.5ex]
&\quad\quad +\widetilde{C}\norm{\tvari}_{L^{\infty}(0,T;L^{2})}+
\widetilde{C} \norm{\tvariii(0)}_{\nhsc{\sigma+1}}
+\eps C(\Omega)\norm{\r\tvard}_{\Kr_{\cpe}^{\sigma+1}(T)}\\[0.5ex]
&\quad\quad + \eps C(\Omega)\norm{(f_{1},f_{2})}_{\Kr_{\cpe}^{\sigma}(T)}
+\cpe \widetilde{C}\norm{f_{3}}_{L^{2}(0,T;H^{\sigma})}.
\end{split}
\end{equation}
\end{lemma}
\begin{remark} This result is an elaboration of Lemma $2.10$ in~\cite{TA}.
\end{remark}
\begin{proof}
\noindent {\textbf{a)}} For further references, we first
prove three estimates.

\smallbreak
\noindent {\textbf{1)}} Let $\index>d/2$ and $\sigma \in [0,\index]$. There exists a
constant $K$ such that for all $\varrho\ge 0$ 
and for all $(u_{1},u_{2})\in H^{\index}\times H^{\sigma}_{\varrho}(\xD)$,
\begin{equation}\label{6.15}
\norm{u_{1}u_{2}}_{H^{\sigma}_{\varrho}}\le K
\norm{u_{1}}_{H^{\index}}\norm{u_{2}}_{H^{\sigma}_{\varrho}}.
\end{equation}
To prove this result, we use a standard Moser's estimate. Since
$\index>d/2$,  the product maps continuously $H^{\index}\times H^{r}$ 
to $H^{r}$ for all $r \in [-\index,\index]$. Hence, one has
\begin{align*}
\norm{u_{1}u_{2}}_{H^{\sigma}_{\varrho}}&\defn\norm{u_{1}u_{2}}_{H^{\sigma-1}}+\varrho\norm{u_{1}u_{2}}_{H^{\sigma}}\\
&\les \norm{u_{1}}_{H^{\index}}
(\norm{u_{2}}_{H^{\sigma-1}}+\varrho\norm{u_{2}}_{\nh{\sigma}}) = \norm{u_{1}}_{H^{\index}}\norm{u_{2}}_{H^{\sigma}_{\varrho}}.
\end{align*}

\smallbreak
\noindent {\textbf{2)}} Let $\sigma\in [0,s]$ and $q\in
[1,+\infty]$. Given a $C^{\infty}$ function $F$, 
there exists a generic function $C$ such that for all $\varrho\ge 0$, 
and for all $u\in L^{q}(0,T;H^{\sigma}_{\varrho})$,
\begin{equation}\label{6.16}
\norm{F(\ter)u}_{L^{q}(0,T;H^{\sigma}_{\varrho})}\les \widetilde{C}\norm{u}_{L^{q}(0,T;H^{\sigma}_{\varrho})},
\end{equation}
where $\ter$ is as defined in~\eqref{defi:phiterTER}. 

In light of~\eqref{6.15}, to prove this estimate we need only show that
\begin{equation}\label{6.16alpha}
\norm{F(\ter)}_{L^{\infty}(0,T;H^{s})}\le \widetilde{C}.
\end{equation}
This will be established (independently) in~\eqref{esti:goodcurlsplit3} below.

\smallbreak
\noindent {\textbf{2')}} Let us infer from~\eqref{6.16} that, 
for all $\sigma\in [0,s]$,
\begin{equation}\label{6:16'}
\norm{F(\ter)u}_{\Kr^{\sigma}_{\cpe}(T)}\le \widetilde{C}\norm{u}_{\Kr^{\sigma}_{\cpe}(T)},
\end{equation}
To see this write
\begin{align*} 
\norm{F(\ter)u}_{\Kr^{\sigma}_{\cpe}(T)}&\defn
\norm{F(\ter)u}_{L^{\infty}(0,T;H^{\sigma-1})}+\cpe\norm{F(\ter)u}_{L^{2}(0,T;H^{\sigma})}\\
&\le
\norm{F(\ter)u}_{L^{\infty}(0,T;H^{\sigma-1})}+\norm{F(\ter)u}_{L^{2}(0,T;H^{\sigma}_{\cpe})}\\
&\le
\widetilde{C}\bigl\{\norm{u}_{L^{\infty}(0,T;H^{\sigma-1})}+\norm{u}_{L^{2}(0,T;H^{\sigma}_{\cpe})}\bigr\}
\le \widetilde{C} \norm{u}_{\Kr^{\sigma}_{\cpe}(T)}.
\end{align*}

\smallbreak
\noindent {\textbf{3)}} Let $\sigma\in [0,s]$, $f\in
H^{\sigma-1}(\xD)$ and $\gamma\colon \xD\rightarrow (0,+\infty)$ be a function bounded from below by a positive
constant. Furthermore, suppose
$\tilde{\gamma}\defn\gamma-\underline{\gamma}\in H^{s}(\xD)$ for some
constant $\underline{\gamma}$. 
We claim that if $u\in H^{\sigma}(\xD)$ satisfies $\cn(\gamma \nabla u)=f$, then 
there exists a constant $K=K(d,s)$ such that
\begin{equation}\label{6.17}
\norm{u}_{\nh{\sigma+1}} \le 
K\norm{\gamma^{-1}}_{L^{\infty}}\bigl\{\norm{f}_{\nh{\sigma-1}}+
\norm{\tilde{\gamma}}_{\nh{s}}\norm{u}_{\nhz}\bigr\}+\norm{u}_{\nhz}.
\end{equation}
Commuting the equation $\cn(\gamma \nabla u)=f$ with $\Fi{}{\sigma}$
and using the commutator estimate~\eqref{commutator:estimateusual} to
bound $\cn([\gamma,\Fi{}{\sigma}] \nabla u)$, we see that the proof
of~\eqref{6.17} can be reduced to the
special case $\sigma=0$. On the other hand, the case $\sigma=0$ is
immediate by usual integration by parts and duality arguments. 

\smallbreak
\noindent {\textbf{b)}} 
Hereafter, RHS denotes the right-hand side
of~\eqref{esti:LFinduction}. To prove that
\begin{equation}\label{6.205}
\norm{\cn\tvard}_{\Kr^{\sigma}_{\cpe}(T)}\le {\rm RHS},
\end{equation}
we begin by showing that
\begin{equation}\label{6.21}
\norm{\cn\tvard}_{\Kr^{\sigma}_{\cpe}(T)}\le
\widetilde{C}\bigl\{\norm{\edtt\tvari}_{\Kr^{\sigma}_{\cpe}(T)}+
\norm{\pe\tvariii}_{\Kr^{\sigma+2}_{\cpe}(T)}+\eps\norm{f_{1}}_{\Kr^{\sigma}_{\cpe}(T)}\bigr\}. 
\end{equation}
Rewrite $\cn\tvard$ as 
\begin{equation*}
-g_{1}(\phi)\edtt\tvari 
+\pe \Er_{1}(\phi)\Delta\tvariii +\pe\Er_{2}(\phi,\nabla\phi)\nabla\tvariii+\eps f_{1},
\end{equation*}
where $\Er_{1}$ and $\Er_{2}$ are as in~\eqref{A:C3split1}. Using~\eqref{6:16'}, one has
\begin{align*}
\norm{g_{1}(\phi)\edtt\tvari}_{\Kr^{\sigma}_{\cpe}(T)}&\le
\widetilde{C} \norm{\edtt\tvari}_{\Kr^{\sigma}_{\cpe}(T)}, \\
\norm{\pe \Er_{1}(\phi)\Delta\tvariii}_{\Kr^{\sigma}_{\cpe}(T)}&\le
\widetilde{C} \norm{\pe\tvariii}_{\Kr^{\sigma+2}_{\cpe}(T)}.
\end{align*}
Therefore, to infer~\eqref{6.21} it remains only to estimate the $\Kr^{\sigma}_{\cpe}(T)$-norm
of $\pe \Er_{2}(\phi,\nabla\phi)\nabla\tvariii$. 
Using the very definition of $\Er_{2}(\phi,\nabla\phi)$ and the
estimate~\eqref{6.16alpha}, one has
$$
\norm{\Er_{2}(\phi,\nabla\phi)}_{L^{\infty}(0,T;H^{s-1})}\les \widetilde{C}.
$$
Since $s-1>d/2$, the estimate~\eqref{6.15} applies with
$\index\defn s-1$. It yields
\begin{equation}
\norm{\Er_{2}(\phi,\nabla\phi)\nabla\tvariii}_{L^{\infty}(0,T;H^{\sigma-1})}
\le
\widetilde{C}\norm{\nabla\tvariii}_{L^{\infty}(0,T;H^{\sigma-1})}
\le \widetilde{C}\norm{\tvariii}_{L^{\infty}(0,T;H^{\sigma})}\label{6.41}.
\end{equation}
Moreover, one can easily verify that 
\begin{align*}
&\cpe\norm{\Er_{2}(\phi,\nabla\phi)\nabla\tvariii}_{L^{2}(0,T;H^{\sigma})}\notag
\\
&\qquad\le
\norm{\cpe\Er_{2}(\phi,\nabla\phi)}_{L^{2}(0,T;H^{s})}\norm{\nabla\tvariii}_{L^{\infty}(0,T;H^{\sigma})}\notag
\\
&\qquad\le
\sqrt{T}\norm{\Er_{2}(\phi,\nabla\phi)}_{L^{\infty}(0,T;H^{s}_{\cpe})}\norm{\tvariii}_{L^{\infty}(0,T;H^{\sigma+1})}\notag
\\
&\qquad\le
\sqrt{T}C(\norm{\phi}_{L^{\infty}(0,T;H^{s+1}_{\cpe})})\norm{\tvariii}_{L^{\infty}(0,T;H^{\sigma+1})}\notag
\\
&\qquad\le
\sqrt{T}C(\Omega)\norm{\tvariii}_{L^{\infty}(0,T;H^{\sigma+1})}\le
\widetilde{C}\norm{\tvariii}_{L^{\infty}(0,T;H^{\sigma+1})}.
\end{align*}
From this together with~\eqref{6.41}, we
conclude that 
$$
\norm{\Er_{2}(\phi,\nabla\phi)\nabla\tvariii}_{\Kr^{\sigma}_{\cpe}(T)}
\le \widetilde{C}\norm{\tvariii}_{L^{\infty}(0,T;H^{\sigma+1})}.
$$

This completes the proof of~\eqref{6.21} since
$\norm{\tvariii}_{L^{\infty}(0,T;H^{\sigma+1})}\le\|\tvariii\|_{\Kr^{\sigma+2}_{\cpe}(T)}$.
So to prove~\eqref{6.205} it remains only to show that 
\begin{equation}\label{desi:pevariii}
\norm{\pe\tvariii}_{\Kr^{\sigma}_{\cpe}(T)}\le {\rm RHS}.
%\widetilde{C}\norm{\tvariii(0)}_{\nhsc{\sigma+1}}
%+C(\Omega)\bigl\{\norm{f_{3}}_{L^{2}(0,T;H^{\sigma})}
%+\eps \norm{f_{1}}_{L^{2}(0,T;H^{\sigma})}\bigr\}.
\end{equation}
To see this, solve the first
equation in~\eqref{system:NS5} 
for~$\cn\vard$ and substitute 
the result in the third equation, to obtain
\begin{equation}\label{Eq:variiiLFinduc}
g_{3}(\phi)\partial_{t}\tvariii
-\pe (B_{3}(\phi)-B_{1}(\phi))\tvariii=f_{3}+f'_{3},
\end{equation}
with $f_{3}'\defn g_{1}(\phi)\edtt\tvari-\eps f_{1}$. Since 
$$
B_{3}(\phi)-B_{1}(\phi)\defn (\chi_{3}(\eps\vari)-\chi_{1}(\eps\vari))\cn(\beta(\variii)\nabla\cdot),
$$
the assumption $\chi_{1}(\wp)<\chi_{3}(\wp)$ (which is (A3) in
Assumption~\ref{assu:structural}) implies that the
equation~\eqref{Eq:variiiLFinduc} is parabolic. Hence, one can use the
following estimate (see~\eqref{esti:hypparusual} below):
\begin{equation}\label{esti:variiipar}
\begin{split}
\norm{\pe\tvariii}_{\Kr^{\sigma+2}_{\cpe}(T)}&\le\norm{\sqrt{\pe}\cpe\tvariii}_{L^{\infty}(0,T;H^{\sigma+1})}
+\sqrt{\pe}\norm{\sqrt{\pe}\cpe\tvariii}_{L^{2}(0,T;H^{\sigma+2})}\\
&\le \widetilde{C}\bigl\{ \norm{\sqrt{\pe}\cpe\tvariii(0)}_{\nh{\sigma+1}}
+\cpe\norm{f_{3}+f_{3}'}_{L^{2}(0,T;H^{\sigma})}\bigr\}.
\end{split}
\end{equation}
From which we easily infer the desired bound~\eqref{desi:pevariii}.

\smallbreak
\noindent {\textbf{c)}} {\em{Estimate for}} $\nabla\tvari$. 
Note that, by combining~\eqref{6:16'} with~\eqref{6.17}, one has
\begin{equation}\label{6.30}
\norm{\nabla\tvari}_{\Kr^{\sigma}_{\cpe}(T)}\le\widetilde{C}\norm{\cn(g_{2}(\phi)^{-1}\nabla\tvari)}_{\Kr^{\sigma-1}_{\cpe}(T)}
+\widetilde{C}\norm{\tvari}_{L^{\infty}(0,T;L^{2})}.
\end{equation}
Starting from 
$\nabla\tvari=-g_{2}(\phi)\edtt\tvard +\eps\r B_{2}(\phi)\tvard +\eps
f_{2}$, it is found that
\begin{equation}\label{6:31}
\begin{aligned}
&\cn(g_{2}(\phi)^{-1}\nabla\tvari)=-\edtt\cn\tvard +\eps\r
\Fr_{1}(\phi)\nabla^{2}\cn\tvard \\
&\quad +\eps\r
\Fr_{2}(\phi,\nabla\phi)\nabla^{2}\tvard+\eps\r\Fr_{3}(\phi,\nabla\phi,\nabla^{2}\phi)\nabla\tvard+\eps \cn(g_{2}(\phi)^{-1} f_{2}), 
\end{aligned}
\end{equation}
for some $C^{\infty}$ functions $\Fr_{1}$, $\Fr_{2}$ and $\Fr_{3}$ vanishing at
the origin.

The most direct estimates show that the
$\Kr^{\sigma-1}_{\cpe}(T)$-norms of the first and the last three terms
in the right-hand side of~\eqref{6:31} are 
bounded by RHS. So we need only concentrate on the second term. Since
$0\le\sigma-1\le s$, 
the estimate~\eqref{6:16'} implies
\begin{equation*}
\norm{\eps\r\Fr_{1}(\phi)\nabla^{2}\cn\tvard}_{\Kr^{\sigma-1}_{\cpe}(T)}
\le \widetilde{C}\norm{\eps\r\nabla^{2}\cn\tvard}_{\Kr^{\sigma-1}_{\cpe}(T)}.
\end{equation*}
The following observation is important since this where the spectral
localization of $\tvar$ enters. Since $\tvard=J_{\eps\cpe/3}\tvard$
and since $\eps\r \nabla^{2}J_{\eps\cpe/3}\les \Fi{}{1}$, one has 
$\norm{\eps\r\nabla^{2}\cn\tvard}_{\Kr^{\sigma-1}_{\cpe}(T)}\les 
\norm{\cn\tvard}_{\Kr^{\sigma}_{\cpe}(T)}$. The piece of information
already determined in~\eqref{6.205} implies that this in turn is $\le{\rm
  RHS}$. Hence, the
right-hand side of~\eqref{6.30} is $\le{\rm RHS}$. 
The proof of~\eqref{esti:LFinduction} is complete.
\end{proof}

Recall that the upshot of the two previous lemma is to estimate the
fast components $\cn J_{\eps\cpe}\vard$ and $\nabla J_{\eps\vard}\vari$.
\begin{lemma}There exists a generic function $C(\cdot)$ such that 
\begin{align}
\norm{J_{\eps\cpe}\vari}_{L^{\infty}(0,T;H^{s})}+\cpe\norm{J_{\eps\cpe}\vari}_{L^{2}(0,T;H^{s+1})}&\le \widetilde{C},\label{esti:vari}\\
\norm{\cn J_{\eps\cpe}\vard}_{L^{\infty}(0,T;H^{s-1})}+\cpe\norm{\cn J_{\eps\cpe}\vard}_{L^{2}(0,T;H^{s})}&\le \widetilde{C}.\label{esti:vard}
\end{align}
%\norm{\vari}_{\Kr^{s+1}_{\cpe}(T)}+\norm{\cn\vard}_{\Kr^{s}_{\cpe}(T)}\le \widetilde{C}.
%$$
\end{lemma}
\begin{proof}
To avoid repetitions, we just give the scheme of the analysis. Set 
$$
X_{\ik}\defn \norm{p_{\ik}}_{\Kr^{s-\ik+1}_{\cpe}(T)}
+\norm{\cn\vard_{\ik}}_{\Kr^{s-\ik}_{\cpe}(T)}.
$$
Setting $\tvar\defn\var_{\ik}$ and $f\defn f^{\ik}_{\LF}-G(\phi)\vard\cdot\nabla\var_{\ik}$, where $G(\phi)$
is as defined in~\eqref{Eq:terest}, we are in the situation of the
previous lemma and hence the estimate~\eqref{esti:LFinduction} easily
implies $X_{\ik}\le \widetilde{C}X_{\ik+1}+Y_{\ik}$ 
where
\begin{align*}
Y_{\ik}&\le \widetilde{C}\bigl\{\norm{U_{\ik}}_{\Hr^{0}_{\pa}(T)}+\norm{\var_{\ik}(0)}_{\Hr^{s-\ik}_{\pa,0}}\bigr\}\\
&\quad
+(\eps+\sqrt{T})C(\Omega)\bigl\{
\norm{f^{\ik}_{\LF}}_{L^{\infty}(0,T;\Hr^{s-1-\ik}_{\pa,0})}+\norm{\var_{\ik}}_{\Hr^{s-\ik}_{\pa}(T)}\bigr\}.
\end{align*}
Gathering the results of the previous lemma, one obtains
\begin{align*}
Y_{\ik}\le
\widetilde{C}\bigl\{\widetilde{C}+C(\Omega_{0})\bigr\}
+(\eps+\sqrt{T})C(\Omega)\bigl\{C(\Omega)+C(\Omega)\bigr\}
\le\widetilde{C}.
\end{align*}
Hence, we end up with 
$X_{\ik}\le \widetilde{C}X_{\ik+1}+\widetilde{C}$. 
By an elementary induction, we therefore obtain 
$X_{0}\le \widetilde{C}X_{s}+\widetilde{C}$. 
Finally, noting that $X_{s}\les
\norm{U_{s}}_{\Hr^{0}_{\pa}(T)}$ and using Lemma~\ref{lemm:induct1} 
leads to the desired bound $X_{0}\le
\widetilde{C}$. 
\end{proof}

\subsection{The slow components}\label{subsection:incompressible}
Having proved the estimates for the fast components, 
we now prove the estimates for the slow components 
$\variii$ and $\curl\vard$. 
To prove both estimates we use Assumption~\ref{assu:compatibility}. 
This furnishes us with a $C^{\infty}$ function $S=S(\vartheta,\wp)$
such that $S(0,0)=0$ and 
\begin{align}
&\diff S (\vartheta,\wp)= g_{3}(\vartheta,\wp)\diff \vartheta -
g_{1}(\vartheta,\wp)\diff\wp , \label{iden:comp1S}\\
&(\vartheta,\wp)\mapsto(S(\vartheta,\wp),\wp)\text{ is a }C^{\infty}
\text{ change of variables.}\label{iden:propCVS}
\end{align}
[Recall that $(\vartheta,\wp)\in\xR^{2}$ is the place holder of $(\variii,\eps\vari)$.]

\begin{nota} Set $\sigma\defn S(\variii,\eps\vari)$.
\end{nota}

One reason it is interesting to introduce the coordinate
$S=S(\vartheta,\wp)$ is that
$\sigma$ is well transported by the flow (see also
Remarks~\ref{rema:sigmawellflowed} and~\ref{rema:sigmainterest}).
\begin{lemma}Given $F\in C^{\infty}(\xR)$ satisfying $F(0)=0$, there
exists a generic function $C(\cdot)$ such that
\begin{align}
\norm{(\partial_{t}+\vard\cdot\nabla)F(\sigma)}_{L^{2}(0,T;H^{s})}\le C(\Omega),\label{esti:goodcurlsplit}\\
\norm{F(\sigma)}_{L^{\infty}(0,T;H^{s})}\le \widetilde{C}.\label{esti:goodcurlsplit2}
\end{align}
\end{lemma}
\begin{proof}
Firstly, we form an evolution equation for
$\sigma$. Directly from the identity~\eqref{iden:comp1S}, one has 
$\partial_{t,x}\sigma=g_{3}(\phi)\partial_{t,x}\variii-\eps g_{1}(\phi)\partial_{t,x}\vari$. 
By combining the first and the last
equations in~\eqref{system:NS2} with this identity, we compute
\begin{equation}\label{Eq:firstsigma}
\partial_{t}\sigma+\vard\cdot\nabla \sigma-\pe \bigl(\chi_{3}(\eps\vari)-\chi_{1}(\eps\vari)\bigr)
\cn(\beta(\variii)\nabla\variii)=0.
\end{equation}
Therefore, one has 
\begin{equation*}
\partial_{t}\sigma+\vard\cdot\nabla\sigma=\pe G_{1}(\phi,\nabla\variii)
+\pe G_{2}(\phi)\Delta\variii,
\end{equation*}
for some $C^{\infty}$ functions $G_{1}$ and $G_{2}$, with $G_{1}(0)=0$. 
Hence, the Moser estimate~\eqref{comp} implies at once
\begin{align*}
\norm{\partial_{t}\sigma+\vard\cdot\nabla\sigma}_{\nh{s}}\le 
C(\norm{\phi}_{\nh{s}})\bigl\{1+
\pe\norm{\nabla\variii}_{\nh{s+1}}\bigr\}.
\end{align*}
Since $\sqrt{\pe}\le\cpe$, one has 
$$
\pe\norm{\nabla\variii}_{\nh{s+2}}\le \sqrt{\pe}\norm{\nabla\variii}_{\nhsc{s+1}},
$$
so that $\norm{\partial_{t}\sigma+\vard\cdot\nabla\sigma}_{\nh{s}}\le
C(R)\{1+R'\}$ where $R$ and $R'$ are as defined in~\eqref{defi:RR'2}. 
Since $\Omega\approx\lA R\rA_{L^{\infty}(0,T)}+\lA R'\rA_{L^{2}(0,T)}$,
we have
\begin{equation}\label{esti:goodcurlsplit0}
\norm{(\partial_{t}+\vard\cdot\nabla)\sigma}_{L^{2}(0,T;H^{s})}\le C(\Omega).
\end{equation}
The chain rule and the rule of product in Sobolev spaces imply
that the left-hand side of~\eqref{esti:goodcurlsplit} is estimated by 
\begin{equation}\label{RHS6s}
K(1+\norm{\widetilde{F'}(\sigma)}_{L^{\infty}(0,T;H^{s})})\norm{(\partial_{t}+\vard\cdot\nabla)\sigma}_{L^{2}(0,T;H^{s})}.
\end{equation}
where $\widetilde{F'}(z)=F'(z)-F'(0)$ where $F'$ denotes the differential of $F$. 
The first term in~\eqref{RHS6s} is bounded by 
$C(\norm{\sigma}_{L^{\infty}(0,T;H^{s})})\le C(\Omega)$ (here is where
we use the hypothesis that $S(0,0)=0$). Hence,
~\eqref{esti:goodcurlsplit} follows from~\eqref{esti:goodcurlsplit0}.

\smallbreak
Moving to the proof of~\eqref{esti:goodcurlsplit2}, we first 
recall an usual estimate in Sobolev spaces for hyperbolic
equations\footnote{See also the estimate~\eqref{esti:hypparusual} with $\eta=0$.}: For all functions $u,\chi \in C^{0}([0,T];H^{s}(\xD))$, 
with $s>1+d/2$, there holds
\begin{equation}
\norm{u}_{L^{\infty}(0,T;H^{s})}\les
e^{TX}\norm{u(0)}_{\nh{s}}+\int_{0}^{T}e^{(T-t)X}\norm{\partial_{t}u+\chi\cdot\nabla u}_{\nh{s}}\,dt,
\end{equation}
where $X\defn K\int_{0}^{T}\norm{\chi}_{\nh{s}}\,dt$ for some constant
$K$ depending only on $(s,d)$. The Cauchy-Schwarz inequality readily
implies that the right-hand side is estimated by
$$
e^{KT\norm{\chi}_{L^{\infty}(0,T;H^{s})}}\bigl\{\norm{u(0)}_{\nh{s}}+
\sqrt{T}\norm{\partial_{t}u+\chi\cdot\nabla u}_{L^{2}(0,T;H^{s})} \bigr\}.
$$
Applying this bound with $(u,\chi)=(F(\sigma),\vard)$, it is found
that the $L^{\infty}(0,T;H^{s})$ norm of $F(\sigma)$ is estimated by
$$
e^{KT\Omega}\bigl\{C(\Omega_{0})+\sqrt{T}\norm{\partial_{t}F(\sigma)+\vard\cdot\nabla
  F(\sigma)}_{L^{2}(0,T;H^{s})}\bigr\}.
$$
We complete the proof of~\eqref{esti:goodcurlsplit2} by using~\eqref{esti:goodcurlsplit}.
\end{proof}

\begin{corollary}\label{coro:estitheta}
Given $F\in C^{\infty}(\xR^{2+d})$ satisfying $F(0)=0$, 
there exists a generic function $C(\cdot)$ such that
\begin{equation}\label{esti:goodcurlsplit3}
\norm{F(\ter)}_{L^{\infty}(0,T;H^{s})}\le \widetilde{C}.
\end{equation}
\end{corollary}
\begin{proof}
The property~\eqref{iden:propCVS} implies that there exists a $C^{\infty}$
function $F^{*}$ such that $F(\vartheta,\wp,\mathsf{v})=F^{*}(S(\vartheta,\wp),\wp,\mathsf{v})$ for all
$(\vartheta,\wp,\mathsf{v})\in\xR\times\xR\times\xR^{d}$. Moreover, the hypothesis $S(0,0)=0$ implies
$F^{*}(0,0,0)=0$. 

Decompose $F(\ter)\defn F(\variii,\eps\vari,\eps\vard)=F^{*}(\sigma,\eps\vari,\eps\vard)$ as 
$$
F^{*}(\sigma,0,0)+\bigl\{F^{*}(\sigma,\eps\vari,\eps\vard)-F^{*}(\sigma,0,0)\bigr\}.$$
The first 
term is estimated by way of~\eqref{esti:goodcurlsplit2}. To bound the
second term, observe that one can factor out $(\eps\vari,\eps\vard)$. Consequently, there exists a
function 
$C(\cdot)$ depending only on $F^{*}$ such that
$$
\norm{F^{*}(\sigma,\eps\vari,\eps\vard)-F^{*}(\sigma,0)}_{\nh{s}}\le \eps
C(\norm{(\variii,\vari,\vard)}_{\nh{s}})\le \eps C(R)\le e^{\eps C(\norm{R}_{L^{\infty}(0,T)})}.
$$
Which completes the proof.
\end{proof}

\begin{remark}\label{rema:sigmawellflowed}
The estimate~\eqref{mlemm:estdescoef} implies
$\norm{\partial_{t}F(\ter)}_{L^{2}(0,T;H^{s}_{\cpe})}\le C(\Omega)$. 
As in~\eqref{boundLinfty}, it yields $\norm{F(\ter)}_{L^{\infty}(0,T;H^{s}_{\cpe})}\le
C(\Omega_{0})+\sqrt{T}C(\Omega)\le \widetilde{C}$. Yet, this is weaker than~\eqref{esti:goodcurlsplit3}
[indeed $H^{s}_{0}=H^{s-1}$].
\end{remark}

Now we really use the special feature of the low frenquency analysis. 
The following lemma states that
$\Fi{\eps\cpe}{-1}\sigma$ satisfies parabolic-type estimates.

\begin{lemma}\label{lemm:sigmaLF}
There exists a generic function $C(\cdot)$ such that
\begin{equation*}
\norm{\Fi{\eps\cpe}{-1}\sigma}_{L^{\infty}(0,T;H^{s+1}_{\cpe})}+
\sqrt{\pe}\norm{\Fi{\eps\cpe}{-1}\sigma}_{L^{2}(0,T;H^{s+2}_{\cpe})}\le \widetilde{C}.
\end{equation*}
\end{lemma}
\begin{proof}
In light of the estimate $\norm{\sigma}_{L^{\infty}(0,T;H^{s})}\le \widetilde{C}$
  (see~\eqref{esti:goodcurlsplit2}), it suffices to prove 
that $\dot{\sigma}\defn \cpe \Fi{\eps\cpe}{-1}\nabla\sigma$ satisfies
\begin{equation}\label{todo:dotsigma}
\norm{\dot{\sigma}}_{L^{\infty}(0,T;H^{s})}+
\sqrt{\pe}\norm{\dot{\sigma}}_{L^{2}(0,T;H^{s+1})}\le \widetilde{C}.
\end{equation}

Let us form a parabolic evolution
equation for $\dot{\sigma}\defn \cpe \Fi{\eps\cpe}{-1}\nabla\sigma$. 
Writing the identity~\eqref{iden:comp1S} in the form 
$
\diff\vartheta
=c_{1}(\vartheta,\wp)\diff
S+c_{2}(\vartheta,\wp)\diff\wp$ with 
$c_{1}\defn 1/g_{3}^{-1}$ and $c_{2}\defn g_{1}/g_{3}$, 
yields 
\begin{equation}\label{iden:variiisigmavari}
\nabla\variii=c_{1}(\phi)\nabla \sigma+ \eps
c_{2}(\phi)\nabla\vari.
\end{equation}
Inserting
this expression for $\nabla\variii$ into the
equation~\eqref{Eq:firstsigma}, yields
\begin{equation}\label{Eq:sigmaend}
\partial_{t}\sigma+\vard\cdot\nabla\sigma-\pe k(\phi)\Delta\sigma
=\pe \mathsf{G}_{3}+\pe\mathsf{G}_{4},
\end{equation}
where $k\defn(\chi_{3}-\chi_{1})\beta c_{1}$ is a smooth positive
function\footnote{By hypotheses, one has $\chi_{1}<\chi_{3}$,
  $0<c_{1}\defn 1/g_{3}$
  and $0<\beta$ (see Assumption~\ref{assu:structural}).} and 
$$
\mathsf{G}_{3}+\mathsf{G}_{4}\defn
G_{3}(\phi,\nabla\phi)
+\eps G_{4}(\phi,\nabla\phi)\Delta\vari,
$$
for some $C^{\infty}$ functions $G_{3}$ and $G_{4}$, with
$G_{3}(0)=0$. 

We compute
\begin{equation}\label{Eq:sigmaend2}
\partial_{t}\dot{\sigma}+\vard\cdot\nabla\dot{\sigma}-\pe k(\phi)\Delta\dot{\sigma}=\Gr,
\end{equation}
where the source term is given by
\begin{align*}
\Gr&\defn -\cpe
\Fi{\eps\cpe}{-1}(\nabla\vard\cdot\nabla\sigma)+\cpe[\vard,\Fi{\eps\cpe}{-1}\bigr]\cdot\nabla\nabla\sigma\\
&\quad +\pe \Fi{\eps\cpe}{-1}(\nabla k(\phi) \Delta\sigma)+\pe\cpe [k(\phi),\Fi{\eps\cpe}{-1}]\Delta\nabla\sigma\\
&\quad +\cpe\pe \Fi{\eps\cpe}{-1}\nabla\mathsf{G}_{3}+ \cpe\pe \Fi{\eps\cpe}{-1}\nabla\mathsf{G}_{4}.
\end{align*}

Let us give the scheme of the analysis. 
In light of the standard estimate ~\eqref{esti:hypparusual} (given in
the appendix below), to prove~\eqref{todo:dotsigma} 
we need only show 
that the source term $\Gr$ can be split as 
$\Gr_{1}+\sqrt{\pe}\Gr_{2}$ with
$$
\norm{\Gr_{1}}_{L^{1}(0,T;H^{s})}\le \widetilde{C}\quad\mbox{and}\quad
\norm{\Gr_{2}}_{L^{2}(0,T;H^{s-1})}\le \widetilde{C}.
$$
To see this, there are only two nontrivial points. 
Firstly, note that
\begin{align*}
&\norm{\cpe\sqrt{\pe}
  \Fi{\eps\cpe}{-1}\nabla\mathsf{G}_{4}}_{L^{2}(0,T;H^{s-1})}\\
&\quad \defn \norm{\cpe\sqrt{\pe}
  \Fi{\eps\cpe}{-1}\nabla \{\eps
  G_{4}(\phi,\nabla\phi)\Delta\vari\}}_{L^{2}(0,T;H^{s-1})} &\\
&\quad \les \norm{\sqrt{\pe}
    G_{4}(\phi,\nabla\phi)\Delta\vari}_{L^{2}(0,T;H^{s-1})}
  &(\text{since } \eps\cpe\nabla \Fi{\eps\cpe}{-1}\les\id)\\
&\quad \les
(1+\norm{\tilde{G_{4}}(\phi,\nabla\phi)}_{L^{\infty}(0,T;H^{s-1})})
\norm{\sqrt{\pe}\vari}_{L^{2}(0,T;H^{s+1})} &(\text{straightforward})\\
&\quad \les \widetilde{C}\widetilde{C}. &(\text{by
}\eqref{esti:goodcurlsplit3}, \eqref{esti:vari})
\end{align*}
Secondly, there is a little trick: set 
$I_{\r,\pe}\defn\cpe/(\sqrt{\r}+\sqrt{\pe})\le 1$ and decompose
$\cpe\nabla\vard\cdot\nabla\sigma$ into two parts:
$$
\cpe\nabla\vard\cdot\nabla\sigma=I_{\r,\pe}\sqrt{\r}\nabla\vard\cdot\nabla\sigma+I_{\r,\pe}\sqrt{\pe}\nabla\vard\cdot\nabla\sigma.
$$ 
Using the
usual Moser estimates as well as the weighted versions
\eqref{esti:weightedprod}-\eqref{esti:weightedcomp}, it is
easily found that
\begin{align*}
\sqrt{\r}\norm{\nabla\vard\cdot\nabla\sigma}_{H^{s}}
&\les
\norm{\nabla\vard}_{H^{s}_{\sqrt{\r}}}\norm{\nabla\sigma}_{H^{s}_{\sqrt{\r}}}\les (R+R')C(R),
\\
\norm{\nabla\vard\cdot\nabla\sigma}_{H^{s-1}}
&\les \norm{\vard}_{H^{s}}\norm{\sigma}_{H^{s}}\les R C(R),
\end{align*}
so that
$$
\sqrt{\r}\norm{\nabla\vard\cdot\nabla\sigma}_{L^{1}(0,T;H^{s})}
+\norm{\nabla\vard\cdot\nabla\sigma}_{L^{2}(0,T;H^{s-1})}\le
\sqrt{T}C(\Omega)\le \widetilde{C}.
$$

With these estimates established, the proof easily follows.
%\qed
\end{proof}

\begin{remark}\label{rema:sigmainterest}
One interesting feature of the equation~\eqref{Eq:sigmaend} for 
$\sigma$ is that it is coupled to the momentum
equation only through the convective term. Indeed, for the purpose of
proving estimates independent of $\pe$, we cannot see the term
$\cn\vard$ [in the equation for $\variii$] as a source term.
\end{remark}

We denote by $\curl \vard$ the matrix with coefficients
$(\curl\vard)_{ij}\defn \partial_{j}\vard_{i}-\partial_{i}\vard_{j}$. 
The basic idea of the forthcoming computations is to apply the $\curl$
operator to the equation for $\vard$ so as to cancel the large term
$\eps^{-1}\nabla\vari$. 
Yet, this requires some preparation because 
the factor $g_{2}(\variii,\eps\vari)$ multiplying the time derivative
of~$\vard$ admits large oscillations in $\Or(1)$. 
To get around this, we follow the analysis of~\cite{MS1}. 
Namely, as in the proof of Corollary~\ref{coro:estitheta}, we decompose
$g_{2}(\variii,\eps\vari)$ 
into two parts: the first which is well transported by the flow, 
the second which admits small oscillations of typical size
$\Or(\eps)$. 

In particular, we do not estimate $\curl\vard$ directly. 
Instead we estimate $\curl(\gamma_{0}\vard)$ where the coefficient $\gamma_{0}$ is
defined as follows.

\begin{nota}
By~\eqref{iden:propCVS}, one can write $g_{2}(\vartheta,\wp)=
\Gamma\bigl(S(\vartheta,\wp),\wp\bigr)$ for some smooth positive function
$\Gamma$. Set
$\Gamma_{0}(\vartheta,\wp)\defn\Gamma\bigl(S(\vartheta,\wp),0\bigr)$
and 
\begin{equation}\label{defi:curlgamma0}
\gamma_{0}\defn\Gamma_{0}(\variii,\eps\vari).
\end{equation}
\end{nota}

\begin{lemma}\label{lemm:curl}
There exists a generic function $C(\cdot)$ such that
\begin{equation}\label{prop:slowcomponent:3}
\norm{\curl(\gamma_{0} \vard)}_{L^{\infty}(0,T;H^{s-1})}
+\sqrt{\r}\norm{\curl(\gamma_{0} \vard)}_{L^{2}(0,T;H^{s})}
\le \widetilde{C}.
\end{equation}
\end{lemma}
\begin{proof}
We begin by computing the equation satisfied by
$\omega\defn\curl(\gamma_{0}\vard)$. It follows from elementary calculus that there 
exists a $C^{\infty}$ function $\Gamma_{1}$ such that for all~$(\wp,\vartheta)\in\xR^{2}$
\begin{equation}\label{split:g2}
\Gamma_{0}(\vartheta,\wp)/g_{2}(\vartheta,\wp)=1+\wp\Gamma_{1}(\vartheta,\wp).
\end{equation}
We first insert the expression for $g_{2}(\variii,\eps\vari)$ given by~\eqref{split:g2} 
into the equation for~$\vard$, thereby obtaining
$$
\gamma_{0}(\partial_{t}\vard+\vard\cdot\nabla\vard)+\eps^{-1}\nabla\vari
-\r B_{2}\vard=F^{\dagger},
$$
with $F^{\dagger}\defn 
-\vari\gamma_{1}\nabla\vari+\eps\r\vari\gamma_{1} B_{2}\vard$ where
$\gamma_{1}\defn\Gamma_{1}(\variii,\eps\vari)$. 
Consequently, the equation for $\vard$ is equivalent to
$$
(\partial_{t}+\vard\cdot\nabla)(\gamma_{0}\vard) +\eps^{-1}\nabla\vari
-\r B_{2}\vard
=F^{\dagger}+(\partial_{t}\gamma_{0}+\vard\cdot\nabla\gamma_{0})\vard.
$$
Using the elementary identity $\curl\nabla=0$, we find
$$
(\partial_{t}+\vard\cdot\nabla\bigr)\curl (\gamma_{0}\vard)-
\r B_{2}\curl\vard=F^{\dagger\dagger},
$$
where
$$
F^{\dagger\dagger}\defn\curl F^{\dagger}
-\bigl[\curl,\vard\bigr]\cdot\nabla(\gamma_{0}\vard)+\r\bigr[\curl,B_{2}\bigr]\vard
+\curl\bigl(\vard(\partial_{t}\gamma_{0}+\vard\cdot\nabla\gamma_{0})\bigr).
$$
Next, write 
$B_{2}\curl=\gamma_{0}^{-1}\gamma_{0}B_{2}\curl=\gamma_{0}^{-1}B_{2}\curl(\gamma_{0}\cdot)
+\gamma_{0}^{-1}\bigl[\gamma_{0},B_{2}\curl\bigr]$ 
to obtain that $\omega=\curl(\gamma_{0}\vard)$ satisfies
\begin{equation}\label{Eq:curlLF}
\partial_{t}\omega+\vard\cdot\nabla\omega-\r
\gamma_{0}^{-1}B_{2}\omega= F,
\end{equation}
with $F\defn
F^{\dagger\dagger}+\r\gamma_{0}^{-1}[\gamma_{0},B_{2}\curl]\vard$. 

To sum up, 
$F=\sum_{1\le i\le } F_{i}$ with
\begin{alignat*}{4}
  F_{1}&\defn -\curl(\vari\gamma_{1}\nabla\vari), 
&&\qquad && F_{4} &&\defn\r[\curl,B_{2}]\vard, \\
  F_{2}&\defn\eps\r\curl(\vari\gamma_{1}B_{2}\vard),
&& \qquad &&F_{5}&&\defn-[\curl,\vard]\cdot\nabla(\gamma_{0}\vard),\\
 F_{3}&\defn\curl\bigl(\vard(\partial_{t}\gamma_{0}+\vard\cdot\nabla\gamma_{0})\bigr),
&& \qquad && F_{6}&&\defn\r\gamma_{0}^{-1}[\gamma_{0},B_{2}\curl]\vard.
\end{alignat*}

\smallbreak
\noindent
{\em Estimate for \/}$F$. 
As in the proof of Lemma~\ref{lemm:sigmaLF}, to
prove~\eqref{prop:slowcomponent:3} it suffices to show that 
one can decompose $F$ as $f_{1}+\sqrt{\r}f_{2}$ with 
\begin{equation}\label{claim:curlLF0}
\norm{f_{1}}_{L^{1}(0,T;H^{s-1})}\le \widetilde{C}\quad\mbox{and}\quad 
\norm{f_{2}}_{L^{2}(0,T;H^{s-2})}\le \widetilde{C}.
\end{equation}
To do so we decompose $F$ as $f_{1}+\sqrt{\r}f_{2}$ with 
\begin{equation}\label{claim:curlLF}
\norm{f_{1}}_{L^{2}(0,T;H^{s-1})}\le C(\Omega)\quad\mbox{and}\quad 
\norm{f_{2}}_{L^{\infty}(0,T;H^{s-2})}\le C(\Omega).
\end{equation}
[Note that~\eqref{claim:curlLF} implies~\eqref{claim:curlLF0} since $\sqrt{T}C(\Omega)\le
\widetilde{C}$.]

Set $f_{2}\defn \eps\sqrt{\r}
p\gamma_{1}\curl(B_{2}\vard) -\sqrt{\r}\gamma_{0}^{-1}\vard\wedge
B_{2}\nabla\gamma_{0}$ and $f_{1}\defn F-f_{2}$. It follows from the
very definition of $\Omega$ that $f_{2}$
satisfies the second estimate in~\eqref{claim:curlLF}. 
As regards $f_{1}$, the key point is the following: Starting from
$$
\norm{F_{3}}_{\nh{s-1}}\les \norm{\vard}_{\nh{s}}\norm{\partial_{t}\gamma_{0}+\vard\cdot\nabla\gamma_{0}}_{\nh{s}},
$$
we obtain 
$\norm{F_{3}}_{L^{2}(0,T;H^{s-1})}\les
\norm{\vard}_{L^{\infty}(0,T;H^{s})}\norm{\partial_{t}\gamma_{0}+\vard\cdot\nabla\gamma_{0}}_{L^{2}(0,T;H^{s})}$.
As a consequence, using the very definition of $\gamma_{0}$
(see~\eqref{defi:curlgamma0}), the estimate ~\eqref{esti:goodcurlsplit} implies that
$\norm{F_{3}}_{L^{2}(0,T;H^{s-1})}\le C(\Omega)$.

\smallbreak
Let us estimate the other terms without repeated uses of the rules of
product. 
To do so set $P\defn (\vari,\nabla\vari)$, $V\defn
(\vard,\nabla\vard)$ and 
$\Theta\defn (\variii,\nabla\variii)$.
Direct computations show that one can write $F_{1}+F_{2}+F_{4}+F_{5}+F_{6}-f_{2}$ as
\begin{equation}\label{form:curlLF}
\sum(\sqrt{\r}\partial_{x})^{\alpha}G_{\alpha,\pa}(P,V,\Theta) \qquad 
(\alpha\in\xN^{d},~\la\alpha\ra\le 1),
\end{equation}
for some 
$C^{\infty}$ functions $G_{\alpha,\pa}$ such that 
the family $\{ G_{\alpha,\pa}\,\arrowvert\,\pa\in\PA,\la\alpha\ra\le 1\}$ is bounded in
$C^{\infty}$, with in addition $G_{0,\pa}(0)=0$. The
Moser estimate~\eqref{esti:weightcomppre} implies
\begin{align*}
\norm{(\sqrt{\r}\partial_{x})^{\alpha}
  G_{\alpha,\pa}(P,V,\Theta)}_{\nh{s-1}}&\le 
\norm{G_{\alpha,\pa}(P,V,\Theta)}_{H^{s}_{\sqrt{\r}}}\\
&\le C(\norm{(P,V,\Theta)}_{L^{\infty}})\norm{(P,V,\Theta)}_{H^{s}_{\sqrt{\r}}}.
\end{align*}
Hence, directly from the definition of $\Omega$, we arise at the
desired bound
$$
\norm{(\sqrt{\r}\partial_{x})^{\alpha}
  G_{\alpha,\pa}(P,V,\Theta)}_{L^{2}(0,T;H^{s-1})}\le C(\Omega) .
$$
Which completes the proof of the claim~\eqref{claim:curlLF}
%\smallbreak
%\noindent {\em Estimate for \/}$\omega$. 
%Let us rewrite~\eqref{Eq:curlLF} in short form
%$L\omega=f_{1}+\sqrt{\r}f_{2}$. Commuting this equation
%with~$\Fi{}{s-1}$ we obtain that
%$\widetilde{\omega}\defn\Fi{}{s-1}\omega$ solves
%$$
%L\widetilde{\omega}=[L,\Fi{}{s-1}]\omega+\Fi{}{s-1}f_{1}+\sqrt{\r}\Fi{}{s-1}f_{2}.
%$$ 
%Multiply this equation by $\widetilde{\omega}$ and integrate over
%$[0,T]\times\xD$. The derivatives of the coefficients involved
%$\nabla\vard$, $\nabla\variii$ and $\nabla\vari$ which are estimated
%in $L^{\infty}$ by $C(\Omega)$. Using the Gronwall's lemma, we obtain
%\begin{align*}
%&\norm{\widetilde{\omega}}_{L^{\infty}(0,T;L^2)}
%+\sqrt{\r}\norm{\widetilde{\omega}}_{L^2(0,T;H^1)}\\
%&\qquad\qquad\le \norm{\widetilde{\omega}(0)}_{\nhz}
%+C(\Omega)\norm{[L,\Fi{}{s-1}]\omega}_{L^{1}(0,T;L^2)}\\
%&\qquad\qquad\quad+C(\Omega)\norm{\Fi{}{s-1}f_{1}}_{L^{1}(0,T;L^2)}
%+C(\Omega)\norm{\Fi{}{s-1}f_{2}}_{L^2(0,T;H^{-1})},
%\end{align*}
%see also~\eqref{esti:hypparusual}. Using~\eqref{claim:curlLF}, the Cauchy Schwarz estimate imply that the
%last two terms in the right-hand side are estimated by
%$\sqrt{T}C(\Omega)$. It remains to prove a similar estimate for
%$[L,\Fi{}{s-1}]\omega$. Firstly, 
%using the commutator estimate~\eqref{commutator:estimateusual}, it 
%is found that $\norm{[L,\Fi{}{s-1}]\omega}_{L^2}\le C(R)\{1+R'\}$, as 
%the reader can verify, yielding the bound~$\norm{[L,\Fi{}{s-1}]\omega}_{L^{1}(0,T;L^2)}\le
%\sqrt{T} C(\Omega)$. This completes the proof of Lemma~\ref{lemm:curl}.
%\qed
\end{proof}

\subsection{End of the proof of Proposition~\ref{prop:LFuniform}}
In order to prove Proposition~\ref{prop:LFuniform} it remains to
estimate $J_{\eps\cpe}\variii$ and $J_{\eps\cpe}\vard$.

\begin{lemma}\label{lemm:antelast}
Given $F\in C^{\infty}(\xR^{2})$ such that $F(0)=0$, 
there exists a generic function $C(\cdot)$ such that
$$
\norm{J_{\eps\cpe}F(\phi)}_{L^{\infty}(0,T;H^{s+1}_{\cpe})}+\sqrt{\pe}\norm{J_{\eps\cpe}\nabla
  F(\phi)}_{L^{2}(0,T;H^{s+1}_{\cpe})}
\le\widetilde{C}.
$$
\end{lemma}
\begin{remark}
With $F(\phi)=\variii$ this is the expected bound for $J_{\eps\cpe}\variii$.
\end{remark}
\begin{proof}
In light of~\eqref{esti:goodcurlsplit2}, it suffices to prove
\begin{equation}\label{todo:antelast}
\cpe\norm{J_{\eps\cpe}\nabla F(\phi)}_{L^{\infty}(0,T;H^{s})}+\cpe\sqrt{\pe}\norm{J_{\eps\cpe}\nabla^{2}F(\phi)}_{L^{2}(0,T;H^{s})}
\le\widetilde{C}.
\end{equation}

By~\eqref{iden:propCVS} one has 
\begin{equation}\label{6.5.1}
\nabla F(\phi)= F_{1}(\phi)\nabla\sigma+ \eps F_{2}(\phi)\nabla\vari,
\end{equation}
for some $C^{\infty}$ functions $F_{1}$ and $F_{2}$. 
Applying Proposition~\ref{prop:Product} with
$\index=s$, $\indexg_{1}=\indexg_{2}=m_{1}=0$ and $m_{2}=1$ leads to
\begin{align*}
\norm{\cpe\nabla \Fi{\eps\cpe}{-1}F(\phi)}_{\nh{s}}&\les \bigl(1+\norm{\tilde{F_{1}}(\phi)}_{\nh{s}}\bigr)\norm{\nabla
\cpe\Fi{\eps\cpe}{-1}\sigma}_{\nh{s}}\\
&\quad+\bigl(1+\norm{\tilde{F_{2}}(\phi)}_{\nh{s}}\bigr)\norm{\nabla
\eps\cpe \Fi{\eps\cpe}{-1}\vari}_{\nh{s}},
\end{align*}
where $\tilde{F}\defn F-F(0)$. Using $\cpe\norm{\cdot}_{\nh{s}}\le \norm{\cdot}_{\nhsc{s}}$,
$\eps\cpe\nabla\Fi{\eps\cpe}{-1}\les \id$
and~\eqref{esti:goodcurlsplit2}, there results
$$
\norm{\cpe\nabla\Fi{\eps\cpe}{-1}F(\phi)}_{L^{\infty}(0,T;H^{s})}\le
  \widetilde{C}\bigl\{\norm{\Fi{\eps\cpe}{-1}\sigma}_{L^{\infty}(0,T;H^{s+1}_{\cpe})}
+\norm{\vari}_{L^{\infty}(0,T;H^{s})}\bigr\}.
$$
As a consequence, from Lemma~\ref{lemm:sigmaLF} and~\eqref{esti:vari}, we get the first
half of~\eqref{todo:antelast}, namely: $\norm{\cpe J_{\eps\cpe}\nabla F(\phi)}_{L^{\infty}(0,T;H^{s})}\le
\widetilde{C}$. 

The technique for obtaining the second half is
similar. Differentiate~\eqref{6.5.1}, to obtain: 
$\partial_{i}\partial_{j}F(\phi)\defn\Fr_{I}+\Fr_{II}$ with 
\begin{align*}
\Fr_{I}&\defn  F_{1}(\phi)\partial_{i}\partial_{j}\sigma
+\eps F_{2}(\phi)\partial_{i}\partial_{j}\vari,\\
\Fr_{II}&\defn \partial_{i}F_{1}(\phi)\partial_{j}\sigma+\eps\partial_{i}F_{2}(\phi)\partial_{j}\vari.
\end{align*}
Exactly as above, one has
\begin{align*}
\cpe\sqrt{\pe}\norm{J_{\eps\cpe}\Fr_{I}}_{L^{2}(0,T;H^{s+2})}&\les
\widetilde{C}\bigl\{\cpe\sqrt{\pe}\norm{\Fi{\eps\cpe}{-1}\sigma}_{L^{2}(0,T;H^{s+2})}
+\sqrt{\pe}\norm{\vari}_{L^{2}(0,T;H^{s+1})}\bigr\}\\
&\le \widetilde{C}\widetilde{C}.
\end{align*}
In order to estimate $\Fr_{II}$, note that $\Fr_{II}\defn
G(\phi,\nabla\phi)$ for some $C^{\infty}$ function $G$ such that
$G(0)=0$. Hence, one has
\begin{align*}
\cpe\sqrt{\pe}\norm{\Fr_{II}}_{L^{2}(0,T;H^{s})}&\les
\sqrt{T}\norm{G(\phi,\nabla\phi)}_{L^{\infty}(0,T,H^{s}_{\cpe})} &&(\text{straightforward})\\
&\les \sqrt{T}C\bigl(\norm{\phi}_{L^{\infty}(0,T;H^{s+1}_{\cpe})}\bigr)
&&(\text{by }\eqref{esti:weightedcomp})\\
&\les\sqrt{T}C(\Omega)\le \widetilde{C} &&(\text{by definition}).
\end{align*}
This completes the proof.
\end{proof}

The next estimate finishes the proof of Proposition~\ref{prop:LFuniform}.
\begin{lemma}
There exists a generic function $C(\cdot)$ such that
$$
\norm{J_{\eps\cpe}\vard}_{L^{\infty}(0,T;H^{s+1}_{\eps\cpe})}+\sqrt{\r}\norm{J_{\eps\cpe}\nabla\vard}_{L^{2}(0,T;H^{s+1}_{\eps\cpe})}
\le\widetilde{C}.
$$
\end{lemma}
\begin{proof}
We need only explain how to combine all the previous Lemma. 

\smallbreak
\noindent\textbf{a)} To shorten notations, given a function $f$ such that $f-\underline{f}\in
H^{\sigma}$ for some constant $\underline{f}$, we denote by
$\norm{f}_{\dot{H}^{\sigma}}$ the norm $\norm{f-\underline{f}}_{\nh{\sigma}}$. The
semi-norm $\norm{\cdot}_{\dot{H}^{\sigma}_{\cpe}}$ is defined in the
same way.

We claim that for all smooth positive function $\gamma$, there exists
a constant $K$ such that for all $(\eps,\cpe)$ and for all vector field $u\in H^{s}(\xD)$, one has
\begin{equation}\label{6.7.4}
\begin{split}
\norm{J_{\eps\cpe} u}_{\nh{s}}\le K M^{2s+1}\bigl\{&
\norm{\cn J_{\eps\cpe} u}_{\nh{s-1}}+\norm{\curl(\gamma
  u)}_{\nh{s-1}}\\
&+\norm{J_{\eps\cpe} u}_{\nhz}+\eps\norm{u}_{\nh{s}}\bigr\},
\end{split}
\end{equation}
where $ M\defn
1+\norm{\gamma^{-1}}_{\dot{H^{s}}}+\norm{\gamma}_{\dot{H}^{s+1}_{\cpe}}$.

To prove this result, we start from the following estimate: for all $\sigma\ge 0$,
\begin{equation}\label{6.7.1}
\norm{u}_{\nh{\sigma+1}}\les
\norm{\cn u}_{\nh{\sigma}}+\norm{\curl u}_{\nh{\sigma}}+\norm{u}_{\nhz}.
\end{equation}

As  $\curl u=\gamma^{-1}\curl(\gamma
u)-\gamma^{-1}\nabla\gamma\wedge u$, 
by the usual rule of product in Sobolev spaces, we infer that for all
$\sigma\in [0,s-1]$,
\begin{align*}
\norm{\curl u}_{\nh{\sigma}}&\les
(1+\norm{\gamma^{-1}}_{\dot{H}^{s}})\{
\norm{\curl(\gamma
  u)}_{\nh{\sigma}}+\norm{\gamma}_{\dot{H}^{s}}\norm{u}_{\nh{\sigma}}\} \\
&\les M \norm{\curl(\gamma u)}_{\nh{\sigma}}+M^{2}\norm{u}_{\nh{\sigma}}.
\end{align*}
Thus, the estimate~\eqref{6.7.1} turns into
\begin{equation*}%\label{6.7.2}
\norm{u}_{\nh{\sigma+1}}\les M^{2}\bigl\{
\norm{\cn u}_{\nh{\sigma}}+\norm{\curl(\gamma u)}_{\nh{\sigma}}+\norm{u}_{\nh{\sigma}}\bigr\}.
\end{equation*}
By induction on $\xN\ni\sigma \le s-1$, it yields
\begin{equation}\label{6.7.3}
\norm{u}_{\nh{s}}\les M^{2s}\bigl\{
\norm{\cn u}_{\nh{s-1}}+\norm{\curl(\gamma u)}_{\nh{s-1}}+\norm{u}_{\nhz}\bigr\}.
\end{equation}
Furthermore, by~\eqref{Friedrichs} applied
with $(m,\index,\indexg)=(1,s+1,s)$, one has
\begin{align*}
\norm{\curl(\gamma J_{\eps\cpe}u)}_{\nh{s-1}}&\le
\norm{J_{\eps\cpe}\curl(\gamma u)}_{\nh{s-1}}+\norm{[\gamma,J_{\eps\cpe}]u}_{\nh{s}}\\
&\les \norm{\curl(\gamma u)}_{\nh{\sigma}}+\eps\cpe\norm{\gamma}_{\dot{H}^{s+1}}\norm{u}_{\nh{s}}.
\end{align*}
From the estimate $\cpe\norm{\gamma}_{\dot{H}^{s+1}}\le
\norm{\gamma}_{\dot{H}^{s+1}_{\cpe}}\le M$ and \eqref{6.7.3} applied
with $u$ replaced by $J_{\eps\cpe} u$, the previous fact implies~\eqref{6.7.4}.

\smallbreak
\noindent\textbf{b)}
Let $\gamma\defn \gamma_{0}$ (as defined in~\eqref{defi:curlgamma0}) 
and let $M\defn
1+\norm{\gamma^{-1}}_{\dot{H^{s}}}+\norm{\gamma}_{\dot{H}^{s+1}_{\cpe}}$. 
We list our bounds:
\begin{align*}
&\text{by Lemma }\ref{lemm:antelast}: 
&&\norm{M}_{L^{\infty}(0,T)}\le\widetilde{C},\\
&\text{by }\eqref{esti:vard}: 
&&\norm{\cn J_{\eps\cpe}\vard}_{L^{\infty}(0,T;H^{s-1})}\le \widetilde{C},\\
&\text{by Lemma }\ref{lemm:curl}: 
&&\norm{\curl(\gamma_{0}\vard)}_{L^{\infty}(0,T;H^{s-1})}\le \widetilde{C},\\
&\text{by  Lemma }\ref{lemm:induct1}:
&&\norm{J_{\eps\cpe}\vard}_{L^{\infty}(0,T;L^{2})}\le\widetilde{C}, \\
&\text{directly from the definitions}:
&&\norm{\eps\vard}_{L^{\infty}(0,T;H^{s})}\le \eps
\Omega\le\widetilde{C}.
\end{align*}

Hence, we deduce from~\eqref{6.7.4} that 
$\norm{J_{\eps\cpe}\vard}_{L^{\infty}(0,T;H^{s})}\le
\widetilde{C}$. This in turn implies 
$\norm{J_{\eps\cpe}\vard}_{L^{\infty}(0,T;H^{s+1}_{\eps\cpe})}\le
\widetilde{C}$, indeed recall that 
$\norm{J_{\eps\cpe}\vard}_{H^{\sigma+1}_{\eps\cpe}}\les \norm{J_{\eps\cpe}\vard}_{\nh{\sigma}}$.

\smallbreak
\noindent\textbf{c)}
By~\eqref{6.7.1}, \eqref{esti:vard} and
Lemma~\ref{lemm:induct1}, 
to prove
$\sqrt{\r}\norm{J_{\eps\cpe}\nabla\vard}_{L^{2}(0,T;H^{s+1}_{\eps\cpe})}\le\widetilde{C}$ 
we need only check that 
$\sqrt{\r}\norm{\curl J_{\eps\cpe}\vard}_{L^{2}(0,T;H^{s})}\le
\widetilde{C}$. For this purpose, write
$$
\norm{\curl J_{\eps\cpe}\vard}_{\nh{s}}\les
(1+\norm{\gamma_{0}^{-1}}_{\dot{H}^{s}})\{\norm{\curl(\gamma_{0} J_{\eps\cpe}\vard)}_{\nh{s}}
+\norm{\gamma_{0}}_{\dot{H}^{s+1}}\norm{\vard}_{\nh{s}}\} ,
$$
and
$$
\norm{\curl(\gamma_{0} J_{\eps\cpe}\vard)}_{\nh{s}}
\les \norm{\curl(\gamma_{0} \vard)}_{\nh{s}}+\norm{\gamma_{0}}_{\dot{H}^{s+1}}\norm{\vard}_{\nh{s}},
$$
where we used the commutator estimate~\eqref{Friedrichs} with
$(m,\index,\indexg)=(0,s+1,s)$. These two facts and the estimate  $\sqrt{\r}\norm{\gamma_{0}}_{L^{\infty}(0,T;\dot{H}^{s+1})}\le
\norm{\gamma_{0}}_{L^{\infty}(0,T;\dot{H}^{s+1}_{\cpe})}\le \widetilde{C}$
(see Lemma~\ref{lemm:antelast}) 
implies 
$$
\sqrt{\r}\norm{\curl J_{\eps\cpe}\vard}_{L^{2}(0,T;H^{s})}
\le
\widetilde{C}\bigl\{\sqrt{\r}\norm{\curl(\gamma_{0}\vard)}_{L^{2}(0,T;H^{s})}
+\norm{\vard}_{L^{2}(0,T;H^{s})}\bigr\}.
$$

Hence, the expected bound directly follows from Lemma~\ref{lemm:curl}
and the estimate 
$\norm{\vard}_{L^{2}(0,T;H^{s})}\le
\sqrt{T}\norm{\vard}_{L^{\infty}(0,T;H^{s})}\le \widetilde{C}$.
\end{proof}

\section{Uniform stability}\label{section:proofoftheouniform}

In this section we prove Theorem~\ref{theo:uniform}. 
Granted the uniform {\em a priori\/} bounds proved in the previous
sections, the proof essentially reduces to establishing a local
existence result for fixed $\pa\in\PA$. 

\subsection{Local existence result}\label{subsection:localexistence}
For fixed $\pa\in\PA$, the fact that the Cauchy problem
for~\eqref{system:NS} is well-posed is
immediate provided that one chooses to work with the unknown~$(\rho,\vard,\mathcal{T})$. 
Indeed, with this choice of dependent variables the system is a coupled
hyperbolic/parabolic system in symmetric form, so that the general
theory applies. We begin by proving that this remains valid for the system~\eqref{system:NS2}.

\smallbreak
Introduce first a definition clarifying the structure of 
the systems we produce. 
Consider a system of nonlinear equations
\begin{equation}\label{defi:MHPsystem}
A_{0}(u)\partial_{t}u +\sum_{\alpha\in\xN^{d},\, 1\le \la\alpha \ra\le 2 }
A_{\alpha}(u)\partial_{x}^{\alpha} u 
= F(u,\nabla u),
\end{equation}
where each $A_{\beta}$ ($0\le\la\beta\ra\le 2$) is a $n\times n$ matrix 
smooth in its arguments and furthermore symmetric:
$A_{\beta}=A_{\beta}^{t}$. We suppose $F$ is smooth in its arguments,
with values in~$\xR^{n}$.
\begin{definition}\label{defi:MHPtype}
The system~\eqref{defi:MHPsystem} is said to be of coupled hyperbolic/parabolic type provided
that there exists
$(n_1,n_2)\in\xN^{2}$, and a splitting of the unknowns
$u=(u^1,u^2)\in\xR^{n_1}\times\xR^{n_2}$ so that 
$ F(u,\nabla u)=(0 , F^{2}(u,\nabla u))^{t}$ and 
$$
A_{0}(u)=\begin{pmatrix} A_{0}^{11}(u) & 0\\0 &A_{0}^{22}(u)
\end{pmatrix},
\quad A_{\alpha}(u)=\begin{pmatrix} 0 &0
\\ 0 & A_{\alpha}^{22}(u)\end{pmatrix} \quad\mbox{for all}~~\la\alpha\ra=2,
$$
where in addition 
$$
\forall (u,\xi)\in\xR^{n}\times\xS^{d-1},\qquad 
A_{0}^{11}(u) > 0,\quad A_{0}^{22}(u) > 0,\quad 
-\sum_{\la\alpha\ra =2}A_{\alpha}^{22}(u)\xi^{\alpha} >0.
$$
(here $M^{j k}$ denotes the sub-blocks of the matrix $M$
which corresponds to the splitting $\xR^{n}\ni
u=(u^1,u^2)\in\xR^{n_{1}}\times\xR^{n_{2}}$).
\end{definition}

The general theory~\cite{Maj,Taylor} 
for hyperbolic or parabolic problems applies for coupled 
hyperbolic/parabolic ones. One can prove the following result.
\begin{proposition}\label{prop:MHPCauchy}
$\mathsf{P}1)$ Let $s>2+d/2$. 
Assume that~\eqref{defi:MHPsystem} is of coupled hyperbolic/parabolic type. 
Then for all initial data $u_{0}$ in $H^{s}(\xD)$ there exists 
a positive $T=T(\norm{u_{0}}_{\nh{s}})$ such that the Cauchy problem 
for~\eqref{defi:MHPsystem} has a unique classical solution 
$\Phi(u_{0})\in C^{0}([0,T];H^{s}(\xD))$ such that $\Phi^{2}(u_{0})\in
L^{2}(0,T;H^{s+1}(\xD))$ where $\Phi^{2}(u_{0})$ are the last $n_{2}$
components of $\Phi(u_{0})$  where $n_{2}$ is as in Definition~\ref{defi:MHPtype}.

$\mathsf{P}2)$ The interval $[0,T^{\star})$, with 
$T^{\star}<+\infty$, is a maximal interval of $H^{s}$ existence 
if and only if $\limsup_{t\rightarrow T^{\star}} \norm{\Phi(u_{0})(t)}_{W^{1,\infty}(\xD)} =+\infty$.
\end{proposition}

\begin{remark}\label{rema:smoothness}
\textbf{i.} If the system is symmetric hyperbolic, it is enough to
assume $s>1+d/2$.

\noindent \textbf{ii.} Let $u\in C^{0}([0,T];H^{s}(\xD))$ be 
a classical solution of~\eqref{defi:MHPsystem}. It follows from 
property $\mathsf{P}2)$ that,  if $u(0)\in H^{\infty}(\xD)$ then $u\in C^{\infty}([0,T];H^{\infty}(\xD))$.
\end{remark}

\smallbreak
Fix~$\pa=(\eps,\r,\pe)\in\PA$. We next show that an appropriate change of variables transforms the
system~\eqref{system:NS2} into a system of coupled hyperbolic/parabolic type. 
Let $\varrho$ denotes the function given by 
Assumption~\ref{assu:compatibility}. Since the
mapping $\xR^{2}\ni(\vartheta,\wp)\mapsto(\vartheta,\varrho(\vartheta,\wp))\in\xR^{2}$ is a
$C^{\infty}$ diffeomorphism, it is equivalent to work
with $(\vari,\vard,\variii)$ or 
with $(\rho,\vard,\variii)$ where $\rho\defn
\varrho(\variii,\eps\vari)$. 

Firstly, we form an evolution equation for
$\rho\defn\varrho(\variii,\eps\vari)$. By combining the first and the last
equations in~\eqref{system:NS2} with the second identity
in~\eqref{iden:comp} written in the form $\chi_{3}\diff
\varrho=\chi_{3}g_{1}\diff \wp-\chi_{1}g_{3}\diff\vartheta$, we get
\begin{equation}
\chi_{3}(\partial_{t}\rho+\vard\cdot\nabla\rho) + (\chi_{3}-\chi_{1})\cn\vard=0.
\end{equation}
The second identity
in~\eqref{iden:comp} written in the form 
$$
\diff \wp=\gamma_{1}\diff\vartheta+\gamma_{2}\diff \varrho\defn
\frac{\chi_{1}g_{3}}{\chi_{3}g_{1}}\diff\vartheta+\frac{1}{g_{1}}\diff\varrho,
$$
yields
$\eps\nabla\vari=\gamma_{1}\nabla\variii+\gamma_{2}\nabla\rho$. Summing up, one has
\begin{equation}\label{system:MHPa}
\left\{
\begin{aligned}
&\chi_{3}(\partial_{t}\rho+\vard\cdot\nabla\rho) + (\chi_{3}-\chi_{1})\cn\vard=0,\\
&g_{2}(\ffp{\vard})+\eps^{-2}\gamma_{1}\nabla\variii+\eps^{-2}\gamma_{2}\nabla\rho
-  \r B_{2}\vard  = 0,\\
&g_{3}(\ffp{\variii})+\cn\vard -  \pe B_{3}\variii=0,
\end{aligned}
\right.
\end{equation}
where all the coefficients are positive by Assumption~\ref{assu:structural}. We also mention that all the coefficients 
are evaluated at $(\variii,\eps\vari)=(\variii,\mathsf{P}(\variii,\rho))$
for some $C^{\infty}$ function $\mathsf{P}$. We next symmetrize the system~\eqref{system:MHPa}. To do so 
multiply the first
equation by $ (\chi_{3}-\chi_{1})^{-1} \gamma_{1}$, the second by
$\eps^{2}$ and the third one by $\gamma_{2}$, to obtain
\begin{equation}\label{system:MHPb}
\left\{
\begin{aligned}
&\delta_{1}(\partial_{t}\rho+\vard\cdot\nabla\rho) + \gamma_{1}\cn\vard=0,\\
&\delta_{2}(\ffp{\vard})+\gamma_{1}\nabla\rho+\gamma_{2}\nabla\variii - \eps^{2}\r B_{2}\vard  = 0,\\
&\delta_{3}(\ffp{\variii})+\gamma_{2}\cn\vard -  \pe \gamma_{2} B_{3}\variii=0,
\end{aligned}
\right.
\end{equation}
with $\delta_{1}\defn
(\chi_{3}-\chi_{1})^{-1}\gamma_{1}\chi_{3}$, $\delta_{2}\defn\eps^{2}g_{2}$ and $\delta_{3}\defn
\gamma_{2}g_{3}$ [compare with the system~\eqref{normalformUr}]. Again, by
Assumption~\ref{assu:structural}, the coefficients $\delta_{i}$'s are positive. 

\begin{lemma}
For all $(\r,\pe)\in [0,1]^{2}$, the symmetric
system~\eqref{system:MHPb} is of coupled hyperbolic/parabolic type.
\end{lemma}
\begin{proof}Set $n=2+d$. If $\r=\pe=0$ set $n_1=n$ and $u^1=u$. If $\r=0$ and $\pe>0$ set 
$n_{1}=1+d$ and $u^1=(\rho,\vard)$. If $\r>0$ and $\pe=0$ set 
$n_1=2$ and $u^1=(\rho,\variii)$. If $\r>0$ and $\pe>0$ set $n_1=1$ 
and $u^1=\rho$.
%\qed
\end{proof}

With this result in hands, we obtain that the Cauchy problem
for~\eqref{system:MHPb} 
is well-posed for all fixed $\pa\in\PA$, so is  the Cauchy problem
for~\eqref{system:NS2} since~\eqref{system:MHPb}
has been deduced from~\eqref{system:NS2} by using a $C^{\infty}$
diffeomorphism. 

We are now prepared to prove Theorem~\ref{theo:uniform}.

\subsection{Proof of Theorem~\ref{theo:uniform}}\label{subsection:endproof1}

Let $s>1+d/2$ and $M_{0}>0$. 
The previous analysis implies that, for all $\pa\in\PA$, there exists $T_{\pa}=T_{\pa}(M_{0})$ such that for
all initial data $U_{0}$ in the ball $B(\Hr^{s}_{\pa,0},M_{0})$, the
Cauchy problem for~\eqref{system:NS2} has a unique classical solution
$U_{\pa}$ in~$\Hr^{s}_{\pa}(T_{\pa})$. Denotes by $T_{\pa}^{\star}$ the 
maximal time of existence of such a classical solution. 
For all $t<T_{\pa}^{\star}$, set 
$\Omega_{a}(t)\defn \norm{U_{\pa}}_{\mathcal{H}^{s}_{\pa}(t)}$. Our 
task is to show that there exists $T>0$ and $M<+\infty$ such that
\begin{equation}\label{task}
\forall \pa \in\PA,\qquad T_{\pa}^{\star}\ge T \quad\mbox{and}\quad
\Omega_{\pa}(T)\le M. 
\end{equation}

Up to replacing the initial data $U_{0}$ by $J_{\delta}U_{0}$ and
letting $\delta$ goes to zero, we can assume that $U_{0}\in
H^{\infty}(\xD)$. In light of Remark~\ref{rema:smoothness}, it means
that we can assume~$U_{\pa}\in
C^{\infty}([0,T_{\pa}^{\star});H^{\infty}(\xD))$. 
By combining Propositions~\ref{prop:HFuniform} and~\ref{prop:LFuniform}, 
we know that there exists a smooth non-decreasing function $C(\cdot)$
such that for all $\pa\in\PA$ and all $t<\min \{T_{\pa}^{\star},1\}$,
there holds
\begin{equation}\label{ineq:final}
\Omega_{\pa}(t)\le C(\Omega_{\pa}(0))e^{(\sqrt{t}+\eps)C(\Omega_{\pa}(t))}.
\end{equation}
Choose f\/irst~$M_{1} > C(M_{0})$ and next~$T_{1}\le 1$ and $\eps_{1}\le 1$ such that 
\begin{equation}\label{defi:Mun}
      C(M_{0})e^{(\sqrt{T_{1}}+\eps_{1})C(M_{1})} < M_{1} .
\end{equation}

Let~$t< \min \{ T_{\pa}^{\star}, T_{1}\}$ and
$\pa=(\eps,\r,\pe)\in\PA$ be such
that $\eps\le\eps_{1}$. By combining the inequalities~\eqref{ineq:final}
and~\eqref{defi:Mun} with the hypothesis $\Omega_{\pa}(0)< M_{0}$, we
infer that $\Omega_{\pa}(t)\neq M_1$. Besides, we can assume without restriction 
that $M_{0} < M_{1}$, so that $\Omega_{\pa}(0)<M_{1}$. Since the function
$\Omega_{\pa}$ is continuous, we infer~$\Omega_{\pa}(t) < M_{1}$. 
Consequently, the continuation principle [which is $\mathsf{P}2)$ in
Proposition~\ref{prop:MHPCauchy}] shows that 
$T_{\pa}^{\star} > T_{1}$. 

On the other hand the [omitted] proof of
Proposition~\ref{prop:MHPCauchy} implies that there exists $T_{2}>0$
and $M_{2}<+\infty$ such that for all $\pa=(\eps,\r,\pe)\in\PA$ with $\eps\ge \eps_{1}$,
one has $T_{\pa}^{\star}\ge T_{2}$ and $\Omega_{\pa}(T_{2})\le M_{2}$. 

We have proved~\eqref{task} with $T\defn \min\{T_{1}/2,T_{2}\}$ and $M\defn\max\{M_{1},M_{2}\}$, which completes the proof of Theorem~\ref{theo:uniform}.

\section{Decay of the local energy}\label{section:decay}

We will consider systems which include~\eqref{system:NS2} as a special case. 
The motivation is to consider a structure 
general enough to include, say, the combustion equations. To do so
we allow the limit constraint on the divergence of the velocity field
to read $\cn\vard=F(D^{m}\psi)$ where $\psi$ is the slow variable
(namely $\partial_{t}\psi=\Or(1)$), $m\ge 0$ is a given integer and 
$$
D^{m}\psi \defn \bigl\{\, \partial_{x}^{\alpha} \psi \,\arrowvert\,
\alpha\in\xN^{d},~\la\alpha\ra\le m\,\bigr\}.
$$ 

\smallbreak

More precisely, we consider system of the form:  
\begin{equation}\label{system:G}
\left\{
\begin{aligned}
&g_{1}(\psi)\partial_{t} \vari + \frac{1}{\eps} \cn\vard = 
\frac{1}{\eps}F_{1}(D^{m}\psi)+F_{2}(D^{m}\psi,D^{m}\vari,D^{m}\vard),\\
&g_{2}(\psi)\partial_{t} \vard + \frac{1}{\eps} \nabla\vari = F_{3}(D^{m}\psi,D^{m}\vari,D^{m}\vard),\\
&\partial_{t}\psi= F_{4}(D^{m}\psi,D^{m}\vari,D^{m}\vard),
\end{aligned}
\right.
\end{equation}
where $(\vari,\vard,\psi)$ is defined on
$[0,T]\times\xR^{d}$ with values
in~$\xR\times\xR^{d}\times\xR^{n}$, $n\ge 0$ and $m\ge 0$ are given integers.

\begin{assu} The functions $g_{i}$ ($i=1,2$) and $F_{j}$ ($j=1,2,3,4$)
  are smooth in their arguments, satisfying
$$
g_{1}(y),g_{2}(y) >0,\quad F_{1}(Y),F_{2}(Y,Z)\in \xR, \quad F_{3}(Y,Z)\in
\xR^{d},\quad F_{4}(Y,Z)\in\xR^{n},
$$
where $y$, $Y$ and $Z$ are the place holders of $\psi$, $D^{m}\psi$
and $(D^{m}\vari,D^{m}\vard)$. Moreover $F_{j}(0)=0$ for $j=1,2,3,4$. 
\end{assu}

Assuming that the slow variable decays sufficiently rapidly at spatial
infinity we prove that the penalized terms converge to $0$.
\begin{proposition}\label{prop:GCV}
Let $T>0$, $m\ge 0$ and $s>2m+d/2$. 
Assume that the functions $(\vari^{\eps},\vard^{\eps},\psi^{\eps})$
satisfy~\eqref{system:G} and 
\begin{equation}\label{assu:propdecayuniform}
\norm{(p^{\eps},\vard^{\eps},\psi^{\eps})}_{\Kr^{s}(T)}\defn
\sup_{t\in [0,T]}\norm{(\vari^{\eps},\vard^{\eps})(t)}_{\nh{s}}+\norm{\psi^{\eps}(t)}_{\nh{s+m-1}}<+\infty.
\end{equation}
Assume further that $\psi^{\eps}$ converges strongly in
$C^{0}([0,T];H^{\sigma}_{loc}(\xR^{d}))$, for some $\sigma >1+d/2$, to a
limit $\psi$ satisfying, for all $(t,x)\in [0,T]\times\xR^{d}$,
\begin{equation}\label{assu:propdecay}
        \la \psi(t,x)\ra\leqslant K \la x \ra^{-1-\gamma}, 
\quad \la \nabla \psi(t,x)\ra \leqslant  K \la x \ra^{-2-\gamma},
\end{equation}
for some positive constants~$K$ and~$\gamma$. Then, 
$\vari^{\eps}$ converges 
strongly to $0$ in $L^{2}(0,T;H^{s'}_{{loc}}(\xR^{d}))$ for all $s'<s$
and $\cn\vard^{\eps}-F_{1}(D^{m}\psi^{\eps})$ converges 
strongly to $0$ in $L^{2}(0,T;H^{s'}_{{loc}}(\xR^{d}))$ for all $s'<s-1$.
\end{proposition}
\begin{proof}
Since $(\vari^{\eps},\vard^{\eps},\psi^{\eps})$ is uniformly bounded in
$\Kr^{s}(T)$ (see~\eqref{assu:propdecayuniform}), using an interpolation argument, 
it is sufficient to prove
that $p^{\eps}$ and $\cn\vard^{\eps}-F_{1}(D^{m}\psi^{\eps})$ converges 
strongly to $0$ in $L^{2}(0,T;L^{2}_{{loc}}(\xR^{d}))$. 
To prove this result we use the following Theorem proved
in~\cite{MS1} [although not explicitly stated in this way]. To
clarify matters, we mention that: (1) all the functions considered in
the following statement are real-valued; (2) all the convergences
considered in this proof are strong.

\begin{theorem}[M\'etivier \& Schochet]\label{theo:MS}
Let~$T>0$ and let $u^{\eps}$ be a bounded sequence in~$C^{0}([0,T];H^{2}(\xR^{d}))$ such that 
$$
\eps^{2}\partial_{t}(a^{\eps}\partial_{t}u^{\eps}) - \cn
(b^{\eps}\nabla u^{\eps}) 
=c^{\eps},
$$
where the source term~$c^{\eps}$ converges to $0$
in~$L^{2}(0,T;L^{2}(\xR^{d}))$. Assume further that, for some $\sigma>1+d/2$, the coefficients~$(a^{\eps},b^{\eps})$ 
are uniformly bounded in $C^{0}([0,T];H^{\sigma}(\xR^{d}))$ and
converges in $C^{0}([0,T];H^{\sigma}_{loc}(\xR^{d}))$ 
to a limit $(a,b)$ satisfying the decay estimate
\begin{alignat*}{3}
&\la a(t,x)-\underline{a}\ra &&\leqslant K \la x \ra^{-1-\gamma}, 
\quad \la \nabla a(t,x)\ra
        &&\leqslant  K \la x \ra^{-2-\gamma},\\
&\la b(t,x)-\underline{b}\ra &&\leqslant K \la x \ra^{-1-\gamma}, 
\quad \la \nabla b(t,x)\ra
       && \leqslant  K \la x \ra^{-2-\gamma},
\end{alignat*}
for some given positive constants $\underline{a}$, $\underline{b}$,
$K$ and~$\gamma$. 

Then, the sequence~$u^{\eps}$ converges to~$0$ 
in~$L^{2}(0,T;L^{2}_{loc}(\xR^{d}))$.
\end{theorem}

We can directly apply Theorem~\ref{theo:MS} to prove the first half of
Proposition~\ref{prop:GCV}, that is the convergence of $p^{\eps}$ to $0$ in
$L^{2}(0,T;L^{2}_{loc}(\xR^{d}))$. Indeed, applying
$\eps^{2}\partial_{t}$ to the first equation in~\eqref{system:G}, we
compute
\begin{equation}\label{eq:pondes}
\eps^{2}\partial_{t}(a^{\eps}\partial_{t}p^{\eps})
-\cn(b^{\eps}\nabla p^{\eps})
=c^{\eps}\defn\eps f^{\eps},
\end{equation}
with $a^{\eps}\defn g_{1}(\psi^{\eps})$, $b^{\eps}\defn1/g_{2}(\psi^{\eps})$ and
\begin{equation}\label{defi:sourcepondes}
f^{\eps}\defn
\partial_{t}F_{1}(Y^{\eps})+\eps\partial_{t}F_{2}(Y^{\eps},Z^{\eps})-\cn\bigl(b^{\eps}F_{3}(Y^{\eps},Z^{\eps})\bigr),
\end{equation}
where $Y^{\eps}\defn D^{m}\psi^{\eps}$ and $Z^{\eps}\defn (D^{m}\vari^{\eps},D^{m}\vard^{\eps})$. 

The equations~\eqref{system:G} imply that one can 
rewrite $f^{\eps}$ as 
$$
f^{\eps} = \mathcal{F}(D^{m+1}p^{\eps},D^{m+1}\vard^{\eps},D^{m+1}\psi^{\eps}),
$$
where $\mathcal{F}$ is a $C^{\infty}$ function such that $\mathcal{F}(0)=0$. 
Using the usual nonlinear estimate in Sobolev spaces~\eqref{comp}, the
hypothesis~\eqref{assu:propdecayuniform} implies that $f^{\eps}$ is 
uniformly bounded in $C^{0}([0,T];H^{s-m-1}(\xR^{d}))$. Consequently, 
$c^{\eps}=\eps f^{\eps}$ converges to $0$ in $L^{2}(0,T;L^{2}(\xR^{d}))$ and
Theorem~\ref{theo:MS} applies.

\smallbreak
To prove the second half of
Proposition~\ref{prop:GCV}, we begin by proving that 
$\widetilde{p}^{\eps}\defn \edtt p^{\eps}$  
converges to $0$ in
$L^{2}(0,T;L^{2}_{loc}(\xR^{d}))$. To do so apply $\edtt$ on the
equation~\eqref{eq:pondes}, to obtain
\begin{equation*}%\label{eq:pondes2}
\eps^{2}\partial_{t}(a^{\eps}\partial_{t}\widetilde{p}^{\eps})
-\cn(b^{\eps}\nabla \widetilde{p}^{\eps})
=\widetilde{c}^{\eps}\defn \eps \widetilde{f}^{\eps},
\end{equation*}
with
\begin{equation*}
\widetilde{f}^{\eps}\defn
\eps\partial_{t}f^{\eps}-\eps\partial_{t}\bigl(\partial_{t} a^{\eps}
\edtt p^{\eps}\bigr) +\cn\bigl(\partial_{t} b^{\eps}\nabla p^{\eps}\bigr),
\end{equation*}
where $f^{\eps}$ is given by~\eqref{defi:sourcepondes}. Again one
can verify that $\widetilde{f}^{\eps}$ is a bounded sequence
in~$C^{0}([0,T];L^{2}(\xR^{d}))$, which proves the desired result. 

To complete the proof, observe that 
$$
\cn\vard^{\eps}-F_{1}(Y^{\eps})=-g_{1}(\psi^{\eps})\edtt\vari^{\eps}
+\eps F_{2}(Z^{\eps}).
$$
Hence, the fact that $\cn\vard^{\eps}-F_{1}(Y^{\eps})$  
converges to $0$ in
$L^{2}(0,T;L^{2}_{loc}(\xR^{d}))$ follows from the previous step
and the fact that
$g_{1}(\psi^{\eps})-g_{1}(0)$ and $F_{2}(Z^{\eps})$ are uniformly bounded in
$C^{0}([0,T];H^{s}(\xR^{d}))$. %\qed 
\end{proof}

\subsection{Proof of Theorem~$\ref{theo:decay}$}

To simplify the presentation we concentrate
on the hardest case when~$\pe$ is a fixed positive constant.

We first prove the convergences  for some
sub-sequence of $(p^{\eps},\vard^{\eps},\variii^{\eps})$.

The equation for $\variii$ implies that $\partial_{t}\variii$ is
bounded in $C^{0}([0,T];H^{s-1}(\xR^{d}))$. Therefore after 
extracting a sub-sequence, we can assume that, for all $s'<s$, 
\begin{equation}\label{limite:2}
\variii^{\varepsilon} \rightarrow \variii \quad  \text{  in  }  C^{0}([0,T];H^{s'+1} _{loc} (\xR^{d})),
\end{equation}
where the limit $\variii$ belongs to $C^{0}([0,T];H^{s'+1} _{loc}
(\xR^{d}))\cap L^{\infty}(0,T;H^{s+1}(\xR^{d}))$.

Since $(p^{\eps},\vard^{\eps})$ is uniformly
bounded in $C^{0}([0,T];H^{s}(\xR^{d}))$, after extracting further sub-sequence,
we can also assume that
\begin{equation}
(p^{\varepsilon} , v^{\varepsilon} ) \rightharpoonup (p , v) \quad \text{weakly } \star \text { in } 
L^{\infty} (0,T;H^{s} (\xR^{d})) \label{limites:1}.
\end{equation}

\smallbreak
Note that the system~\eqref{system:G} includes
system~\eqref{system:NS2} as a special case where
$\psi^{\eps}\defn(\variii^{\eps},\eps\vari^{\eps})$ and $m=2$. It follows from the
very definition of~$\psi^{\eps}$ and the norms $\norm{\cdot}_{H^{\sigma}_{\varrho}}$ that 
\begin{align*}
\norm{(\vari^{\eps},\vard^{\eps})}_{\nh{s}}+\norm{\psi^{\eps}}_{\nh{s+m-1}}
&\le\norm{(\vari^{\eps},\vard^{\eps})}_{\nh{s}}+\eps\norm{\vari^{\eps}}_{\nh{s+1}}+\norm{\variii^{\eps}}_{\nh{s+1}}\\
&\le
\norm{(\vari^{\eps},\vard^{\eps})}_{H^{s+1}_{\eps}}+\norm{\variii^{\eps}}_{\nh{s+1}}\\
&\les \norm{(\vari^{\eps},\vard^{\eps})}_{H^{s+1}_{\eps\cpe}}+\norm{\variii^{\eps}}_{H^{s+1}_{\cpe}},
\end{align*}
where the implicit constant depends only on the
fixed positive value of $\cpe\defn\sqrt{\r+\pe}$. It follows that 
$\norm{(p^{\eps},\vard^{\eps},\psi^{\eps})}_{\Kr^{s}(T)}\les
\norm{(\vari^{\eps},\vard^{\eps},\variii^{\eps})}_{\Hr_{(\eps,\r,\pe)}^{s}(T)}$.
As a consequence the first hypothesis of Proposition~\ref{prop:GCV} is
satisfied. Therefore, in order to apply Proposition~\ref{prop:GCV} to the
equations~\eqref{system:NS2} so as to obtain Theorem~\ref{theo:decay},
it only remains to show that $\psi=(\variii,0)$ satisfies 
the crucial assumption~\eqref{assu:propdecay}. 

\smallbreak
Multiplying the first equation in~\eqref{system:NS2} by $\eps$ and
passing to the weak limit (in the sense of distributions) shows that
$\cn\vard^{\eps}-\chi_{1}(\eps\vari^{\eps})\cn(\beta(\variii^{\eps})\nabla\variii^{\eps})$
converges to $0$ in that sense. Therefore, the limit $\variii$ satisfies
\begin{equation}\label{eq:parabolic}
g_{3}(\variii,0)(\partial_{t}\variii+\vard\cdot\nabla\variii) -\pe(\chi_{3}(0)-\chi_{1}(0))\cn(\beta(\variii)\nabla\variii)=0.
\end{equation}
The assumption $\chi_{3}>\chi_{1}$ [see $(\text{A}3)$ in
Assumption~\ref{assu:structural}] implies that~\eqref{eq:parabolic} is
parabolic. The desired spatial decay thus follows from estimates in weighted
Sobolev spaces. Indeed, introduce
$\variii_{\delta}\defn \L{x}^{\delta}\variii$ for some given $\delta>2$. Then
$$
\partial_{t}\variii_{\delta}-\cn\bigl(k(\variii)\nabla\variii_{\delta}\bigr)=\sum
V_{\alpha}(t,x)\partial_{x}^{\alpha}\variii_{\delta}  \qquad (\alpha\in\xN^{d},~ \la\alpha\ra\le 1),
$$
for some smooth positive function $k$ and some coefficients
$V_{\alpha}$ such that $\norm{V_{\alpha}}_{\nh{s}}\le
C(\norm{(\vard,\variii)}_{\nh{s}})$. Commuting $\Fi{}{s}$ with the
equation and applying the usual $L^{2}$ estimate, we obtain
\begin{equation}\label{esti:wpara}
\norm{\variii_{\delta}}_{L^{\infty}(0,T;H^{s})}\le C
e^{TC}\norm{\variii_{\delta}(0)}_{\nh{s}},
\end{equation}
for some constant $C$ depending only on
$\norm{(\vard,\variii)}_{L^{\infty}(0,T;H^{s})}$
(see~\eqref{esti:hypparusual} with $\eta=1$). The convergence~\eqref{limite:2} implies that $\variii(0)$ is the limit $\variii_{0}$ of the initial
data~$\variii^{\eps}(0)$. Using the hypothesis $\L{x}^{\delta}\variii_{0}\in H^{s}(\xR^{d})$, it follows
from~\eqref{esti:wpara} that
$$
\norm{\variii_{\delta}}_{L^{\infty}(0,T;H^{s})} <+\infty.
$$
Since $s>1+d/2$, the Sobolev Theorem implies that
$H^{s}(\xR^{d})\hookrightarrow W^{1,\infty}(\xR^{d})$. Hence, the
assumption~\eqref{assu:propdecay} is satisfied.

\smallbreak
A moment's thought shows that the convergences 
holds for the full sequence. 
The proof of Theorem~\ref{theo:decay} is complete.

\begin{appendix}
\section{Estimates for hyperbolic and parabolic systems}
Here we briefly recall some standard estimates for hyperbolic or
parabolic systems we used throughout the paper. We state estimates for the solutions
$u=u(t,x)\in\xR^{n}$ of systems having the form:
\begin{equation}\label{Eq:hypparusual}
\partial_{t}u+\sum_{\la\alpha\ra=1}A_{\alpha}\partial_{x}^{\alpha}u-\eta
\sum_{\la\alpha\ra=2}A_{\alpha}\partial_{x}^{\alpha} u =f_{1}+\sqrt{\eta}f_{2},
\end{equation}
where each $A_{\alpha}=A_{\alpha}(t,x)$ is a $n\times n$ matrix
valued-function smooth in its arguments and furthermore symmetric:
$A_{\beta}=A_{\beta}^{t}$; also, the source terms
$f_{i}=f_{i}(t,x)$ are smooth functions with values in $\xR^{n}$. 
In particular, when $\eta=0$, the system \eqref{Eq:hypparusual} is
symmetric hyperbolic. We assume that for $\eta>0$ the system is
parabolic: for all triple $(t,x,\xi)$ with $\xi\neq 0$, 
$$
\sum_{\la\alpha\ra=2}A_{\alpha}(t,x)\xi^{\alpha}>0.
$$

\begin{lemma}Let $0\le\sigma\le s \in (1+d/2,+\infty)$. There exists a function $C(\cdot)$ such
  that for all $T\in [0,1]$, all $u\in C^{1}([0,T];H^{\infty}(\xD))$ 
satisfying~\eqref{Eq:hypparusual} and all $\eta\in [0,1]$,
\begin{equation}\label{esti:hypparusual}
\begin{split}
&\norm{u}_{L^{\infty}(0,T;H^{\sigma})}+\sqrt{\eta}\norm{u}_{L^{2}(0,T;H^{\sigma+1})}\\
&\qquad\le K e^{C(R)T}
\left(\norm{u(0)}_{\nh{\sigma}}+\norm{f_{1}}_{L^{1}(0,T:H^{\sigma})}+\norm{f_{2}}_{L^{2}(0,T;H^{\sigma-1})}\right),
\end{split}
\end{equation}
where $K=K(s,d)$ and 
$R\defn\sum_{\la\alpha\ra\le 2}\norm{\tilde{A}_{\alpha}}_{L^{\infty}(0,T;H^{s})}$ 
where $\tilde{A}=A-A(0)$.
\end{lemma}
\begin{proof}
Commuting the equation~\eqref{Eq:hypparusual} with $\Fi{}{\sigma}$, using
the commutator estimate~\eqref{commutator:estimateusual} 
and the Gronwall's lemma, the proof of~\eqref{esti:hypparusual}
can be reduced to the special case $\sigma=0$. 
The latter case relies upon usual integrations by
parts and duality arguments.
\end{proof}

\end{appendix}

%%%%%%%%%%%%%%%%%%%%%%%%%%%%%%%%%%%%%%%%%%%%%%%%%%%%%%%%%%%%%%%%%%%%%%%%%%%%%%%%%%%%%%%%%%%%%%%%%%%%%%%%%%%%%%%%%%%%%%%%%%%%%%%%%%%%
%%%%%%%%%%%%%%%%%%%%%%%%%%%%%%%%%%%%%%%%%%%%%%%%%%%%%%%%%%%%%%%%%%%%%%%%%%%%%%%%%%%%%%%%%%%%%%%%%%%%%%%%%%%%%%%%%%%%%%%%%%%%%%%%%%%%
%%%%%%%%%%%%%%%%%%%%%%%%%%%%%%%%%%%%%%%%%%%%%%%%%%%%%%%%%%%%%%%%%%%%%%%%%%%%%%%%%%%%%%%%%%%%%%%%%%%%%%%%%%%%%%%%%%%%%%%%%%%%%%%%%%%%
%%%%%%%%%%%%%%%%%%%%%%%%%%%%%%%%%%%%%%%%%%%%%%%%%%%%%%%%%%%%%%%%%%%%%%%%%%%%%%%%%%%%%%%%%%%%%%%%%%%%%%%%%%%%%%%%%%%%%%%%%%%%%%%%%%%%
%%%%%%%%%%%%%%%%%%%%%%%%%%%%%%%%%%%%%%%%%%%%%%%%%%%%%%%%%%%%%%%%%%%%%%%%%%%%%%%%%%%%%%%%%%%%%%%%%%%%%%%%%%%%%%%%%%%%%%%%%%%%%%%%%%%%

{\small \noindent\textbf{Acknowledgments.} 
I would like to express my sincere
gratitude to G.~M\'etivier for introducing me to the
analysis of the low Mach number limit as well as for helpful 
discussions. I warmly thanks D.~Bresch and D.~Lannes for many advices.
I am most appreciative for the hospitality shown to me by the Rennes university during the beginning of this work.
}

\bibliographystyle{smfplain}
\bibliography{bib_thesis2}

\providecommand{\bysame}{\leavevmode ---\ }
\providecommand{\og}{``}
\providecommand{\fg}{''}
\providecommand{\smfandname}{\&}
\providecommand{\smfedsname}{\'eds.}
\providecommand{\smfedname}{\'ed.}
\providecommand{\smfmastersthesisname}{M\'emoire}
\providecommand{\smfphdthesisname}{Th\`ese}
\begin{thebibliography}{10}

\bibitem{TA}
{\scshape T.~Alazard} -- {\og Incompressible limit of the nonisentropic {E}uler
  equations with solid wall boundary conditions\fg}, \emph{Adv. in Differential
  Equations} \textbf{10} (2005), p.~19--44.

\bibitem{Bony}
{\scshape J.-M. Bony} -- {\og Calcul symbolique et propagation des
  singularit\'es pour les \'equations aux d\'eriv\'ees partielles non
  lin\'eaires\fg}, \emph{Ann. Sci. \'Ecole Norm. Sup. (4)} \textbf{14} (1981),
  no.~2, p.~209--246.

\bibitem{BDGL}
{\scshape D.~Bresch, B.~Desjardins, E.~Grenier {\normalfont \smfandname} C.-K.
  Lin} -- {\og Low {M}ach number limit of viscous polytropic flows: formal
  asymptotics in the periodic case\fg}, \emph{Stud. Appl. Math.} \textbf{109}
  (2002), no.~2, p.~125--149.

\bibitem{BRDL}
{\scshape D.~Bresch, B.~Desjardins {\normalfont \smfandname} C.-K. Lin} -- {\og
  On some compressible fluid models: {K}orteweg, lubrication, and shallow water
  systems\fg}, \emph{Comm. Partial Differential Equations} \textbf{28} (2003),
  no.~3-4, p.~843--868.

\bibitem{CGM}
{\scshape C.~Cheverry, O.~Gu{\`e}s {\normalfont \smfandname} G.~M{\'e}tivier}
  -- {\og Oscillations fortes sur un champ lin\'eairement d\'eg\'en\'er\'e\fg},
  \emph{Ann. Sci. \'Ecole Norm. Sup. (4)} \textbf{36} (2003), no.~5,
  p.~691--745.

\bibitem{Dan}
{\scshape R.~Danchin} -- {\og Global existence in critical spaces for
  compressible {N}avier-{S}tokes equations\fg}, \emph{Invent. Math.}
  \textbf{141} (2000), no.~3, p.~579--614.

\bibitem{DanE1}
\bysame , {\og Global existence in critical spaces for flows of compressible
  viscous and heat-conductive gases\fg}, \emph{Arch. Ration. Mech. Anal.}
  \textbf{160} (2001), no.~1, p.~1--39.

\bibitem{Dan1}
\bysame , {\og Zero {M}ach number limit for compressible flows with periodic
  boundary conditions\fg}, \emph{Amer. J. Math.} \textbf{124} (2002), no.~6,
  p.~1153--1219.

\bibitem{Dan2}
\bysame , {\og Zero {M}ach number limit in critical spaces for compressible
  {N}avier-{S}tokes equations\fg}, \emph{Ann. Sci. \'Ecole Norm. Sup. (4)}
  \textbf{35} (2002), no.~1, p.~27--75.

\bibitem{DG}
{\scshape B.~Desjardins {\normalfont \smfandname} E.~Grenier} -- {\og Low
  {M}ach number limit of viscous compressible flows in the whole space\fg},
  \emph{R. Soc. Lond. Proc. Ser. A Math. Phys. Eng. Sci.} \textbf{455} (1999),
  no.~1986, p.~2271--2279.

\bibitem{DGLM}
{\scshape B.~Desjardins, E.~Grenier, P.-L. Lions {\normalfont \smfandname}
  N.~Masmoudi} -- {\og Incompressible limit for solutions of the isentropic
  {N}avier-{S}tokes equations with {D}irichlet boundary conditions\fg},
  \emph{J. Math. Pures Appl. (9)} \textbf{78} (1999), no.~5, p.~461--471.

\bibitem{DL}
{\scshape B.~Desjardins {\normalfont \smfandname} C.-K. Lin} -- {\og A survey
  of the compressible {N}avier-{S}tokes equations\fg}, \emph{Taiwanese J.
  Math.} \textbf{3} (1999), no.~2, p.~123--137.

\bibitem{Gallagher}
{\scshape I.~Gallagher} -- {\og R\'esultats r\'ecents sur la limite
  incompressible\fg}, \emph{S\'eminaire Bourbaki 2003--2004} \textbf{926}
  (2003), no.~926.

\bibitem{Grenier}
{\scshape E.~Grenier} -- {\og Oscillatory perturbations of the
  {N}avier-{S}tokes equations\fg}, \emph{J. Math. Pures Appl. (9)} \textbf{76}
  (1997), no.~6, p.~477--498.

\bibitem{HeVW}
{\scshape D.~R. van~der Heul, C.~Vuik {\normalfont \smfandname} P.~Wesseling}
  -- {\og A conservative pressure-correction method for the {E}uler and ideal
  {MHD} equations at all speeds\fg}, \emph{Internat. J. Numer. Methods Fluids}
  \textbf{40} (2002), no.~3-4, p.~521--529, ICFD Conference on Numerical
  Methods for Fluid Dynamics, Part II (Oxford, 2001).

\bibitem{Hoff}
{\scshape D.~Hoff} -- {\og The zero-{M}ach limit of compressible flows\fg},
  \emph{Comm. Math. Phys.} \textbf{192} (1998), no.~3, p.~543--554.

\bibitem{HorL}
{\scshape L.~H{\"o}rmander} -- \emph{Lectures on nonlinear hyperbolic
  differential equations}, Math\'ematiques \& Applications, vol.~26,
  Springer-Verlag, Berlin, 1997.

\bibitem{Iso}
{\scshape H.~Isozaki} -- {\og Singular limits for the compressible {E}uler
  equation in an exterior domain\fg}, \emph{J. Reine Angew. Math.} \textbf{381}
  (1987), p.~1--36.

\bibitem{Iso2}
\bysame , {\og Singular limits for the compressible {E}uler equation in an
  exterior domain. {II}. {B}odies in a uniform flow\fg}, \emph{Osaka J. Math.}
  \textbf{26} (1989), no.~2, p.~399--410.

\bibitem{KS}
{\scshape S.~Kawashima {\normalfont \smfandname} Y.~Shizuta} -- {\og On the
  normal form of the symmetric hyperbolic-parabolic systems associated with the
  conservation laws\fg}, \emph{Tohoku Math. J. (2)} \textbf{40} (1988), no.~3,
  p.~449--464.

\bibitem{KY}
{\scshape S.~Kawashima {\normalfont \smfandname} W.-A. Yong} -- {\og
  Dissipative structure and entropy for hyperbolic systems of balance laws\fg},
  \emph{Arch. Ration. Mech. Anal.} (To appear), p.~1--20.

\bibitem{KM1}
{\scshape S.~Klainerman {\normalfont \smfandname} A.~Majda} -- {\og Singular
  limits of quasilinear hyperbolic systems with large parameters and the
  incompressible limit of compressible fluids\fg}, \emph{Comm. Pure Appl.
  Math.} \textbf{34} (1981), no.~4, p.~481--524.

\bibitem{KM2}
\bysame , {\og Compressible and incompressible fluids\fg}, \emph{Comm. Pure
  Appl. Math.} \textbf{35} (1982), no.~5, p.~629--651.

\bibitem{Klein}
{\scshape R.~Klein} -- {\og Semi-implicit extension of a {G}odunov-type scheme
  based on low {M}ach number asymptotics. {I}. {O}ne-dimensional flow\fg},
  \emph{J. Comput. Phys.} \textbf{121} (1995), no.~2, p.~213--237.

\bibitem{David}
{\scshape D.~Lannes} -- {\og Sharp estimates for pseudo-differential operators
  with symbols of limited smoothness and commutators\fg}, \emph{Preprint
  Bordeaux1} (2004), p.~1--31.

\bibitem{Lions}
{\scshape P.-L. Lions} -- \emph{Mathematical topics in fluid mechanics. {V}ol.
  1}, Oxford Lecture Series in Mathematics and its Applications, vol.~3.

\bibitem{LionsMas}
{\scshape P.-L. Lions {\normalfont \smfandname} N.~Masmoudi} -- {\og
  Incompressible limit for a viscous compressible fluid\fg}, \emph{J. Math.
  Pures Appl. (9)} \textbf{77} (1998), no.~6, p.~585--627.

\bibitem{Maj}
{\scshape A.~Majda} -- \emph{Compressible fluid flow and systems of
  conservation laws in several space variables}, Applied Mathematical Sciences,
  vol.~53, Springer-Verlag, New York, 1984.

\bibitem{MN}
{\scshape A.~Matsumura {\normalfont \smfandname} T.~Nishida} -- {\og The
  initial value problem for the equations of motion of viscous and
  heat-conductive gases\fg}, \emph{J. Math. Kyoto Univ.} \textbf{20} (1980),
  no.~1, p.~67--104.

\bibitem{MS2}
{\scshape G.~M{\'e}tivier {\normalfont \smfandname} S.~Schochet} -- {\og Limite
  incompressible des \'equations d'{E}uler non isentropiques\fg}, in
  \emph{S\'eminaire: \'Equations aux D\'eriv\'ees Partielles, 2000--2001}.

\bibitem{MS1}
\bysame , {\og The incompressible limit of the non-isentropic {E}uler
  equations\fg}, \emph{Arch. Ration. Mech. Anal.} \textbf{158} (2001), no.~1,
  p.~61--90.

\bibitem{MS3}
\bysame , {\og Averaging theorems for conservative systems and the weakly
  compressible {E}uler equations\fg}, \emph{J. Differential Equations}
  \textbf{187} (2003), no.~1, p.~106--183.

\bibitem{Meyer}
{\scshape Y.~Meyer} -- {\og Remarques sur un th\'eor\`eme de {J}.-{M}.
  {B}ony\fg}, in \emph{Proceedings of the Seminar on Harmonic Analysis (Pisa,
  1980)}, no. suppl. 1, 1981, p.~1--20.

\bibitem{Munz}
{\scshape C.-D. Munz} -- {\og Computational fluid dynamics and aeroacoustics
  for low {M}ach number flow\fg}, in \emph{Hyperbolic partial differential
  equations (Hamburg, 2001)}, Vieweg, Braunschweig, 2002, p.~269--320.

\bibitem{SBGK}
{\scshape T.~Schneider, N.~Botta, K.~J. Geratz {\normalfont \smfandname}
  R.~Klein} -- {\og Extension of finite volume compressible flow solvers to
  multi-dimensional, variable density zero {M}ach number flows\fg}, \emph{J.
  Comput. Phys.} \textbf{155} (1999), no.~2, p.~248--286.

\bibitem{SchoE}
{\scshape S.~Schochet} -- {\og The compressible {E}uler equations in a bounded
  domain: existence of solutions and the incompressible limit\fg}, \emph{Comm.
  Math. Phys.} \textbf{104} (1986), no.~1, p.~49--75.

\bibitem{SchoFast}
\bysame , {\og Fast singular limits of hyperbolic {PDE}s\fg}, \emph{J.
  Differential Equations} \textbf{114} (1994), no.~2, p.~476--512.

\bibitem{Secchi}
{\scshape P.~Secchi} -- {\og On slightly compressible ideal flow in the
  half-plane\fg}, \emph{Arch. Ration. Mech. Anal.} \textbf{161} (2002), no.~3,
  p.~231--255.

\bibitem{Se}
{\scshape P.~Secchi} -- {\og On the singular incompressible limit of inviscid
  compressible fluids\fg}, \emph{J. Math. Fluid Mech.} \textbf{2} (2000),
  no.~2, p.~107--125.

\bibitem{Taylor}
{\scshape M.~E. Taylor} -- \emph{Pseudodifferential operators and nonlinear
  {PDE}}, Progress in Mathematics, vol. 100, Birkh\"auser Boston Inc., Boston,
  MA, 1991.

\bibitem{Ukai}
{\scshape S.~Ukai} -- {\og The incompressible limit and the initial layer of
  the compressible {E}uler equation\fg}, \emph{J. Math. Kyoto Univ.}
  \textbf{26} (1986), no.~2, p.~323--331.

\bibitem{VT}
{\scshape V.~N. Vatsa {\normalfont \smfandname} E.~Turkel} -- {\og Choice of
  variables and preconditioning for time dependent problems\fg}, in \emph{AIAA,
  16th Computational Fluid Dynamics Conference (Orlando)}, 2003.

\bibitem{ZeytounianE}
{\scshape R.~K. Zeytounian} -- \emph{Theory and applications of nonviscous
  fluid flows}, Springer-Verlag, Berlin, 2002.

\bibitem{ZeytounianNS}
\bysame , \emph{Theory and applications of viscous fluid flows},
  Springer-Verlag, Berlin, 2004.

\end{thebibliography}
\end{document}